\documentclass{amsart}
\usepackage{relsize}
\usepackage{url}
\usepackage{mathdots}
\usepackage{amsmath}
\usepackage{fancyhdr}
\usepackage{amsfonts}
\usepackage{amssymb}
\usepackage{amscd}
\usepackage{amsthm}
\usepackage{verbatim}
\usepackage{mathrsfs}

\newtheorem{thm}{Theorem}[section]

\newtheorem{cor}[thm]{Corollary}
\newtheorem{lemma}[thm]{Lemma}
\newtheorem{prop}[thm]{Proposition}

\newcommand{\be}{\begin{equation}}
\newcommand{\ee}{\end{equation}}
\newcommand{\bea}{\begin{eqnarray}}
\newcommand{\eea}{\end{eqnarray}}

\newcommand{\boldell}{\mbox{\boldmath$\ell$\unboldmath}}
\newcommand{\boldlambda}{\mbox{\boldmath$\lambda$\unboldmath}}

\newcommand{\boldop}{\mbox{\boldmath$($\unboldmath}}
\newcommand{\boldcp}{\mbox{\boldmath$)$\unboldmath}}

\begin{document}
\title[Residual-Data Eisentstein Series, I]{
Analytic Properties of Residual Eisenstein Series, I}
\author{Eliot Brenner}
\address{University of Minnesota, School of Mathematics\newline
206 Church Street, SE\newline
Minneapolis USA, 55455}
\email{brenn240@umn.edu}
\begin{abstract}
We partially generalize the results of \cite{kudlarallisfest} on the
poles of degenerate, Siegel-parabolic Eisenstein series
to residual-data Eisenstein series.  In particular, for $a,b$ integers
greater than $1$,
we show that poles of the Eisenstein series induced
from the Speh representation $\Delta(\tau,b)$ on the Levi $\mathrm{GL}_{ab}$
of $\mathrm{Sp}_{2ab}$ are located in the `segment'
of half integers $X_{b}$ between a `right endpoint' and its negative, inclusive
of endpoints.
The right endpoint is $\pm b/2$, or $(b-1)/2$, depending on the analytic properties
of the automorphic $L$-functions attached to $\tau$.
We study the automorphic forms $\Phi_{i}^{(b)}$ obtained as residues
at the points $s_i^{(b)}$ (defined precisely in the paper)
by calculating their
cuspidal exponents in certain cases.
In the case of the `endpoint' $s_0^{(b)}$ and `first interior point' 
$s_1^{(b)}$ in the segment
of singularity points, we are able to determine a set containing \textit{all possible}
cuspidal exponents of $\Phi_0^{(b)}$ and $\Phi_1^{(b)}$ precisely for all $a$ and $b$.
In these cases, we use the result of the calculation to deduce that
the residual automorphic forms lie in
$L^2(G(k)\backslash G(\mathbf{A}))$.  In a more precise
sense, our result establishes a relationship between, on the one hand,
the actually occurring cuspidal exponents of $\Phi_i^{(b)}$,
residues at interior points which lie to the right of the origin,
and, on the other hand, the `analytic properties' of the original residual-data Eisenstein
series at the origin.   MSC Numbers: 11F70 (22E55)
\end{abstract}
\maketitle

\section{Introduction}
\label{sec:introduction}

In this paper, we initiate a systematic study of the singularities of the
Eisenstein series of split classical groups which have non-cuspdial
(or not necessarily cuspidal), but discrete-spectrum automorphic forms as data.  
Since the literature already contains
extensive work (see, \textit{e.g.}, \cite{krannals}, \cite{vtan}) on the degenerate-data
case (induced from a character on the Borel), we will concentrate our efforts
on the complementary case when the \textit{cuspdial support} of the data
is a non-minimal parabolic.  We obtain one complete result, Theorem \ref{thm:definitive},
giving an inductive formula for the `controlling' constant
term of the Eisenstein series, and describing  bounded, finite set of points at which all \textit{possible} points
of the (normalized) Eisenstein series must be located.  As one of several
possible applications of this formula, we use it to obtain some
partial--but for endpoints and `first-interior points' definitive--results
on the non-vanishing and square integrability of residues of the Eisenstein
series, summed up in Theorem \ref{thm:main}.

\subsection{General Context}
\label{subsec:generalcontext}
For this paper, at least, we will restrict our investigations to
one example of a wider context.
Nevertheless, it is useful first to make a brief
survey of the wider theory of Eisenstein series in the context of spectral
expansions. Such a discussion reveals
that the example considered in this paper is much more representative of the general
case than might appear at first sight.  

Let $G=G_n$ be a
semisimple algebraic group defined
over a global field $k$ (in this paper always a number field), of rank $n$.  For
simplicity, the reader
may consider only the cases of $G=\mathrm{GL}_n$
the general linear group, and the split classical groups, namely
$\mathrm{Sp}_{2n}$, $\mathrm{SO}_{2n}$, $\mathrm{SO}_{2n+1}$, and 
$U_{n,n}$.  Our primary object in this discussion is the representations of the adelic points
$G(\mathbf{A})$ of the group, so we often refer to $G(\mathbf{A})$
simply as $G$, but we also sometimes refer to the local group
$G_v=G(k_v)$, the points of $G$ over the completion of a single valuation $v\in\Omega(k)$.
The overarching purpose of the Langlands program may be described
as ``relating different ways of parametrizing" the irreducible admissible
representations of $G_v$ (local version), respectively 
the irreducible automorphic representations of $G$ (global version).
Of the different ways of parametrizing these representations,
perhaps the central one has turned out to be \textit{parabolic induction},
which involves Eisenstein series as an additional element in the global, or automorphic, case. 
The essential data for an (parabolically) induced admissible representation
of $G_v$, is a standard parabolic $P=MU$ of $G$ and an irreducible admissible representation
 $\pi_v$ of the Levi component
$M_v$; while for an Eisenstein series on $G$, it is the parabolic $P$ and a unitary
automorphic representation $\pi$ of $M$.  In addition to this basic data, there is
also a character of the center of $M$ (more precisely, of its `$k$-split component')
which serves as a ``parameter".  The `interesting' values of the parameter,
in the local case, are the points of reducibility of the induced representation
and in the global case, the poles of the Eisenstein series.
These facts are the basis of a somewhat imprecise but extremely useful analogy
between the theory of (local) induced admissible representations and (global)
Eisenstein series, that we will refer to repeatedly in this introduction.

As a first manifestation of this analogy, we recall that in the local
theory a general and very simple principle that goes
by the names of `transitivity of induction' or `induction in stages'
says that the family of representations of $G$ which are induced from
representations $\pi$ of $M$ which are \textit{themselves parabolically}
from a representation $\pi'$ of a standard Levi $M'$ in $M$ (hence standard
Levi of $G$ as well), are best viewed as 
the family of induced representations from data the $(M',\pi')$.  
Analogously, the much more involved general theory of Langlands \cite{langlands76},
as exposited in Chapter V of \cite{mwbook}, implies that for $\pi$
belonging to the orthogonal complement of the
discrete part $L^2(M)_0$ of the spectrum of $M$,
the data is itself expressible as a certain integral of Eisenstein
series on $M$.  Then one is best treating Eisenstein series on $G$
as a certain integral of Eisenstein series induced from the smaller parabolic.  
(See Proposition \ref{prop:residualdataeisseries} and its proof
for fully realized illustration of this kind of thinking ``in reverse".)

 As is well-known from the study of root systems and Dynkin diagrams,
a standard proper parabolic of $G$ has Levi component $M$
which decomposes as product of general linear groups and (possibly) a group
`of the same type', with the rank of all these factors adding up to $n$.
For example, when as in this paper, we take $G_n=\mathrm{Sp}_{2n}$, 
the most general $M$ has the form
\[
M\cong \mathrm{GL}_{n_1}\times \mathrm{GL}_{n_r}\times \mathrm{Sp}_{m}\; \text{with}
\sum_{i=1}^{r} n_i+ m=n,
\]
allowing the possibility that $m=0$ and there is no $G_m$ in the product.
Then the discrete data factors as a tensor product
$\pi\cong \pi_1\cdots\otimes\cdots \pi_r\otimes \sigma$
of discrete data on the general linear factors and on classical-group factor $G_m$.
Further $P$ is maximal if and only if $r=1$.  Also,
the `parameter' of the Eisenstein series is in $\mathbf{C}^r$,
so when we seek ``residues" of these Eisenstein series, we really seek
$r$-fold singularities in the sense of complex analysis of several variables.
When $r$ is not equal to $1$ and 
the number of general linear factors is greater than $1$, we associate
with $P$, $P_{\rm max}$, the standard maximal parabolic
containing $P$, with $M_{\rm max}\cong \mathrm{GL}_{n-m}\times \mathrm{G}_m$.
We can invoke a second simplifying argument by a parallel with local representation theory.
Namely, it is best to
think of the Eisenstein series induced from cuspidal data on $M$
as instead being Eisenstein series with data $P_{\rm max}$ and
Eisenstein series on $M_{\rm max}$ induced from $P\cap M_{\rm max}$. 
Then the ``residues" of the $P$-induced Eisenstein series
are for our purposes better conceived of residues--now, in the ordinary
sense of one-variable complex analysis--of residual-data Eisenstein series.
This is precisely the outcome of Proposition \ref{prop:residualdataeisseries},
which is not a particular deep result, and can be expected to generalize
to all classical groups.
The utility of this approach only really becomes apparent
when we consider the much deeper result of \cite{moeglinwaldspurgerens}
describing the residual spectrum of the general linear group completely (recalled
in \ref{thm:mwmain} below):
the upshot is that we need only consider the case when all $\mathrm{n_i}$ are equal
to a single value $b$ and all $\pi_i$ are self-dual isomorphic, say equal to the single
representation $\tau$.  

Therefore, if we are just interested in residues of Eisenstein
of `maximal' order we lose nothing by considering the case when $P$
is already maximal (so $r=1$) and the data is, not in general cuspidal,
but rather residual-spectrum, and of the form $\pi\otimes\sigma$, with
\[
\pi\cong \Delta(\tau,b)=J_{P(\mathbf{A})}^{\mathrm{GL}_{ab}(\mathbf{A})}(\tau^{\otimes b},\Lambda_b),
\]
and $\sigma$ a discrete-data representation
on $G_{n-m}$.  This notation, explained in full detail \S\ref{subsec:inductioninstages},
refers to the global Speh representation (special case
of the Langlands quotient) associated to $\tau$ on $\mathrm{GL}_{ab}$.
We do make two major simplifying assumptions for the purposes
of this paper: first that $L(\tau,\frac{1}{2})\neq 0$, and second
\textit{that $m=0$ so that the parabolic $P$ is the Siegel
parabolic}.  As a consequence of the second assumption the cuspidal representation of $G_m$
`disappears', or, more precisely, is replaced by the trivial representation of the trivial group.
This assumption, while admittedly a substantive one, is a sensible one to make
at the outset, because it allows us to avoid the complications
arising from the non-genericity of certain cuspidal representations
of the classical groups.

\subsection{Result contained in this paper.}
The main unconditional result of this paper, Theorem \ref{thm:definitive}(b)
that the residues of the 
normalized Eisenstein series $E^*(\Delta(\tau,b),s)$ (defined rigorously in \S\ref{subsec:residues}, normalized in 
\S\ref{sec:principalnormconst}) are located precisely at the points $s\in X_b$,
where $X_b$ is the segment of points on the real line stretching from an ``endpoint"
$s_0^{(b)}$ to its opposite $-s_0^{(b)}$, and at integer distance from each endpoint.  The value of $s_0^{(b)}$ is either
$b/2$ or $(b-1)/2$, depending on the `type' of $\tau$ in the sense described
 in \S\ref{subsec:generalnotation} below.  Further, the possible cuspidal exponents 
 of the residues $\Phi_0^{(b)}$ at the endpoints and $\Phi_1^{(b)}$ at the first interior
 points are calculated and, as a corollary, these residues are shown to be square integrable.
The square integrability of the family of automorphic forms $\Phi_0^{(b)}$ is already expected and essentially known
to experts, but that of the entire family of $\Phi_1^{(b)}$ is apparently new (but known from
the special calculation contained in \cite{pogge} in the case $a=b=2$).

However, as will become evident on a closer reading, the more fundamental main
main result is really the inductive
formula for the `principal constant term', found in its  normalized form in
Theorem \ref{thm:definitive}(a) again in representation-theoretic form in
 \eqref{eqn:inductiveformulareps}.  We expect that various pieces of useful
 information can be extracted from this formula.  One direction
 explored in this paper is that it can be used to extract
 some information concerning square-integrability of the residues,
 because it sets up a dependence of the cuspidal exponents of the $\Phi_i^{(b)}$ for $0<i<b/2$,
(\textit{i.e.}, for the interior points of the segment of poles on the positive real axis)
on the \textit{order of zero} of a non-normalized Eisenstein series $E^{a\cdot b'}$ series at the origin.
Here $E^{a\cdot b'}$ is an Eisenstein series of the same type but of smaller
rank $b'=i+1$ or $b'=i+2$, depending on the type of $\tau$, and with data
the image of $\Delta(\tau,b)$ under a certain normalized intertwining operator.
The impediment to our extending our result, in the form of a recursive
description of the cuspidal exponents, to all `interior residues' $\Phi_i^{(b)}$
is lack of complete knowledge of the \textit{order of zeros} of these Eisenstein series at the origin.
Ad hoc calculations are related to the case of $b'=2$ and $b'=3$ carried
out in \S\ref{subsec:partialresults}, and are sufficient to narrow
down the possibilities for cuspidal exponents in the case $i=1$
to allow us to prove our ``main result" on square-integrability cited
in the previous paragraph.

\subsection{Relation to other results in the literature.}
Such a plethora of papers in the literature relate in some way
to the example considered here that we cannot
attempt to do them justice, and will leave a full discussion
for the sequel paper.  We will only mention the previous
papers most directly related and some survey articles
containing a more comprehensive list of references.

The two papers most directly related to our results are those of Kudla-Rallis
in \cite{kudlarallisfest}, and Pogge in \cite{pogge}.
The example considered in this paper is \textit{almost}
a generalization of the results of the first section of Kudla-Rallis,
in the sense that the degenerate case may be identified
as the case $a=1$ in our setup.  In fact Proposition 1.2.1, in
what they call Case 1 (symplectic group), directly generalizes to our setup in the guise
of the inductive formula \eqref{eqn:inductiveformulareps}
for the constant term.  Nevertheless, we prefer to think
of their degenerate case as a related case alongside
our ``generalized principle series" case, because of
the fact that the the set of residue points $X_n$,
which would be called $X_b$ in our notation, has a different endpoint.
This results, ultimately, from the fact that the standard $L$-function
of the cusp form in their formulas is replaced by the zeta function.
The case that Pogge considers in \cite{pogge} lines
up exactly with what we would call the case of $a=b=2$.
Referring to Figure 1 in his paper, the reader
may verify that the ``real axis" in the parameter $s$
in our paper corresponds exactly to the affine hyperplane $S_1$
in his paper.  However, because he only considers the case $b=2$
there are nothing but ``endpoints" in his case.  Although
$\beta$ and $\gamma_4$ are plotted on the same affine hyperplane $S_1$,
they are in fact both ``endpoints" $s_0^{(2)}$ in our notation,
with $\gamma$ corresponding to the $\tau$-symplectic
and $\beta_4$ corresponding to the $\tau$-orthogonal case.

In addition to the very extensive literature
partially computing the residual spectrum
by construction of residues for many split classical
groups by Kim, Moeglin, Shahidi and others Moeglin in \cite{moeglinmanuscripta}
has a very precise description of the residual spectrum
of split classical groups, including $\mathrm{Sp}_n$ as a special
case.  This description is not unconditionally
proven, but based on the Arthur's Conjectures.
We mention this in order to point out that this work
does not remove the need for analysis such as undertaken
in our work.  The reason is that only certain of the
square-integrable residues of the Eisenstein series
are needed to construct the entire residual spectrum,
so work taking the Arthur conjectures as their starting point
do not analyze the cuspidal exponents and square-integrability of 
\textit{all} possible residues of the Eisenstein series.
The properties of all the residues do however become
important in such applications as analysis of $L$-functions via
integral representations involving the Eisenstein series, for example
the Rankin-Selberg method.

\subsection{Further Developments and Generalizations.}
The most straightforward generalization would be to consider
the case of other $k$-split or quasi-split classical groups.  The previous
experience of \cite{vtan} and \cite{krannals}
in the degenerate Eisenstein series case suggests
that this should not essentially take any new ideas. 

The next most natural extension is to complete
the treatment of the question of cuspidal exponents
and square integrability of all the $\Phi^{(b)}_i$.  In order to
accomplish this, various
methods, such as certain formulas generalizing the computations in \cite{krs}, have to be combined
with the irreducibility results of \cite{tadicisrael}.
This will be the subject of the immediate sequel to this paper.

For further possible applications of the inductive formula,
in particular to `first-term identities', which are relations
between the residues $\Phi_{i}^{(b)}$ for varying values of
both $i$ and $b$, one may see \cite{jiangfirstterm} and
\cite{jiangmemoirs}.

In the sequel to this paper, we intend to continue to work on the classification of the representations
generated by $\Phi^{(b)}_i$.  We wish to  establish results in parallel to the first few sections of Kudla and Rallis's paper
\cite{krannals}.  As will be seen from the examples considered towards the end 
of the paper, the further study of these representations depends on
some delicate \textit{local} questions.   For the case of degenerate principal series representations
(induced from a character on the Siegel parabolic), Kudla and Rallis
settled the local questions (of reducibility and constituents at points of reducibility) in 
\cite{krisrael}.  Certain existing papers in the literature
settle certain local questions in the case when the inducing data
is a Speh representation.  We particularly mention that in \cite{tadicisrael}
Tadi\'{c} has completely determined the points of reducibility in the cases
that interest us.  However, we do not know a full classification of the constituents
at points of reducibility, and it seems that at least some properties of these constituents
(such as non-singularity in the sense of Howe) must be established.

\vspace*{.3cm}\noindent
\textbf{Acknowledgments.}
The author heartily thanks Prof.\hspace*{-.15mm} Dihua Jiang, who suggested
the problem and provided constant encouragement and advice during
the work that lead to this paper.  He also thanks Lei Zhang, Ben Rosenfield, and 
Prof.\hspace*{-.15mm} Paul Garrett for stimulating conversations related to the paper.

\section{Generalities on Automorphic Forms, Eisenstein
Series, and Square-integrability}

\subsection{General notation.}
\label{subsec:generalnotation}
\noindent\textbf{Matrix model for the group $G_n=\mathrm{Sp}_{2n}$.}
Although we expect to generalize our results
to arbitrary split classical groups, we restrict ourselves
in this paper to the symplectic group $G_n$ of rank $n$.
Throughout, we use the following matrix model of the group
\[
G_n=\{ g\in\mathrm{GL}_{2n}\;|\; ^t g J_{2n}g=J_{2n} \},
\]
where
\[
J_{2n}=\begin{pmatrix}0& j_n\\
-j_n&0
\end{pmatrix}, j_n=\begin{pmatrix}
&&1\\
&\iddots&\\
1&&
\end{pmatrix}\in\mathrm{GL}(n)
\]

Henceforth, for $h\in\mathrm{GL}_r$, we will denote
$j_r{}^th^{-1}j_r$ by $\tilde{h}$. 

When denoting diagonal matrices we will write simply
$\mathrm{diag}(m_1,m_2,\ldots m_b)$ when the dimension
of the diagonal blocks $m_i$ is clear.  In some calculations
we will abbreviate further by using 
bold parenthesis, so that for example
$\boldop t,h\tilde{t}\boldcp:=\mathrm{diag}(t,h,\tilde{t})$ 

\vspace*{0.3cm}\noindent
\textbf{Parabolics and Levi Decompositions}.  
Let $a$ be a positive integer dividing $n$, so that $n=ab$.
Within $G=G_n=G_{ab}=\mathrm{Sp}_{2ab}$, we consider the following parabolic subgroups.  In order to define the notion of a  ``standard parabolic", we fix the Borel subgroup $B=P_0=M_0U_0$, where
\[
M_0=T_0=\{\mathrm{diag}(\lambda_1,\ldots, \lambda_n,\lambda_n^{-1},\ldots,\lambda_1^{-1})
\;|\; \lambda_i\in\mathrm{GL}_1\},
\]
and $U_0$, the unipotent radical of the Borel, consists of all upper triangular
matrices in $G_n$.  A \textbf{standard parabolic} is any
parabolic containing $P_0$.

For example, $P=P^{ab}_{ab}$ denotes the standard Siegel parabolic
in $\mathrm{Sp}_{2ab}$, with standard Levi
decomposition $P=MU$.  Here, $M\cong\mathrm{GL}_{ab}$,
because
\[
M=\{\mathrm{diag}(m,\tilde{m})\;|\; m\in\mathrm{GL}_{ab}\}.
\]  
Likewise $Q=P^{ab}_a$ denotes the standard
maximal Levi subgroup with Levi decomposition
\[
Q=L V,\; \text{where}\; L\cong \mathrm{GL}_a\times
G_{ab-a},
\]
and $V$ is the unipotent radical.  The standard Levi
decomposition is fixed, so that
\[
L=\{\mathrm{diag}(t,h,\tilde{t})\;|\; t\in\mathrm{GL}_a,\,
h\in G_{ab-a}\}.
\]

We also have need to consider non-maximal
parabolic subgroup of $G_n$.  For example, 
we also have the standard maximal parabolic
$P_{a,ab-a}$ in $\mathrm{GL}_{ab}$ with
Levi decomposition
\[
P_{a,ab-a}=M_{a,ab-a}N,\;\text{with}\; M_{a,ab-a}\cong \mathrm{GL}_a\times 
\mathrm{GL}_{ab-a}.
\]
By replacing the standard Levi $M_{ab}$ of the Siegel
parabolic by the parabolic $P_{a,ab-a}$, we obtain
the non-maximal standard parabolic $P_{a,ab-a}^{ab}$
of $G_{ab}$
with standard Levi decomposition
\[
P_{a,ab-a}^{ab}=M_{a,ab-a}^{ab}U_{a,ab-a}^{ab}\;\text{where}\; M_{a,ab-a}^{ab}\cong \mathrm{GL}_a\times 
\mathrm{GL}_{ab-a}.
\]
We also have the standard non-maximal parabolic
$P_{a^b}$ in $\mathrm{GL}_{ab}$ with
standard Levi decomposition
\[
P_{a^b}=M_{a^b}N_{a^b},\;\text{with}\; M_{a^b}\cong
M_{a^b}^{ab}\cong (\mathrm{GL}_a)^{\times b}, 
\]
giving rise, in the same manner to a non-maximal
parabolic
\[
P_{a^b}^{ab}=M_{a^b}^{ab}U_{a^b}
\]
in $G_{ab}$.

For readability, we denote by the same notation
for example, $P_{a,ab-a}$ and $P_{a^b}^{ab}$, pairs of parabolic subgroups
sharing the same standard Levi factor, but having
different nilpotent radicals depending
on whether they are sitting in the general linear,
or the classical groups.  This
practice has certain conveniences,
but may cause confusions when discussing
modulus characters.
To avoid such confusions, where they may arise, we denote
the modulus characters by the full notation
$\delta(G_n,P_n)$.  Other times,
we will use abbreviations such as $\delta_n$ for $\delta(G_n,P_n)$
and $\delta_{a^b}$ for $\delta(\mathrm{GL}_n,P_{a^b})$

Let $k$ be a number field and $\mathbf{A}$ its ring
of adeles.

We fix a maximal compact subgroup $K$ of $G(\mathbf{A})$
with the requirements of \S I.1--2 of \cite{mwbook} satisfied.
This allows us in particular to write the Langlands (generalized
Iwasawa) decomposition $g=muk$ of an element of $G$
with respect to any standard parabolic $P$ of $G$.
We frequently use notations such as $m(g)$ for the $M=M_{ab}^{ab}$-component of $g\in G_{ab}$.
In the full Langlands decomposition with respect to $P_{ab}^{ab}$,
we write  $g=m(g)u(g)k(g)$,
with $m\in M_{ab}^{ab}$, $u\in U_{ab}^{ab}$, and $k\in K$.   For other parabolic subgroups and their associated Langlands decompositions,
we often add indices such as $m_{a,ab-a}(g)$, and so on.

\vspace*{0.3cm}\noindent
\textbf{Cuspidal representations of $\mathrm{GL}_a$.}
Let $\tau$ be a irreducible \textit{unitary}, cuspidal automorphic
representation of $\mathrm{GL}_a(\mathbf{A})$.  We recall a classification of $\tau$ in terms of the existence of poles of the attached $L$-functions.  The classification will allow us, below,
to give precise conditions for the non-vanishing of the possible residues of the residual-data
Eisenstein series.  We follow pages 680--681 of \cite{gjrjams},
and the reader is referred to this paper for further context
and references.

The Rankin-Selberg $L$-function $L(s,\tau\times\tau)$
has a pole at $s=1$ if and only if $\tau\cong \tau^{\vee}$,
\textit{i.e.} $\tau$ is self-dual; and in that
case the pole is simple.  Therefore, we will
henceforth assume that $\tau$ is self-dual.  Further,
we have a factorization
\begin{equation}\label{equation:RSfactorization}
L(s,\tau\times\tau)=L(s,\tau,\wedge^2)\cdot L(s,\tau,\vee^2),
\end{equation}
of the Rankin-Selberg $L$-function into the product
of the exterior square and symmetric square $L$-functions
attached to $\tau$.  Exactly one of the two $L$-functions
$L(s,\tau,\wedge^2)$ and $L(s,\tau,\vee^2)$ has a simple
pole at $s=1$ and the other one is holomorphic and non-vanishing
at $s=1$.  If the exterior square $L$-function
$L(s,\tau,\wedge^2)$ has a simple pole at $s=1$,
which can only be the case if $a$ is even, then we say
that $\tau$ is \textbf{symplectic}, while
if the symmetric square $L$-function $L(s,\pi,\vee^2)$
has a simple pole at $s=1$, we say
that $\tau$ is \textbf{orthogonal}.  The reasons
for this terminology are explained in terms of Langlands
functorial lifts on page 680 of \cite{gjrjams}.
A full set of references to the literature found there.

We will also assume throughout that the central value
of the standard $L$-function, $L\left(\tau,\frac{1}{2}\right)$, is nonzero.

The $b$-fold tensor product $\tau^{\otimes b}$ will be
abbreviated by $\pi^b$ or even by $\pi$ when the context
is clear.

\vspace*{.3cm}\noindent
\textbf{Concerning the rational points of a variety $V$.}
  In a context in which both $V(k)$ and $V(\mathbf{A})$
arise, in order to distinguish them, we will often drop the specific reference to $k$
and just write $V$ for $V(k)$, while retaining
the specific reference to the field in $V(\mathbf{A})$.
In contexts where it may easily be discerned which
one we are talking about, we will drop all reference
to the field and write simply `$V$'.

\vspace*{0.3cm}\noindent
\textbf{Odds and Ends.}
Concerning discrete spectrum, noncuspidal representations,
we use the notation $\Delta(\tau,b)$ for the representations that Moeglin, in \textit{e.g.} \cite{moeglinmanuscripta}
refers to as $\mathrm{Speh}(\tau,b)$ (and also discussed 
more extensively in \cite{moeglinwaldspurgerens}).
Sometimes, when the context is clear
we abbreviate by $\Delta^b$, $\Delta^{b-1}$, and so on,
or even $\Delta$.
We will review the construction of these representations in more
detail \S\ref{subsec:residues}.

Generally speaking vectors will be represented
by bold typeface and numbers by roman.   In the context
where $s\in \mathbf{C}$ is a number, we will denote
by $\mathbf{s}$ the vector in $\mathbf{C}^b$ all whose entries are $s$.  Further when $\mathbf{s}$ is a general
vector in $\mathbf{C}^b$, $\mathbf{s}'$ will represent
a truncation of $\mathbf{s}$ to an element of $\mathbf{C}^{b-1}$.
Whether the truncation is achieved by dropping the first
element or the last ($b^{\rm th}$) element will be
made clear in the specific context.

The ``residue point" $\Lambda_b\in\mathbf{C}$ is special vector with (all real)
entries
\[
\Lambda_b=\left(\frac{b-1}{2},\frac{b-3}{2},\ldots, \frac{1-b}{2}\right).
\]
Alternatively, one can also denote $\Lambda_b$
by the `increasing segment' notation $\left[\frac{b-1}{2},\frac{1-b}{2}\right]$, where it is understood that a segment $[\alpha,\alpha+m]$
consists of the numbers lying between $\alpha$ and $\alpha+m$ differing from the endpoints by integer values.

\subsection{Cuspidal support, cuspidal exponents
and a criterion for square integrability of automorphic
forms.} For those parts of the general theory of automorphic
forms that are
directly relevant for stating and proving our principle
results, we now refer the reader to various parts of Chapter II of \cite{mwbook}.
We will follow the notation laid out in that book precisely
whenever possible.

\label{subsec:cuspidalexponents}
See \S I.3.3,  \S I.3.5 of \cite{mwbook}
for the notion of the \textbf{cuspidal support} $\Pi_0(M,\phi)$ of the
automorphic form $\phi$ along $M$, and the notion that
$\phi$is  \textbf{concentrated on} the set of standard 
parabolics $\{P\}$.

\vspace*{0.3cm}
\noindent \textbf{Spectral decomposition of $L^2(G(k)\backslash G(\mathbf{A}))$
according to the cuspidal support.}
We now recall main result of Chapter II of \cite{mwbook}, contained in
\S II.2.4.  Let $\xi$ be a central character of the center $Z_G$ of $G(\mathbf{A})$.  
Denote by $(M,\mathfrak{P})$, $(M',\mathfrak{P}')$ a pair of $X^G_M$-orbits
of cuspidal data (see \eqref{eqn:ourcuspidaldatabasepoint}, below, for more detail)
We consider $(M,\mathfrak{P})$ to be equivalent to $(M',\mathfrak{P}')$, and write
$(M,\mathfrak{P})\sim (M',\mathfrak{P}')$, if and only if there exists a $\gamma\in G(k)$
such that
\[
\gamma M \gamma^{-1}=M'\;\text{and}\; \gamma\mathfrak{P}=\mathfrak{P}'.
\]
An equivalence class of data 
with central character equal to $\xi$ is denoted by $\mathfrak{X}$.
For $\xi$ fixed, denote by $\mathfrak{E}$ the set of equivalence classes $\mathfrak{X}$.  For fixed
$\mathfrak{X}\in\mathfrak{E}$, denote by $L^2(G(k) Z_G\backslash G(\mathbf{A}))_{\mathfrak{X}}$
the closed subspace of $L^2(G(k)\backslash G(\mathbf{A}))$ generated by
\textit{pseudo-Eisenstein series} $\theta_{\phi}$, where $\phi$ runs over the Paley-Wiener space
$P_{(M',\mathfrak{X})}$ and $(M',\mathfrak{P}')$ runs over $\mathfrak{X}$.
(See \S\S II.1.2 and \S\S II.1.10 of \cite{mwbook} for definitions).  The decomposition
is stated in the following Proposition (II.2.4 of \cite{mwbook}).
\begin{prop} \label{prop:Edecomposition} With the preceding notations, one has
\[
L^2(G(k)\backslash \mathbf{G})_{\xi}=\hat{\bigoplus}_{\mathfrak{X}\in\mathfrak{E}}L^2(G(k)\backslash G(\mathbf{A})
_{\mathfrak{X}},
\]
where the sum is orthogonal, and the symbol $\hat{\cdot}$ indicates the Hilbert space completion.
\end{prop}
From the construction of the pseudo-Eisenstein series it follows that the cuspidal
support of any element $\phi\in L^2(G(k)\backslash G(\mathbf{A}))_{\mathfrak{X}}$,
lies in $\mathfrak{X}$.  Further, from Theorem II.1.12 (or more precisely, the proof
of the theorem), we can derive the opposite inclusion, so an alternate characterization
of the spaces in Proposition \ref{prop:Edecomposition},
\begin{equation}\label{eqn:l2spacefixedcuspidalsupport}
L^2(G(k)\backslash G(\mathbf{A})
_{\mathfrak{X}}=\left\{\phi\in L^2(G(k)\backslash G(\mathbf{A})\;\left|\; \bigcup_M \Pi_0(M,\phi)\subseteq \mathfrak{X}
\right.\right\}.
\end{equation}
Further, the decomposition of $L^2(G(k)\backslash \mathbf{G})_{\xi}$ can be summed over
the central characters $\xi$ and intersected with $L^2_d$, the discrete part of $L^2$,
(\textit{i.e.}, the part that is a direct sum of irreducible $G(\mathbf{A})$-modules) to yield the decomposition
\begin{equation}\label{eqn:discretepartdecomposition}
L^2_d(G(k)\backslash G(\mathbf{A}))=\bigoplus_{\mathfrak{X}\in \mathfrak{E}} L^2_d(G(k)\backslash G(\mathbf{A})
)_{\mathfrak{X}}.
\end{equation}
Langlands' Theory of Residues of Eisenstein series, recounted in Chapters V and VI of \cite{mwbook},
tells us how to produce, from cuspidal data $(M,\mathfrak{P})$,
certain automorphic forms,
which have cuspidal support in $\mathfrak{X}$ (the equivalence class of $(M,\mathfrak{P})$), which are orthogonal to both
the Eisenstein integrals (continuous part of the spectrum) and the cusp forms.  We wish to use this method to 
make the decomposition \eqref{eqn:discretepartdecomposition} explicit
in specific cases.
A crucial tool is a criterion for determining when the automorphic forms
constructed in this manner indeed lie in $L^2$.

\vspace*{.3cm}
\noindent  \textbf{Cuspidal Exponents
and the Criterion for square integrability.} 
Let  $\phi$ be an arbitrary automorphic
form on $G(k)\backslash G(\mathbf{A})$.  Then for each parabolic
$P$ of $G$ such that $\phi$ has nonzero cuspidal
support along $P$, the set of real-valued characters $\mathrm{Re}\pi\in\mathrm{Re} X_M^G$, 
\[
\{\mathrm{Re}\pi\}=\{\mathrm{Re}\chi_{\pi}\}=\{|\chi_{\pi}|\},
\; \text{as $\pi$ ranges over $\Pi_0(M,\phi)$}
\]
is called the \textbf{set of cuspidal exponents of $\phi$ along
$P$}.  For the statement of the criterion
for square integrability of an automorphic form in terms
of its cuspidal exponents see \S I.4.11 of \cite{mwbook}.
We will state and use the specific cases of this criterion as they arise
in our study of the residues.

\subsection{Eisenstein Series and their Meromorphic Continuation}
We first recall some of the basic notation and results of
\S\S I.3 and II.1 of \cite{mwbook}, and then specialize
to the cases that interest us.
Let $\pi$ be an unitary automorphic irreducible representation $M$.  
Denote by $(M,\mathfrak{P})$, the $X_{M}^G$-orbit 
of the ``datum" $(M,\pi)$.  
Using the notation of \S 1.3 of \cite{mwbook}, the $\mathrm{Re}X_M^G$-orbit of $(M,\pi)$ 
can be characterized as
\begin{equation}\label{eqn:ourcuspidaldatabasepoint}
\{(M,\rho)\in(M,\mathfrak{P}),\; \;|\; \mathrm{Im}\rho=\pi\}.
\end{equation}
Let $\phi_{\pi}$ denote an element of $A(U(\mathbf{A})M(k)\backslash G(\mathbf{A}))_{\pi}$,
obtained by multiplying a section of $\pi$ by 
the ``normalizing factor" $\delta^{\frac{1}{2}}(G,M)$
and then extending ``trivially" to a left $U(\mathbf{A})$-invariant
function on $G$.  For every $\lambda\in X_M^G$, one denotes by
\begin{equation}\label{eqn:xmgaction}\text{
$\lambda\phi_{\pi}$
the element $\lambda\circ m_P\cdot \phi_{\pi}\in A(U(\mathbf{A})M(k)\backslash G)_{\pi\otimes\lambda}$.}
\end{equation}

\begin{prop}\label{prop:convergencegeneral} (\textbf{II.1.5} of \cite{mwbook}.)
There is an open positive cone of $\mathrm{Re} X_{M}^G$ such that
for every element of $\mathfrak{P}$ with $\mathrm{Re}\pi$ in the cone, the series
\begin{equation}\label{eqn:eisensteingeneralseriesform}
E(\lambda\phi_{\pi},\pi\otimes\lambda)(g):=\sum_{\gamma\in P(k)\backslash G(k)}\lambda\phi_{\pi}(\gamma g)
\end{equation}
 converges absolutely for $g$ in a
compact set and $\lambda$ in a neighborhood of $0$ in $X_{M}^G$.  So the series
defines an automorphic form on $G(k)\backslash G$.  Supposing that $\mathfrak{P}$
is formed of cuspidal representations then the cone is
\begin{equation}\label{eqn:conecuspidalgeneral}
\mathfrak{c}_{\mathfrak{P}}:=
\{\lambda\in\mathrm{Re}X_M^G\;|\; \langle\lambda,\alpha^{\vee}\rangle>
\langle \rho_P,\alpha^{\vee}\rangle, \;\forall \alpha\in\Delta^+(T_M,G)\rangle\}.
\end{equation}
\end{prop}

\vspace*{0.3cm}
\noindent
\textbf{Holomorphy of $E$ in the cone $\mathfrak{c}_{\mathfrak{P}}$
and meromorphic continuation to $\mathfrak{P}$.}
Now fix a character $\xi$ of $Z_G$ and a finite
set $\mathfrak{F}$ of $K$-types which is assumed
stable under passage to the contragredient.
Denote by
\[
A(M,\pi)^{\mathfrak{F}}:=A(U(\mathbf{A})M(k)\backslash 
G(\mathbf{A}))^{\mathfrak{F}}_{\pi},
\]
the space of automorphic forms of type $\pi$ transforming
under $K$ according to $\mathfrak{F}$.
Define $A^{\mathfrak{F}}_{\xi}=A(G(k)\backslash G(\mathbf{A}))
^{\mathfrak{F}}_{\xi}$, so that 
$A(M,\pi)^{\mathfrak{F}}\subset A^{\mathfrak{F}}_{\xi}$.

Now assume that $\xi$ is the central character of
each element of $(M,\mathfrak{P})$.

For each $\pi\in\mathfrak{c}_{\mathfrak{P}}$ and $\varphi\in A(M,\pi)^{\mathfrak{F}}$
(see \S IV.1.1 for notation), let $E(P,\varphi,\pi)$
be defined by the series.  One has $E(P,\varphi,\pi)\in A^{\mathfrak{F}}_{\xi}$
and so a fortiori $E(P,\pi)\in L^{2,\mathfrak{F}}_{\xi,\mathrm{loc}}$.  This defines on $\mathfrak{c}_{\mathfrak{P}}$
a function
\[
E(P): \pi\mapsto E(P,\pi)\in\mathrm{Hom}_{\mathbf{C}}(A(M,\pi)^{\mathfrak{F}},L_{\xi,\mathrm{loc}}^{2,\mathfrak{F}}).
\]

\vspace*{0.3cm}
\noindent

The definition of \textit{holomorphic} in this context is given in \S IV.1.3 of \cite{mwbook} as follows.  For $\lambda\in X^G_{M}$,
the action \eqref{eqn:xmgaction} defines an isomorphism of vector space
\[
\underline{\lambda}: A(M,\pi)^{\mathfrak{F}}\rightarrow A(M_,\pi\otimes\lambda)^{\mathfrak{F}}.
\]
Let $U$ be an open subset of $\mathfrak{P}$, $H$ a Fr\'{e}chet space and $E$ a function
defined on $U$ such that for every $\pi\in U$,
\begin{equation}\label{eqn:Edefn}
E(\pi)\in \mathrm{Hom}_{\mathbf{C}}(A(M,\pi)^{\mathfrak{F}},H).
\end{equation}
Since $A(M,\pi)^{\mathfrak{F}}$ is finite-dimensional, this homomorphism
space is also a Fr\'{e}chet space.  We will say that $E$ is holomorphic
if for every $\pi\in U$, the function
\[
\lambda\mapsto E(\pi\otimes \lambda)\otimes \underline{\lambda}
\]
defined in a neighborhood of zero in $X^G_{M_{a^b}}$ is holomorphic.  
\vspace*{0.3cm}
\begin{lemma} \label{eqn:holomorphicityofthesum} The function
\[
E(P): \pi\mapsto E(P,\pi)\in \mathrm{Hom}_{\mathbf{C}}(A(M_{a^b},\pi)^{\mathfrak{F}},L_{\xi,\rm loc}^{2,\mathfrak{F}})
\]
is holomorphic on $\mathfrak{c}_{\mathfrak{P}}$.
\end{lemma}
\noindent\textbf{Comments on Proof.}  In practice, as noted
in the discussion in \S I.4.9, the holomorphy
that has to be proved is that for every $\pi\in\mathfrak{c}_{\mathfrak{P}}$, $\varphi\in A(M,\pi)^{\mathfrak{F}}$
compact subset $C\subset\mathfrak{P}$, and $\psi\in L^2_{\xi,C}$, the mapping
\[
\lambda\rightarrow\int_{Z_G(k)\backslash G}\overline{\psi}(g) E(\lambda\varphi,\pi\otimes\lambda)(g)\,
\mathrm{d} g
\]
is holomorphic at $0$.  

This holomorphy follows from the absolute convergence of the series \eqref{eqn:eisensteingeneralseriesform}
on $\mathbf{c}_{\mathfrak{P}}$.
\qed

Now, following p. 137 of \cite{mwbook}, we define
\textit{meromorphic} in this context.  Suppose that
$E$ is only defined almost everywhere on $U$.  We
say that $E$ is meromorphic if for every $\pi\in U$
and every sufficiently small neighborhood $V$ of $\pi$
in $U$, there exist two holomorphic functions
\[
d: V\rightarrow\mathbf{C},\quad d\neq 0,
\]
and $E_1$ defined on $V$, and taking values in $\mathrm{Hom}_{\bf C}(A(M,\pi)^{\mathfrak{F}},H)$ such that
\[
d(\pi')E(\pi')=E_1(\pi')
\]
for every $\pi'\in V$ where $E(\pi')$ is defined.
\begin{thm} \label{eqn:meromorphiccontinuationgeneral}(\textbf{IV.1.8}(a) in \cite{mwbook}.)  Let $(M,\mathfrak{P})$ be as above.  Then
$E(P)$ extends in a unique manner to a meromorphic function on $\mathfrak{P}$.
\end{thm}
The proof of this and several related results, some of which we will also need, is the subject
of Chapter IV of \cite{mwbook}.  
The principle further properties of the meromorphically continued
Eisenstein series are that it still, just as at the points of the convergence of the series,
is a function of moderate growth (hence an automorphic form) and is
orthogonal to cusp forms.  More precisely,

\begin{prop}\label{prop:propsofeisseries}  
\textbf{IV.1.9(b) of \cite{mwbook}}
Let $(M,\mathfrak{P})$ be as above. Let $U$ be the open set of $\mathfrak{P}$ where
$E(P)$ is holomorphic.
\begin{itemize}
\item[(i)]  For $\pi\in U$ and $\varphi\in A(M,\pi)^{\mathfrak{F}}$,
$E(P,\varphi,\pi)\in A^{\mathfrak{F}}_{\xi}$.
\item[(ii)]  For $\pi\in U$ and $\varphi\in A(M,\pi)^{\mathfrak{F}}$,
$E(P,\varphi,\pi)^{\rm cusp}=0$ if $M\neq G$.
\end{itemize}
\end{prop}

\vspace*{0.3cm}
\noindent
\textbf{Examples of Cuspidal-Data Eisenstein series and their cones of convergence.}
The above immediately applies to two examples of Eisenstein series
that interest us, and allows us to compute the cones of convergence precisely.
Let $\tau$ be a fixed unitary, irreducible cuspidal automorphic representation
of $\mathrm{GL}_a$.  Generally it will be assumed that $\tau$ self-dual, but
this is not necessary in the present context.  Then let $\pi=\pi^b=\tau^{\otimes b}$
be the $b$-fold tensor product of $\tau$, a representation of the same
type, of $M_{a^b}\cong \mathrm{GL}_a^{\times b}$.

Let the $f_i$, $i=1,\cdots, b$ be the natural
set of coordinates on $\mathrm{Re}\mathfrak{a}^*_{M_{a^b}}$,
in terms of which
\begin{equation}\label{eqn:restrictedrootsexplicit}
\Delta^+(T_{M_{a^b}},G)=\{f_1-f_2,f_2-f_3,\ldots, f_{b-1}-f_b,2f_b\}.
\end{equation}
In terms of the coordinates imposed by the $f_i$, one calculates that
\begin{multline}\label{eqn:rhoabcalculation}
\rho_b^{(a)}:=\rho(G,P_{a^b})=\left(a\left(b-\frac{1}{2}\right)+\frac{1}{2},
\left(\rho_{b-1}^{(a)}\right)_1,\ldots \left(\rho_{b-1}^{(a)}\right)_{b-1}\right)=\\
\left(a\left(b-\frac{1}{2}\right)+\frac{1}{2}, a\left(b-\frac{3}{2}\right)+\frac{1}{2},\ldots, \frac{a+1}{2}\right).
\end{multline}
Now define the cone
\[\label{eqn:coneforcuspdataeisseries}
\mathfrak{c}_{\mathfrak{P}}:=\left\{r_1f_1+\cdots r_b f_b\in \mathrm{Re}X_{M_{a}^b}^G \;\left|\;
r_b>\frac{a+1}{2};\quad r_{i}-r_{i+1}>a,\; \text{for}\; 1\leq i\leq b-1\right.\right\}.
\]
From \eqref{eqn:rhoabcalculation}, it is easy to see that
$\mathfrak{c}_{\mathfrak{P}}$ is the special case of \eqref{eqn:conecuspidalgeneral} 
for $G=G_n$ and $P=P_{a^b}$.  Therefore, denoting by $\mathfrak{P}$
the $X_{M_{a^b}}^G$-orbit of the datum $(M_{a^b},\pi)$, Proposition \ref{prop:convergencegeneral} applies
to say that for any $\pi\otimes\lambda\in \mathfrak{P}$ such that $\pi\otimes\lambda\in \mathfrak{c}_{\mathfrak{P}}$, 
the sum 
\[
E(P_{a^b},\lambda\phi_{\pi},\pi\otimes\lambda)(g):=\sum_{\gamma\in P_{a^b}(k)\backslash G(k)}\lambda\phi_{\pi}(\gamma g)
\]
converges absolutely and uniformly on compact sets.  Further, by Lemma \ref{eqn:holomorphicityofthesum},
Theorem \ref{eqn:meromorphiccontinuationgeneral} and Proposition \ref{prop:propsofeisseries},
the associated homomorphism-valued function $E(P_{a^b})$ on $\mathfrak{c}_{\mathfrak{P}}$ has a unique
meromorphic continuation to the complement $U$ of a finite set of $\mathfrak{P}$, such that
\[
 \pi\in U,\; \varphi\in A(M_{a^b},\pi)^{\mathfrak{F}} \;\Rightarrow\;
E(P_{a^b},\varphi,\pi)\in (A_0(G(k)\backslash G)^{\mathfrak{F}}_{\xi})^{\perp}.
\]
Here, the orthogonality symbol indicates the orthogonal complement 
(in this case, of the cuspidal subspace) in $A^{\mathfrak{F}}_{\xi}$.

Now view $M_{a^b}$ as a Levi for the parabolic subgroup $P_{a^b}$ of $M\cong \mathrm{GL}_{ab}$,
as discussed above.  
In terms of the same coordinates $f_i$, the restricted root
system $R(T_{M_{a^b}},M)$ has simple system
\[
\Delta^+(T_{M_{a^b}},M)=\{f_i-f_{i+1}\;|\; i=1,\ldots b-1\}.
\]
One calculates that
\begin{equation}\label{eqn:rhosubabcalculation}
\rho(\mathrm{GL}_{ab},P_{a^b}):=\rho_{a^b}=\left(a\frac{b-1}{2},a\frac{b-3}{2},\ldots, a\frac{1-b}{2} \right).
\end{equation}
For later use, we note the relation
\begin{equation}\label{eqn:normfactorsrel}
\rho_b^{(a)}=\rho_{a^b}+\mbox{\boldmath $\rho_{ab}$\unboldmath}.
\end{equation}

Now we define the cone
\begin{equation}\label{eqn:generallinearcondefn}
\mathfrak{c}_{\mathfrak{P}_{a^b}}:=\{\lambda_1f_1+\cdots \lambda_b f_b\in\mathrm{Re} X_{M_{a^b}}^M\;|\;
\lambda_i-\lambda_{i+1}>a\;\text{for}\; i=1,\ldots b-1\}.
\end{equation}
From \eqref{eqn:rhoabcalculation}, it is easy to see that
$\mathfrak{c}_{\mathfrak{P}_{a^b}}$ is the special case of \eqref{eqn:conecuspidalgeneral} 
for $G=M$ and $P=P_{a^b}$.  Therefore, denoting by $\mathfrak{P_{a^b}}$
the $X_{M_{a^b}}^M$-orbit of the datum $(M_{a^b},\pi)$, Proposition \ref{prop:convergencegeneral} applies
to say that for any $\pi\otimes\lambda\in \mathfrak{P}_{a^b}$ such that 
$\mathrm{Re}\pi\otimes\lambda\in \mathfrak{c}_{\mathfrak{P}_{a^b}}$, 
the sum 
\begin{equation}\label{eqn:McuspidalEissum}
E^{M}(\lambda\phi_{\pi},\pi\otimes\lambda,g)=\sum_{\gamma\in P_{a^b}\cap M(k)\backslash M(k)} \lambda\phi(\gamma g)\;
\end{equation}
converges absolutely and uniformly on compact sets.  Further, by Lemma \ref{eqn:holomorphicityofthesum},
Theorem \ref{eqn:meromorphiccontinuationgeneral} and Proposition \ref{prop:propsofeisseries},
the associated homomorphism-valued function $E^M(P_{a^b})$ on $\mathfrak{c}_{\mathfrak{P}_{a^b}}$ has a unique
meromorphic continuation to the complement $U$ of a finite set of $\mathfrak{P}_{a^b}$, such that
\begin{equation}\label{eqn:propseisseriesglapp}
 \pi\in U,\; \varphi\in A(M_{a^b},\pi)^{\mathfrak{F}} \;\Rightarrow\;
E^M(P_{a^b},\varphi,\pi)\in (A_0(M(k)\backslash M)^{\mathfrak{F}}_{\xi})^{\perp}.
\end{equation}
 
\vspace*{0.3cm}\noindent
\textbf{An Eisenstein Series with non-cuspidal, non-discrete data} 
 By imposing right-$U(\mathbf{A})$-invariance, we can extend $E^{M}(P_{a^b})$ to a function
on $G=\mathrm{Sp}_{2ab}$.
  Then we formally define the $G$-Eisenstein series
\begin{equation}\label{eqn:eisseriesdataeisseriesdefined}
E(P,E^M(P_{a^b},\lambda\phi,\pi\otimes\lambda),g,s):=
\sum_{\gamma\in P\backslash G} E^M(P_{a^b},\lambda\phi,\pi\otimes\lambda)(\gamma g)m_{P}(\gamma g)^{\rho_{ab}+s}.
\end{equation}
Although this is an Eisenstein series, defined with a datum that is itself an Eisenstein
series, thus in the orthogonal complement of cusp forms, is not as familiar
as Eisenstein series with data from the discrete spectrum, we claim
that all the main results concerning convergence, holomorphy, and meromorphic
continuation still pertain.  In order to state the appropriate form of these results,
first define the following subset of $\mathrm{Re}X_{M}^G\cong \mathbf{C}$, 
where $M$ is the Siegel parabolic $M_{ab}^{ab}$ of $G_{ab}$:
\[
\mathfrak{c}_{\Delta}=\left\{ s\in\mathbf{C}\;\left|\; \mathrm{Re}s>\rho_{ab}=\frac{ab+1}{2}\right.\right\}.
\] 
The reason for the notation $\Delta$ will become clear, below.

\begin{prop} \label{prop:eisdataeisseries} Let $\lambda_0\pi \in U$, the dense open subset of $\mathfrak{P}_{a^b}$ where $E^M(P_{a^b})$
is holomorphic.  Set $\lambda=\lambda_0+\lambda'$.  For $s$ fixed in $\mathfrak{c}_{\Delta}$
and $\lambda\phi\in A(M_{a^b},\lambda\pi)^{\mathfrak{F}}$,
the series 
\[
E(P,E^M(P_{a^b},\lambda\phi,\pi\otimes\lambda),g,s)
\]
defined by \eqref{eqn:eisseriesdataeisseriesdefined} converges absolutely and uniformly for $g$
in a compact set and $\lambda'$ in a neighborhood of $0$ in $X^M_{M_{a^b}}$.  Further, the series
defines an automorphic form on $G(k)\backslash G(\mathbf{A})$.
\end{prop}

\textbf{Proof.}  
Implicit in the conclusion of Proposition \ref{prop:propsofeisseries} (i), we have 
the moderate growth of the function
\[
m\in \mathrm{GL}_{ab} \rightarrow E^M(P_{a^b},\lambda\phi,\pi\otimes\lambda).
\]
Therefore, we have the moderate growth of the function
\[
m \in\mathrm{GL}_{ab}\rightarrow m^{s} E^M(P_{a^b},(\lambda\phi,\pi\otimes\lambda)).
\]
In terms of Siegel domains the discussion of \S I.2.3 of \cite{mwbook}
says that the moderate growth of the function may be expressed as follows:
if $\mu$ is \textit{any} element of $\mathrm{Re}_{M_0}^M$, just so long
as $\mu$ is sufficiently positive, then there exists $c$ such that for all $k\in K$ and for
all $m\in M^1\cap S^{\mathrm{GL}_{ab}}$,
\begin{equation}\label{eqn:moderategrowthest}
|m^{s+\rho_{ab}} E^M(P_{a^b},\lambda\phi,\pi\otimes\lambda)(mk)|\leq
c m_{P_0}(m)^{\mu+\rho_0}.
\end{equation}
Here, $\rho_0$ is an abbreviation for $\rho_{P_0}=\rho(\mathrm{GL}_{ab},P_0)$.
In terms of the coordinates $e_i,\, i=1,\ldots, ab$,
\begin{equation}\label{eqn:rhozeroexplicit}
\rho_0=\left(\frac{ab-1}{2},\frac{ab-3}{2},\ldots,\frac{1-ab}{2}\right).
\end{equation}
(Remark that since $m\in M^1$, the first factor $m^{s+\rho_{ab}}$ in \eqref{eqn:moderategrowthest} is actually $1$.)
Let
\[
y=y_{\mu}: G\rightarrow \mathbf{R}_+^*
\]
be a function such that for all $u\in U_0(\mathbf{A})$, all $m\in M_0$,
all $g\in G$ and $k\in K$,
\begin{equation}\label{eqn:ydefinition}
y(umg)=y_{\mu}(umg)=m^{\mu+\rho_0}y(g)\;\text{and}\; 1<y(k).
\end{equation}
Define
\[
\mathfrak{c}(M,P_0)=\{\mu\in\mathrm{Re}_{M_0}^M\;|\; \langle\mu,\alpha^{\vee}  \rangle>\langle\rho_0,\alpha^{\vee}  
\rangle,\, \;\alpha=e_i-e_{i+1},\, 1\leq i< ab\}.
\]
From \eqref{eqn:rhozeroexplicit}, it follows that this cone is described explicitly in coordinates as
\begin{equation}\label{eqn:glnminimalconeexplicit}
\mathfrak{c}(M,P_0)=\{\mu\in\mathrm{Re}_{M_0}^M\;|\; \mu_i-\mu_{i+1}>1\}.
\end{equation}
\textit{Assume from now on that $\mu\in\mathfrak{c}(M,P_0)$}.
For all $g\in G$, set
\[
E^M(P_{a^b},y,g):=\sum_{\gamma\in M(k)\cap P_0(k)\backslash M(k)}
y(\gamma g)
\]
By Theorem 3 of Godement's Bourbaki article, this series
converges and defines an element of $A(U(\mathbf{A})M(k)\backslash G)$.
Moreover, $E^M(P_{a^b},y,g)$ is a series with only positive terms, one of
which is $y(g)$, so
\begin{equation}\label{eqn:yeisseriestriviallowerbound}
E^M(P_{a^b},y,g)\geq y(g),\;\text{for all}\; g\in G.
\end{equation}
Now we have the crucial estimate, for $m\in M^1\cap S^{\mathrm{GL}_{ab}}$,
\[
m_{P_0}(m)^{\mu+\rho_0}=cy(mk)y(k)^{-1}\leq c y(mk)\leq c E^M(P_{a^b},y,mk).
\]
where the equality and the first inequality follow directly from 
\eqref{eqn:ydefinition}.  The second inequality follows from 
\eqref{eqn:yeisseriestriviallowerbound}.  Therefore, by 
\eqref{eqn:moderategrowthest} we have
\begin{equation}\label{eqn:eisserieseisseriesest}
|m^s E^M(P_{a^b},\lambda\phi,\pi\otimes\lambda)(mk)|\leq
 c E^M(P_{a^b},y,mk)
 \end{equation}
for all $m\in M^1\cap S^{\mathrm{GL}_{ab}},\;
 \text{and}\; k\in K$.

Using the left-invariance of under $M(k)$ of both sides
of this inequality, we obtain \eqref{eqn:eisserieseisseriesest}
for all $m\in M(k)(M^1\cap S^{\mathrm{GL}_{ab}})$.  Recall
that there is a compact set $\omega$ of  
$P_0$ (the fixed minimal parabolic of $\mathrm{GL}_{ab}$) involved in the
definition of the Siegel set $S^{\mathrm{GL}_{ab}}$.  Taking $\omega$ sufficiently
large, we can arrange to have
\[
M(k)(M^1\cap S^{\mathrm{GL}_{ab}})=M^1,
\] 
so we obtain \eqref{eqn:eisserieseisseriesest} for all $m\in M^1$
and all $k\in K$. 

Now, let $z\in Z(M)$, the center of $M$.  According to \eqref{eqn:propseisseriesglapp}
\begin{equation}\label{eqn:centralchargenlineisseries}
E^M(P_{a^b},\pi,\phi_{\pi})\in A^{\mathfrak{F}}_{\xi},
\end{equation}
where $\xi$ is the restriction of the central character of any
element of $\mathfrak{P}_{a^b}$ to $Z(M)$.
By \eqref{eqn:centralchargenlineisseries},
\[
E^M(P_{a^b},\pi,\phi_{\pi},zmk)=z^{\xi} E^M(P_{a^b},\pi,\phi_{\pi},mk).
\]
In particular, for $z\in A_M\subset Z_M$ we have the estimate
\begin{multline}\label{eqn:cuspidaldatadroppingcenter}
|(zmk)^{s+\rho_{ab}}E^M(P_{a^b},\pi,\phi_{\pi},zmk)|\leq z^{\mathrm{Re}(\mathrm{s})+\rho_{ab}+\mathrm{Re}\xi}|m^sE^M(P_{a^b},\pi,\phi_{\pi},mk)|\\=z^{\mathrm{Re}s+\rho_{ab}}|m^sE^M(P_{a^b},\pi,\phi_{\pi},mk)|,
\end{multline}
the latter equality arising because $\xi$ is a unitary character.
On the other hand, from the series definition of $E^M(P_{a^b},y)$ above and 
the defining conditions \eqref{eqn:ydefinition} of $y$, we have
\begin{equation}\label{eqn:degeisseriesaddingcenter}
E^M(P_{a^b},y,zmk)=z^{\mu+\rho_0}E^M(P_{a^b},y,mk)=E^M(P_{a^b},y,mk),
\end{equation}
where the second equality follows because by construction $\mu, \rho_0\in X_{M_0}^M$.  Togehter, the estimates 
\eqref{eqn:eisserieseisseriesest}, \eqref{eqn:cuspidaldatadroppingcenter},
and \eqref{eqn:degeisseriesaddingcenter} imply that
\[
|(zmk)^{s+\rho_{ab}}E^M(P_{a^b},\pi,\phi_{\pi},zmk)|\leq c z^{\mathrm{Re}s+\rho_{ab}}E^M(P_{a^b},y,zmk).
\]
 Because
of the factorization $M=A_M M^1$ (in number field case: see top of p. 20 in \cite{mwbook}),
this inequality yields the inequality
\[
|(mk)^{\rho_{ab}+s}E^M(P_{a^b},\pi,\phi_{\pi},mk)|\leq c m^{\mathrm{Re}s+\rho_{ab}}E^M(P_{a^b},y,mk),\;\text{for all}\;
m\in M\, k\in K.
\]
Each side is also invariant by $U(\mathbf{A})$, the $y$-Eisenstein series, because
it is the sum of terms which are left-$U(\mathbf{A})$-invariant by definition,
and the cuspidal-data Eisenstein series, because it is \textit{defined} that
way as an element of $A(U(\mathbf{A})M(k)\backslash G)$.  Therefore,
for all $g\in G$ we have
\[
|m_P(g)^{s+\rho_{ab}}E^M(P_{a^b},\pi,\phi_{\pi},g)|\leq c m_P(g)^{\mathrm{Re}s+\rho_{ab}}E^M(P_{a^b},y,g).
\]
Because of this estimate we have
\begin{multline*}
\sum_{\gamma\in P\backslash G} |E^M(P_{a^b},\lambda\phi,\pi\otimes\lambda)(\gamma g)m_{P}(\gamma g)^{\rho_{ab}+s}|\\
\leq
c \sum_{\gamma\in P\backslash G}  m_P(\gamma g)^{\mathrm{Re}s+\rho_{ab}}E^M(P_{a^b},y_{\mu},\gamma g)=
\\
c \sum_{\gamma\in P(k)\backslash G(k)}  m_P(\gamma g)^{\mathrm{Re}(s)+\rho_{ab}}
\sum_{\gamma'\in M(k)\cap P_0(k)\backslash M(k)}
y_{\mu}(\gamma' \gamma g)=\\
\sum_{\gamma\in P_0\backslash G} m_{P}(\gamma g)^{\mathrm{Re}s+\rho_{ab}}y_{\mu}(\gamma g)=
\sum_{\gamma\in P_0\backslash G} y_{\mu'}(\gamma g)
\end{multline*}
where
\[
\mu'=\mu+\rho_{0}+\mathrm{Re}s+\rho_{ab}-\rho(G,P_0)
\]
Applying \eqref{eqn:rhoabcalculation} to calculate $\rho(G,P_0)=\rho^{(1)}_{ab}$,
and using \eqref{eqn:rhozeroexplicit} to see that $\rho_0+\rho_{ab}=\rho_{ab}^{(1)}$, 
we calculate that actually
$\mu'=\mu+\mathrm{Re}s$.  Together, the hypotheses of the proposition that $s\in\mathfrak{c}_{\Delta}$,
and the assumption just after \eqref{eqn:glnminimalconeexplicit} 
that $\mu\in \mathfrak{c}(M,P_0)$ imply precisely that $\mu'\in \mathfrak{c}(G,P_0)$!
Therefore, the right-hand side of the estimate converges to the minimal-parabolic
(degenerate) Eisenstein series, and 
\[
\sum_{\gamma\in P\backslash G} |E^M(P_{a^b},\lambda\phi,\pi\otimes\lambda)(\gamma g)m_{P}(\gamma g)^{\rho_{ab}+s}|
\leq E^G(P_0,y_{\mu'},g).
\]
The convergence and uniformity in $g$ now follow from Theorem 3 of \cite{godement}.
The uniformity in $\lambda'$ follows from elementary properties of holomorphic functions.
\qed

The convergence established in Proposition \ref{prop:eisdataeisseries}, together with
the general meromorphic continuation result, Theorem \ref{eqn:meromorphiccontinuationgeneral},
 imply through a simple calculation in the cone of absolute convergence,
 that although \textit{a priori} the $P$-Eisenstein series
under consideration above is non-cuspidal data, it can actually be identified with a cuspidal-data
Eisenstein series from the smaller parabolic $P_{a^b}$, with an appropriate parameter shift.

\begin{prop}  Let $s\in\mathfrak{c}_{\Delta}$.  Let $U$ be the open subset of ${P_{a^b}}$
where $E^M(P_{a^b})$ is holomorphic.  Then
\begin{equation}\label{eqn:eisensteinidentity}
E(P,E^M(P_{a^b},\lambda\phi,\pi\otimes\lambda),g,s)=E(P_{a^b},(\lambda+\mathbf{s})\phi,\pi\otimes(\lambda+\mathbf{s}),g).
\end{equation}
\end{prop}
\textbf{Proof}  First suppose that $\lambda\in\mathfrak{c}_{P_{a^b}}$.  Then $E^M(P_{a^b},\lambda\phi,\pi\otimes\lambda)$
is given by the sum \eqref{eqn:McuspidalEissum}, which we can then substitute into
the right-hand side of \eqref{eqn:eisseriesdataeisseriesdefined}.  
\[
E(P,E^M(P_{a^b},\lambda\phi,\pi\otimes\lambda),g,s):=
\sum_{\gamma\in P\backslash G} \sum_{\gamma'\in M_{a^b}\backslash M}
\lambda\phi(m_P(\gamma'\gamma g)) m_{P}(\gamma g)^{\rho_{ab}+s}
\]
We can take the representatives of $M_{a^b}\backslash M$ to have
determinant one (without loss of generality).  We further
apply the definition of $A(N_{a^b}(\mathbf{A})M_{a^b}(k)\backslash\mathrm{GL}_{ab}(\mathbf{A}))$ to obtain
that the right-hand side of the above equals
\[
\sum_{\gamma\in P\backslash G} \sum_{\gamma'\in M_{a^b}\backslash M}
\lambda\phi(m_{P_{a^b}}\gamma'\gamma g)
m_{P_{a^b}}(\gamma'\gamma g)^{\rho_{a^b}} m_{P}(\gamma'\gamma g)^{\rho_{ab}+s}.
\]
Combining the two sums, we then obtain
\begin{multline*}
\sum_{\gamma\in M_{a^b}\backslash G}
\lambda\phi(m_{P_{a^b}}(\gamma g))m_{P_{a^b}}(\gamma g)^{\mbox{\boldmath ${\scriptstyle \rho_{ab}}
$\unboldmath}+\mathbf{s}+\rho_{a^b}}=
\sum_{\gamma\in M_{a^b}\backslash G}
\lambda\phi(m_{P_{a^b}}(\gamma g))m_{P_{a^b}}(\gamma g)^{\mathbf{s}+\rho_b^{(a)}}=
\\
\sum_{\gamma\in M_{a^b}\backslash G}\lambda\phi(\gamma g)=
E(P_{a^b},(\lambda+\mathbf{s})\phi,\pi\otimes(\lambda+\mathbf{s}),g).
\end{multline*}
We have thus established the identity \eqref{eqn:eisensteinidentity} for 
$\lambda\in\mathfrak{c}_{\mathfrak{P}_{a^b}}$.  The identity for all $\lambda\in\mathfrak{P}_{a^b}$
 then follows
from the uniqueness of the meromorphic continuation of the left
hand side of \eqref{eqn:eisensteinidentity} to all of $\mathfrak{P}_{a^b}$ and of the right hand
side to the affine subspace $\mathbf{s}+\mathfrak{P}_{a^b}$ of $\mathfrak{P}$.
\qed

\subsection{Residues of Meromorphically Continued Eisenstein Series}
\label{subsec:residues}
We continue to follow the notation \cite{mwbook}, now mostly \S V.1,
but without recalling all the details of the definitions.  Consider
the complete flag of affine subspaces $\mathfrak{D}=\mathfrak{D}_{\Delta,\Lambda_b}$
(here $\Delta=\Delta(\mathrm{GL}_{ab},P_{a^b})$ denotes
the simple restricted root system and $\Lambda_b$ the intersection point) of $\mathfrak{P}_{a^b}$
defined as follows.
Let $\pi_0$ be the fixed element of $\mathfrak{P}_{a^b}$ defined as
\begin{equation}\label{eqn:basepointdefn}
\pi_0=\pi^b\otimes\Lambda_b.
\end{equation}
Then set 
\[
\mathfrak{D}_{\Delta,\Lambda_b}=\{\mathfrak{S}_0=\{\pi_0\}\subset \mathfrak{S}_1\subset \cdots\subset \mathfrak{S}_{b-1}=\mathfrak{P}_{a^b}\},
\]
where for each $i$ from $1$ to $b-1$, the vector parts of the affine hyperplanes are defined by
\begin{equation}\label{eqn:hyperplanesdefn}
\mathfrak{S}_{i-1}^0:=\mathfrak{S}_{i}^0 \cap H_{\alpha_{b-i}^*,\pi_0},\;\text{where}\; \alpha_{b-i}:=e_{b-i}-e_{b-i+1}\in
\Delta(\mathrm{GL}_{ab},P_{a^b}).
\end{equation}
Note that in particular we have
\begin{equation}\label{eqn:codimensiononestatement}
\mathrm{dim}_{\mathbf{R}}(\mathrm{Re}\mathfrak{S}_i/\mathrm{Re}\mathfrak{S}_{i-1})=1.
\end{equation}

We now define the ``residue datum" from $\mathfrak{P}_{a^b}$ to $\Lambda_{b}$ taken along $\mathfrak{D}$
to be the $(b-1)$-fold composition of certain more elementary residue data
\[
\mathrm{Res}_{\Lambda_b}^{\mathfrak{P}_{a^b}}:=\mathrm{Res}_{b-1}\circ \mathrm{Res}_{b-2}\circ
\cdots\circ \mathrm{Res}_{1}.
\]
Here, $\mathrm{Res}_i$, or more fully, $\mathrm{Res}^{\mathfrak{S}_{b-i}}_{\mathfrak{S}_{b-i-1}}$ are the operators
defined using the method of \S V.1.3 of \cite{mwbook}.  First, the operator
$\mathrm{Res}_i$ takes a certain space of meromorphic functions $A_i$
on $\mathfrak{S}_i$ to a certain space of meromorphic functions $A_{i-1}$ on $\mathfrak{S}_{i-1}$.
The further condition that these function satisfy is called \textit{having polynomial
singularities $S_{\mathfrak{X}}$}.  For the exact definition of this condition
see \cite{mwbook}.  All that concerns us here is that if we set
\[
A_{b-1}(\pi)=E^M(P_{a^b},\pi),
\]
then $A_{b-1}$ has polynomial singularities on $S_{\mathfrak{X}}$, and that consequently, for $i=2,\ldots b-1$, if we set
\[
A_{b-i}=\mathrm{Res}_{b-i+1} A_{b-i+1}.
\]
Then $A_{b-i}$ likewise has polynomial singularities $S_{\mathfrak{X}}$.

Proceeding with the definition of $\mathrm{Res}_i$, fix a nonzero element
\[
\epsilon_i\in (\mathrm{Re}\mathfrak{S}_{i-1}^0)^{\perp}\cap\mathrm{Re}\mathfrak{S}_i^0,\; i=1,\ldots, b-1,
\]
which by \eqref{eqn:codimensiononestatement} is uniquely determined up to scalar.  Let $z$ be a complex
variable.  Then $\mathrm{Res}_iA_i$ is defined by setting
\begin{equation}\label{eqn:basicresiduedefn}
\pi\in\mathfrak{S}_{i-1}\mapsto (\mathrm{Res}_iA_i)(\pi):=
(P_{\pi_0,\alpha_{b-i}^{\vee}}A_i)(\pi\otimes z\epsilon_i)|_{z=0}
\end{equation}
  It is well known (and follows
for example from the formula of Gindinkin-Karpelevic for the 
intertwining operators appearing in the $P_{a^b}$-constant term of $E^M(P_{a^b}))$)
that the collection singular hyperplanes passing through $\pi$ of the function
\[
\pi\in\mathfrak{S}_{i-1}\mapsto (P_{\pi_0,\alpha_{b-i}^{\vee}}A_i)(\pi\otimes z\epsilon_i)
\]
is the collection singular hyperplanes of $A_i$ passing through $\pi$ with the exception
of $\mathfrak{S}_{i-1}$.  That is to say, the positive power $n$ in the definition
in \S V.1.3 of \cite{mwbook} for $\mathrm{Res}_iA_i$ may be taken simply
to be $1$, and the linear combination of powers $P_{\pi_0,\alpha_i}^m$, $m\in[1,n]$,
called ``$Q_{n,\alpha_i^{\vee}}$'', may be taken simply to be $P_{\pi_0,\alpha_{b-i}^{\vee}}$,
As a result of these observations, the relatively simple definition \eqref{eqn:basicresiduedefn} suffices in our
case but not in general.

We can more explicitly calculate $\mathrm{Res}_{b-1}$, say, by choosing
\begin{equation}\label{eqn:epsilondefn}
\epsilon_{b-1}=\frac{1}{2}(1,-1,0,\cdots,0).
\end{equation}
Applying \eqref{eqn:hyperplanesdefn}, any $\pi\in\mathfrak{S}_{b-2}$ can be written in the form in the form
\[
\pi=\pi_0\otimes \mu,\;\text{where}\; \mu=(\mu_1,\mu_1,\mu_3,\ldots, \mu_b),\;\mu_i\in\mathbf{C}\; \text{such that}\;
2\mu_1+\sum_{i=3}^b\mu_i=0
\]
Then the mapping
\[
\mathrm{Res}_{b-1}E^M(P_{a^b}): \pi\mapsto \mathrm{Res}_{b-1}E^M(P_{a^b},\pi)\in
\mathrm{Hom}_{\mathbf{C}}(A(M_{a^b},\pi)^{\mathfrak{F}},L_{\xi,\mathrm{loc}}^{2,\mathfrak{F}})
\]
is given by
\[
\pi=\pi_0\otimes\mu\mapsto \left.\left(  P_{\pi_0,\alpha_1^{\vee}}(\pi_0\otimes\mu\otimes z\epsilon_i)E^{M}(P_{a^b},\pi_0\otimes\mu\otimes
z\epsilon_i) \right)\right|_{z=0}
\]
One readily computes that
\[
P_{\pi_0,\alpha_1^{\vee}}(\pi_0\otimes\mu\otimes z\epsilon_i)=\langle \mu\otimes z\epsilon_i,\alpha_1^{\vee} \rangle=z.
\]
Further, using \eqref{eqn:basepointdefn} and \eqref{eqn:epsilondefn}, we write $\pi\otimes z\epsilon_i$ in the form
\[
\pi\otimes z\epsilon_i=\pi^b\otimes \lambda_0(\pi)\otimes \lambda'(z),
\]
where
\[
\lambda_0(\pi):=\left(\frac{b-1}{2}+\mu_1,\frac{b-3}{2}+\mu_1,\ldots \frac{1-b}{2}+\mu_b\right),\;\text{and}\,
\lambda'(z)=\left(\frac{z}{2},-\frac{z}{2},0,\ldots,0\right).
\]
Then in terms of the variable $z$,
\[
\mathrm{Res}_{b-1}E^M(P_{a^b},\pi)=zE^M(P_{a^b},\pi^b\otimes(\lambda_0(\pi)+\lambda'(z)))|_{z=0}
\]
which is the residue, $\mathrm{Res}_{z=0}E^M(P_{a^b},\pi^b\otimes(\lambda_0(\pi)+\lambda'(z)))$, of a meromorphic 
complex function of one-variable at the origin.
Thus, the classical form of the residue theorem gives,
\begin{equation}\label{eqn:residuetointegral}
\mathrm{Res}_{b-1}E^M(P_{a^b},\pi)=\frac{1}{2\pi\mathbf{i}}\int_{\eta}
E^M(P_{a^b},\pi^{b}\otimes(\lambda_0(\pi)+\lambda'(z)))\,\mathrm{d}z,\;\text{for all}\; \pi\in\mathfrak{S}_{b-2},
\end{equation}
for $\eta$ a sufficiently small circle winding counterclockwise around the origin.

\begin{lemma}  \label{lem:residuesumcommutation}  Let $s\in\mathfrak{c}_{\Delta}$ be fixed. 
Then we have
\[
\sum_{\gamma\in P\backslash G} 
\mathrm{Res}^{\mathfrak{P}_{a^b}}_{\Lambda_b}E^M(P_{a^b})(\gamma g)m_{P}(\gamma g)^{\rho_{ab}+s}=
\mathrm{Res}^{\mathfrak{P}_{a^b}}_{\Lambda_b}E(P,E^M(P_{a^b}),g,s).
\]
\end{lemma}
\textbf{Proof.}  We are supposed to show that the operator 
$\mathrm{Res}^{\mathfrak{P}_{a^b}}_{\Lambda_b}$
commutes with the summation in \eqref{eqn:eisseriesdataeisseriesdefined}.
Clearly, it will suffice to show that each of the ``elementary residue" operators 
$\mathrm{Res}_{b-1},\,\mathrm{Res}_{b-2},\,\ldots, \mathrm{Res}_{1}$
commutes with the summation.

Let us show this for $\mathrm{Res}_{b-1}$.  For each $\gamma\in P\backslash G$,
define the operator of ``left-shift by $\gamma$", $\ell_{\gamma}$ as an operator
on any space of functions on $G$.  For example
\[
(\ell_{\gamma} E^M)(P_{a^b},\lambda\phi_{\pi^b},\lambda\pi^{b},g)=
E^M(P_{a^b},\lambda\phi_{\pi^b},\lambda\pi^{b},\gamma g).
\]
Then by \eqref{eqn:residuetointegral}, for all $\pi\in\mathfrak{S}_{b-2}$,
\begin{multline*}
\sum_{\gamma\in P\backslash G} 
((\mathrm{Res}_{b-1}E^M)(P_{a^b})(\pi))m_{P}^{\rho_{ab}+s}=\\
\sum_{\gamma\in P\backslash G} \int_{\eta}\ell_{\gamma}(E^M(P_{a^b})\cdot m_P^{\rho_{ab}+s})
(\pi^b\otimes (\lambda_0(\pi)+\lambda'(z)))\,\mathrm{d}z=\\
\int_{\eta}\sum_{\gamma\in P\backslash G}\ell_{\gamma}(E^M
(P_{a^b})
\cdot 
m_P^{\rho_{ab}+s})
(\pi^b\otimes (\lambda_0(\pi)+\lambda'(z)))\,\mathrm{d}z=\\
\mathrm{Res}_{b-1}E(P,E^M(P_{a^b},\pi),s).
\end{multline*}
The proof continues in a similar way until we obtain the commutation
statement of the Lemma.
\qed

Set
\begin{equation}\label{eqn:deltadefinition}
\Delta(\tau,b)(g):=\mathrm{Res}^{\mathfrak{P}_{a^b}}_{\Lambda_b}E^M(P_{a^b})(g).
\end{equation}
In \cite{jacquetarticle}, Jacquet used Langlands' criterion (the statement in \S I.4.11 of \cite{mwbook})
to show that $\Delta(\tau,b)$ is an square-integrable
automorphic representation on $\mathrm{GL}_{ab}$, or more precisely,
\[
\phi_{\pi^b}\in A_{\pi^b}\;\Rightarrow\; (\Delta(\tau,b))(\phi)\in L^2_{\rm disc}(\mathrm{GL}_{ab}(k)\backslash
\mathrm{GL}_{ab}(\mathbf{A})).
\]
Jacquet further conjectured, and Moeglin and Waldspurger
later proved, that all automorphic non-cuspidal discrete
spectrum representations of the general linear group are of this type.

\begin{thm}\cite{moeglinwaldspurgerens}
As $b$ ranges over the divisors of $n$, $n=ba$,
$\tau$ ranges over the irreducible, cuspidal, automorphic representations
of $\mathrm{GL}_a$, and $\phi_{\tau}$ ranges over  $A(N_{a^b}(\mathbf{A})M_{a^b}(k)
\backslash \mathrm{GL}_{ab}(\mathbf{A}))$, the automorphic forms $\Delta(\tau,b)(\phi_{\pi^b})$ span
$L^2_{\rm res}(\mathrm{GL}_{n}(k)\backslash\mathrm{GL}_n(\mathbf{A}))$.
\end{thm}
In this paper we do not specifically use this deep result, although it is important
for understanding the context, as explained in the introduction.

In order to state the following result, we will have to consider taking
a residue along a flag
\[
\mathcal{D}'=\mathcal{D}_{\Delta,\mathfrak{S}_0'}
\]
not to a point but to a one-dimensional affine
subspace $\mathfrak{S}_0'$ of $\mathfrak{P}$, which we now define.
Let $\mathfrak{S}_{b-1}'$ equal the whole space $\mathfrak{P}$.
Recall the simple system $\Delta(G,P_{a^b})$ of restricted roots 
introduced at \eqref{eqn:restrictedrootsexplicit}.  Using the notation
\[
\alpha_i: f_i-f_{i+1},\, \text{for}\; i=1,\ldots b-1,\; \alpha_b=2f_b,
\]
We define the flag by the same equations, formally as \eqref{eqn:basepointdefn}
and \eqref{eqn:hyperplanesdefn}, but now taking place in $\mathfrak{P}$,
instead of $\mathfrak{P}_{a^b}$.  Therefore, the intersection
is an affine space of dimension $1$, namely, setting $\mathbf{1}_b$
equal to the the $b$-vector having $1$ in all coordinates, and setting
$V$ equal to the one dimensional
vector space $s\mathbf{1}_a$
\[
\mathfrak{S}_0'= \pi_0\otimes V=\pi_0\otimes\left(\bigcap_{i=1}^b H_{\alpha_i^*}\right)=
\pi^b\otimes(\Lambda_b+\{s\mathbf{1}\;|\; s\in\mathbf{C}\})\subset \mathfrak{P}.
\]
is the $b$-tuple having $1$'s in all coordinates,
Then we can define the residue data $\mathrm{Res}^{\mathfrak{P}}_{\mathfrak{S}_0'}$.
By the general discussion of Chapter V of \cite{mwbook}, the image
under this operator of a meromorphic function on $\mathfrak{P}$
with polynomial singularities on $S_{X}^{\mathfrak{F}}$ is a meromorphic
function on $\mathfrak{S}_0'$, which can be identified with $\mathbf{C}$.
By Theorem \ref{eqn:meromorphiccontinuationgeneral},
the cuspidal-data Eisenstein series $E(P_{a^b})$ has a meromorphic continuation
to $\mathfrak{P}$. 
So in effect we have
\begin{equation}\label{eqn:meromorphicityresidualeisseries}
(\mathrm{Res}^{\mathfrak{P}}_{\mathfrak{S}_0'}E(P_{a^b},\pi^b))(s)\;\text{is a meromorphic function on}\; \mathbf{C}.
\end{equation}

From the results previously established above we can deduce
\begin{prop}\label{prop:residualdataeisseries}  Let $\Delta(\tau,b)(\phi_{\pi^b})$ be defined
as above.
\begin{itemize}
\item[(a)]  Let $s$ range over $\mathfrak{c}_{\mathfrak{P}_{a^b}}$.  Then the Eisenstein
series defined by
\[
E(P,\Delta(\tau,b)(\phi_{\pi^b}),g,s)
:=\sum_{\gamma\in P\backslash G}
(\Delta(\tau,b)(\phi_{\pi^b}))(\gamma g)m_P(\gamma g)^{s+\rho_{ab}},
\]
converges absolutely and uniformly for $g,s$ contained in compact subsets.
\item[(b)]  There is a unique meromorphic continuation in the complex variable $s$ of the Eisenstein
series $E(P,\Delta(\tau,b)(\phi_{\pi^b}),s)$
to the complement $U$ of a finite set of points in $\mathbf{C}$, and for all $s\in U$.
On $U$, this residual-data Eisenstein series can be identified with the residue of
the cuspidal-data Eisenstein series considered above, as follows,
 \begin{equation}\label{eqn:rescuspdatarel}
E(P,\Delta(\tau,b),s)=(\mathrm{Res}^{\mathfrak{P}}_{\mathfrak{S}_0'}E(P_{a^b},\pi^b))(s).
\end{equation}
\end{itemize}
\end{prop}
\textbf{Proof.}  By definition, the sum in part (a) is the sum in Lemma \ref{lem:residuesumcommutation}.
So we can apply Proposition \ref{prop:eisdataeisseries} to obtain part (a).

For part (b), first make the substitution of \eqref{eqn:eisensteinidentity} to obtain the equality
\eqref{eqn:rescuspdatarel} valid for all $s\in\mathfrak{c}_{\Delta}$.  
By \eqref{eqn:meromorphicityresidualeisseries}, we can define the meromorphic continuation of
$E(P,\Delta(\tau,b),s)$ as $\mathrm{Res}^{\mathfrak{P}}_{\mathfrak{S}_0'}$
applied to this meromorphically continued Eisenstein series.  By uniqueness,
this is the only possible meromorphic continuation.  Thus \eqref{eqn:eisensteinidentity}
is true for all but finitely many values of $s$, the poles of $(\mathrm{Res}^{\mathfrak{P}}_{\mathfrak{S}_0'}E(P_{a^b},\pi^b))(s)$.
\qed

\section{The Principal Non-vanishing Constant Term}
The present section contains the calculations that form the technical heart of
the paper.  Nevertheless, the reader may wish to skip this section
on a first reading and take the result as a generalization of (1.2.14)
of \cite{kudlarallisfest}.  At any rate, the result will be restated in the course
of deriving the main result in \S\ref{sec:residues}.

\subsection{Automorphic Forms and Induced Representations}
Our principle tool for describing the cuspidal support and exponents
of the automorphic forms defined in the introduction will be an ``inductive
formula" for the constant term
\[
\mathrm{CT}_Q E(P,\Delta(\tau,b),s)
\]
of the residual-data Eisenstein series along a the non-Siegel maximal parabolic
subgroup $Q=P_{a}^{ab}$.
In order to make this computation, it will be convenient to express the Eisenstein
series in the more ``classical" notation of induced representations.   For the following
subsection only, we follow certain notational conventions in \cite{soudryannals}.

Let $\Delta$ be $\Delta(\tau,b)$ as above. For the following discussion,
including Lemma \ref{lem:sectionofinduces}, $\Delta$ could be any square-integrable automorphic
representation of $\mathrm{GL}_n$.  Let
\[
\phi_{\Delta,s}\in\mathrm{Ind}_P^G(\Delta,s).
\]
This means that $\phi_{\Delta,s}: G\rightarrow A(\mathrm{GL}_n(k)\backslash\mathrm{GL}_n(A))_{\Delta}$
is a family of functions the complex variable $s$, varying holomorphically with $s$, each funcion satisfying a certain transformation law.  For $s\in\mathbf{C}$ fixed,
we may consider such a $\phi$ as a function of two variables, one in $G_n$ and the other in $\mathrm{GL}_n$,
and write the transformation law as
\begin{equation}\label{eqn:phideltastransflaw}
\phi_{\Delta,s}(umg;r)=|\mathrm{det} m|^{s+\rho_{ab}}\phi_{\Delta,s}(g;rm),\text{for all}\; r,m\in\mathrm{GL}_n(\mathbf{A}),\,
g\in G(\mathbf{A}), \; u\in U(\mathbf{A}).
\end{equation}
We set
\begin{equation}\label{eqn:datanewform}
f_{\Delta,s}^{\phi}(g)=\phi_{\Delta,s}(g;1),\text{for all}\; g\in G_n.
\end{equation}
where $1$ denotes the identity element in $\mathrm{GL}_n$.  The point
of the following lemma is to establish a precise equivalence
between between two
frameworks for constructing Eisenstein series,
with the `$f_{\Delta,s}^{\phi}$' representing the
framework of ``automorphic forms with respect to parabolics"
and the $\phi_{\Delta,s}$ representing the framework of
``holomorphic section of induced representations".

\begin{lemma}\label{lem:sectionofinduces}  We have
\[
f_{\Delta,s}^{\phi}\in A(U(\mathbf{A})M(k)\backslash G(\mathbf{A})_{\Delta\otimes s},
\]
and conversely, every element of $f\in A(U(\mathbf{A} M(k)\backslash G(\mathbf{A}))_{\Delta\otimes s}$
arises in this way, specifically as the $f^{\phi}_{\Delta,s}$ associated 
via \eqref{eqn:datanewform} to $\phi_{\Delta,s}$ defined by 
\begin{equation}\label{eqn:phiassoctof}
(g;r)\in G_n\times \mathrm{GL}_n\mapsto \phi_{\Delta,s}(g;r):= f(rg)|\det r|^{-s-\rho{ab}}.
\end{equation}
\end{lemma}
\textbf{Proof}  The proof consists of a few routine calculations, which we only sketch.
For the first statement, in order to show 
that $f_{\Delta,s}$ belongs to  $A(U(\mathbf{A} M(k)\backslash G)$,
we calculate that
\[
f_{\Delta,s}^{\phi}(g)=f(u(g)m(g)k(g))=\phi_{\Delta,s}(u(g)m(g)k(g);1)=\phi_{\Delta,s}(k(g);m(g))|\det m(g)|^{s+\rho_{ab}}
\]
In order to complete the proof of the first statement, by \S I.2.17 of \cite{mwbook},
we are to show that for each $k\in K$,
\begin{equation}\label{eqn:subkedveroff}
(f)_k(m):= |\det m|^{-\rho_{ab}}f_{\Delta,s}^{\phi}(mk)\in A(\mathrm{GL}_n(k)\backslash 
\mathrm{GL}_n(\mathbf{A}))_{\Delta\otimes s}.
\end{equation}
The above calculation implies that
\[
(f_{\Delta,s})_k(m)=\phi_{\Delta,s}(k;m)|\det m|^s.
\]
One completes the proof of \eqref{eqn:subkedveroff} by
using the part of the definition of $\Phi_{\Delta,s}$
that says that
$\phi_{\Delta,s}^{\phi}(k;\cdot)\in A(\mathrm{GL}_n(k)\backslash\mathrm{GL}_n(\mathbf{A} ))_{\Delta}$.

For the converse, let $f\in A(U(\mathbf{A}) M(k)\backslash G(\mathbf{A}))_{\Delta\otimes s}$
be given.  Define $\phi_{\Delta,s}$ by \eqref{eqn:phiassoctof}, and note this definition implies the relations
\[
\phi_{\Delta,s}(g;r_1)=|m(g)|^{\rho_{ab}}|\det r_1|^{-s}f_{\Delta,s, k(g)}(r_1m(g)).
\]
and
\[
\phi_{\Delta,s}(r_2 g; r_1)= |r_2 m(g)|^{\rho_{ab}}|\det r_1|^{-s} f_{\Delta,s,k(g)}(r_1r_2m(g)),
\]
for all $g\in G_n$ and $r_1,\,r_2\in \mathrm{GL}_{n}$.
Comparing the above two expressions, we obtain the transformation law,
\[
\phi_{\Delta,s}(r_2 g; r_1)=|\det r_2|^{s+\rho_{ab}}\phi_{\Delta,s}(g,r_2 r_2).
\]
We have shown that $\phi_{\Delta,s}$ belongs to $\mathrm{Ind}^G_P(\Delta,s)$, as required.\qed

\subsection{Non-normalized constant term}
\label{subsec:nonnormalizedconstterm}
Because of Lemma \ref{lem:sectionofinduces},
we have that the automorphic form $\Delta(\tau,b)(\phi_{\pi^b})$
may be identified as $f_{\Delta,s}^{\phi}$ for suitable
$\phi_{\Delta,s}\in \mathrm{Ind}^G_P(\Delta,s)$
as above.  kor the remainder of this section, we write $E(g,f_{\Delta,s}^{\phi})$
for $E(P,\Delta(\tau,b)(\phi),s)$.

We are calculating the $Q$-constant term of the following
Eisenstein series,
\[
E(g,f_{\Delta,s}^\phi)=\sum_{\gamma\in P\backslash G}
f^{\phi}_{\Delta,s}(\gamma g)=
\sum_{\gamma\in P\backslash G}\phi_{\Delta,s}(\gamma g;1).
\]
Initially, the integral for the constant term is
\[
\int\limits_{V(k)\backslash V(\mathbf{A})}E(vg,f^{\phi}_{\Delta,s})=
\int\limits_{V(k)\backslash V(\mathbf{A})}\sum_{\gamma\in P\backslash G}
f^{\phi}_{\Delta,s}(\gamma vg)\,\mathrm{d}v.
\]
Since the Bruhat double coset decomposition in this case is
\[
P\backslash G/V\cong\bigsqcup_{i=0}^a P\backslash P w_i L=
\bigsqcup_{i=0}^a w_i \; w_i^{-1}Pw_i\backslash w_i^{-1}Pw_i L=
\bigsqcup_{i=0}^a w_i\;  w_i^{-1}Pw_i\cap L\backslash L.
\]
Here $w_i$ is the Weyl element reversing the sign
of the first $i$ coordinates, with $w_0=\mathrm{Id}$.
We initially have, according to the standard unfolding,
\[
E_Q(g,f_{s,\Delta}^{\phi})=\sum_{i=0}^a\sum_{\gamma\in
 w_i^{-1}Pw_i\cap L\backslash L}
 \int_{V_{w_i}\backslash V_{w_i}(\mathbf{A})}
 \int\limits_{V^{w_i}(\mathbf{A})}f^{\phi}_{\Delta,s}
 (w_i\gamma v'v''g)\,\mathrm{d}v'\,\mathrm{d}v''.
\]

It is not difficult to verify the following.
\begin{lemma}  \label{lem:cuspidalityvanishing} The term corresponding to $w_i$ for $i=1,\ldots, a-1$
vanishes by cuspidality of $\tau$.
\end{lemma}

So we have two terms, first
\[
E_Q(g,f_{s,\Delta}^{\phi})_{\rm Id}=
\sum_{\gamma\in P\cap L\backslash L}
\int_{V(k)\backslash V(\mathbf{A})}
f^{\phi}_{\Delta,s}(\gamma v'g)\,\mathrm{d}v',
\]
where $v'$ ranges over
\begin{equation}\label{eqn:vprimenotation}
v'(X,Y,Z)=\begin{pmatrix}1_a&X&Y&Z\\
0&1_{a(b-1)} &0 &Y'\\
0 & 0&  1_{a(b-1)} & X'\\
0 &0 &0& 1_a\end{pmatrix},\;\text{and}\; X'=-j{}^t X j,\,
Y'=j{}^t Y j,
\end{equation}
and $Y$ and $X$ range over 
$k^{a\times ab-a}\backslash\mathbf{A}^{a\times ab-a}$,
and $Z$ over $k^{\frac{a(a+1)}{2}}\backslash\mathbf{A}^{\frac{a(a+1)}{2}}$, with the symmetry of $Z$ occurring around the second diagonal.
Second,
\begin{multline*}
E_Q(g,f_{s,\Delta}^{\phi})_{w}=
\sum_{\gamma\in w^{-1}Pw\cap L\backslash L}
 \hspace*{.3cm}\int\limits_{V_w(\mathbf{k})\backslash V_w(\mathbf{A})}
\int\limits_{V^w(\mathbf{A})} f^{\phi}_{\Delta,s}
(w\gamma v' v'' g)\,\mathrm{d}v'\,\mathrm{d}v''=\\
\sum_{\{1_a                                                                                                                                                                                                                                                                                                                                                                                                                                                                                                                                                                                                                                                                                                                                                                                                                                                                                                                                                                                                                                                                                                                                                        \}\times P_{ab-a}\backslash G_{ab-a}}
\int_{Y\in k^{a^2(b-1)}\backslash \mathbf{A}^{a^2(b-1)}}
\int\limits_{\begin{matrix}{\scriptstyle X\in \mathbf{A}^{a^2(b-1)}}\\
{\scriptstyle Z\in \mathbf{A}^{\frac{a(a+1)}{2}}}
\end{matrix}}
f^{\phi}_{\Delta,s}
(w\gamma v'(Y) v''(X,Z) g)\,\mathrm{d}Y\,\mathrm{d}X\,\mathrm{d}Z
\end{multline*}
In these two terms, the `$\gamma$' commutes
elementwise with the factor $v'(0,0,Z)$ (respectively, $v''(0,Z)$).
It does not commute elementwise, however, with the other factor $v'(X,Y,0)$ (resp., the factors $v'(Y)$, $v''(X,0)$).  Because
conjugation by $\gamma$ is a unimodular, rational transformation,
using an appropriate change of variables, we see that $\gamma$
does commute with the entire integral over $X,Y$, after appropriate
change of variable.

Further, the integration of $Y,Z$ in the identity term clearly
contributes a constant factor, which in the
usual normalization of measures is $1$.  So we have
\begin{multline*}
E_Q(g,f_{s,\Delta}^{\phi})=
E_Q(g,f_{s,\Delta}^{\phi})_{\rm Id}+E_Q(g,f_{s,\Delta}^{\phi})_{w}=\\
\sum_{\gamma\in 1_a\times P_{n-a}\backslash G_{n-a}}
\int\limits_{X\in (k\backslash\mathbf{A})^{a^2(b-1)}
}f^{\phi}_{\Delta,s}(v'(X,0,0)\gamma g)\,\mathrm{d}v'+\\
\sum_{\gamma\in 1_a\times P_{n-a}\backslash G_{n-a}}
\hspace*{.3cm}
\int\limits_{Y\in k^{a^2(b-1)}\backslash \mathbf{A}^{a^2(b-1)}}
\int\limits_{\begin{matrix}{\scriptstyle X\in \mathbf{A}^{a^2(b-1)}}\\
{\scriptstyle Z\in \mathbf{A}^{\frac{a(a+1)}{2}}}
\end{matrix}}f^{\phi}_{\Delta,s}(wv'(Y)v''(X,Z)\gamma g)\,\mathrm{d}X
\,\mathrm{d}Y\,\mathrm{d}Z.
\end{multline*} 
After we restrict to a $g\in G$ of the form $g=\mathrm{diag}(t,h,\tilde{t})\in L$, we claim that each term is of the form $|\det t|$
to an appropriate exponent, times an Eisenstein series
on $G_{n-a}$ on a section of the induced representation
$\mathrm{Ind}_{P_{n-a}}^{G_{n-a}}(\Delta(\tau,b-1),s')$.
We now prove this, in the process determining the
exponent of $|\det t|$ and the value of the parameter
$s'$ precisely.

\vspace*{0.3cm}\noindent
\textbf{Identity Term.}  Fix a number $q$, to be specified later.  We \textit{define}
\begin{equation}\label{eqn:identitytermsummand}
\psi_{\Delta}(t;h;\ell)=|\det t|^{-q}\int\limits_{X\in k^{a^2(b-1)}\backslash \mathbf{A}^{a^2(b-1)}}
\phi_{\Delta,s}(v'(X,0,0)\mathrm{diag}(t,h,\tilde{h},\tilde{t});\mathrm{diag}(1_a,\ell))\,\mathrm{d}x,
\end{equation}
where $t\in\mathrm{GL}_{a}$, $h\in G_{a(b-1)}$, and $\ell\in\mathrm{GL}_{a(b-1)}$.
We sometimes write $\boldell$ for $\mathrm{diag}(1_a,\ell)$.

Since
\[
f^{\phi}_{\Delta,s}(g)=\varphi^{\phi}_{s,\Delta}(g;1),\; \text{where}\;
g\in G,\, 1\in \mathrm{GL}_n,
\]
we have
\[
E_Q(\mathrm{diag}(t,h,\tilde{t}),f^\phi_{\Delta,s})_{\rm Id}=|\det t|^q
\sum_{\gamma\in P_{n-a}\backslash G_{n-a}}\psi_{\Delta}(t;\gamma h;1)
\]
We claim that for appropriate choice of $q$, one has
\begin{equation}
\psi_{\Delta}\in \tau\otimes \mathrm{Ind}_{P_{n-a}}^{G_{n-a}}\Delta(\tau,b-1)\cdot |\det \cdot|^{s+\frac{1}{2}}.
\end{equation}
In order to verify the claim, we have to compute
\[
\psi_{\Delta}(t;\mathrm{diag}(r,\tilde{r})h;\boldell),\;\text{for}\; r\in\mathrm{GL}_{a(b-1)}.
\]
Note that
\begin{equation}\label{eqn:productrearrangement}
\mathrm{diag}(1_a,r,\tilde{r},1_a)\mathrm{diag}(t,h,\tilde{t})=\mathrm{diag}(t,r,\tilde{r},\tilde{t})
\mathrm{diag}(1_a,h,1_a).
\end{equation}

Also, note that
\[
v'(X,0,0)=\mathrm{diag}\left(n_{a,ab-a}(X),\widetilde{{n}_{a,ab-a}(X)}\right).
\]
Thus, also using the transformation law
\eqref{eqn:phideltastransflaw},
\eqref{eqn:identitytermsummand} can be rewritten
\begin{multline*}
\int\limits_{X\in k^{a^2(b-1)}\backslash  \mathbf{A}^{a^2(b-1)}}
{\varphi}^{\phi}_{\Delta,s}(\mathrm{diag}\left(n_{2,2}(X),\widetilde{{n}_{2,2}(X)}\right)
\mathrm{diag}(t,r,\tilde{r},\tilde{t})
\mathrm{diag}(1_a,h,1_a);\ell)\,
\mathrm{d} X=
\\ |\det t|^{s+\rho_{ab}}|\det r|^{s+\rho_{ab}}
\int\limits_{n\in N_{a,ab-a}(\mathbf{A})}
{\varphi}^{\phi}_{\Delta,s}(\mathrm{diag}(1_a,h,1_a);\boldell
n\,
\mathrm{diag}(t,r))\,
\mathrm{d} X.
\end{multline*}
Now we use the description of $\phi_{\Delta}$ from 
\cite{moeglinwaldspurgerens}, saying that for each $g\in G$, 
\[
\phi_{\Delta}(g)=\mathrm{Res}_{\mathbf{s}=\Lambda_b}
E(f^{\phi_g}_{\pi,\mathbf{s}})=\mathrm{Eval}\prod_{i=1}^{b-1}
(s_{i+1}-s_{i}+1)\sum_{w\in W}M(w,\mathbf{s})
f^{\phi_g}_{\pi,\mathbf{s}},
\]
where the sum is over $W=W(\mathrm{GL}_{ab},P_{a^b})$.
Recall the normalization of the intertwining operator
\[
M(w,\mathbf{s})=r(w,\mathbf{s})R(w,\mathbf{s}),
\]
where $r(w,\mathbf{s})$ is a certain ratio of L-functions
determined by the formula of Gindinkin-Karpelevich, and meromorphic, and $R(w,\mathbf{s})$ is a holomorphic operator.
It is not difficult to see that, at $\mathbf{s}=\Lambda_b$,
the only summand which has a singularity of total order $b-1$ is
the term $r(w_{\sigma},\mathbf{s})$ corresponding to the longest Weyl element $w_{\sigma}$.
Here $\sigma$ is the permutation reversing the order of the
$b$ coordinates.  So we obtain that, up to a constant
(namely the residue of $r(w_{\sigma},\mathbf{s})$ at the point
$\mathbf{s}=\Lambda_b$),
\begin{equation}\label{eqn:spehasresidue}
\phi_{\Delta}(g;\cdot)=\mathrm{Eval}_{s=\Lambda_b}f^{\phi_g}_{\pi,
\mathbf{s}}(\cdot)=\mathrm{Eval}_{s=\Lambda_b}\phi^{g}_{\pi,\sigma
\mathbf{s}}(\cdot;\mathbf{1}^b),
\end{equation}
where $\mathbf{1}^b$ is the identity element in $(\mathrm{GL}_a)^b$.
So we have
\begin{multline*}
\psi_{\Delta}(t;r;h)=|\det t|^{s+\rho_{ab}-q}|\det r|^{s+\rho_{ab}}
|\det m(\mathbf{h})|^s\\
\mathrm{Eval}_{s=\Lambda_b}\int\limits_{n\in N_{a,ab-a}(k)\backslash N_{a,ab-a}(\mathbf{A})}
\phi_{\pi,\sigma\mathbf{s}}^{\mathbf{h}}(\ell n\,\mathrm{diag}(t,r);\mathbf{1}^b),
\end{multline*}
where we have used the abbreviation $\mathbf{h}=\mathrm{diag}(1_a,h,1_a)$.

Since $N_{a,ab-a}\subseteq N_{a^b}$ and $\phi^{\mathbf{h}}_{\pi,\sigma \mathbf{s}}(g;h)$ is invariant by 
$N_{a^b}$ in the first variable, the entire integral evaluates to $1$ because of the normalization
of the measures.  Since $\boldell=\mathrm{diag}(1_a,\ell)$ we have
the commutation relation
\[
\boldell(t,1_a,\ldots,1_a)=(t,1_a,\ldots,1_a)\boldell
\]
Therefore,
\begin{multline*}
\psi_{\Delta}(t;\mathrm{diag}(r,\tilde{r})h;\ell)=|\det r|^{s+\rho_{ab}}|\det t|^{s+\rho_{ab}-q}
|\det m(\mathbf{h})|^s\\
\mathrm{Eval}_{\mathbf{s}=\Lambda_b}\phi^{\mathbf{h}}_{\pi,\sigma\mathbf{s}}
(\mbox{\boldmath$($\unboldmath}t,1_a,\ldots,1_a\boldcp\boldell\boldop1_a,r\boldcp;1),
\end{multline*}
where we have eliminated the ``diag" to improve readability, and replace it with bold parentheses.
Here, we have
\[
\phi^{\mathbf{h}}_{\pi,\sigma\mathbf{s}}\in \mathrm{Ind}_{P_{a^b}}^{G_{ab}}(\pi, \sigma\mathbf{s}).
\]
meaning that, first, as a function of the first variable, in $G_{ab}$, $\phi^{\mathbf{h}}_{\pi,\sigma\mathbf{s}}$
 takes values in the space of $\pi^b$,
and, second, it satisfies the ``transformation law"
\[
\phi_{\pi,\sigma\mathbf{s}}^{\mathbf{h}}(\mathbf{m}g;\mathbf{r})=
|\det \mathbf{m}|^{\sigma\mathbf{s}}\delta^{\frac{1}{2}}_{a^b}
(\mathbf{m})\phi^{\mathbf{h}}_{\pi,\sigma\mathbf{s}}(g;\mathbf{rm}),
\]
where $\mathbf{m},\, \mathbf{r}\in M_{a^b}\cong (\mathrm{GL}_a)^b$.

Therefore,
\begin{multline*}
\phi^{\mathbf{h}}_{\pi,\sigma\mathbf{s}}
(\boldop t,1_a,\ldots,1_a\boldcp\boldop1_a,r\boldcp;1)=\\
|\boldop t,1_a,\ldots,1_a\boldcp|^{\sigma\mathbf{s}}
\delta^{\frac{1}{2}}_{a^b}(\boldop t,1_a,\ldots,1_a\boldcp)
\phi^{\mathbf{h}}_{\pi,\sigma\mathbf{s}}(\boldop1_a,r\boldcp;\boldop t,1_a,\ldots,1_a\boldcp)=\\
|\det t|^{\mathbf{s}_b+(\rho_{a^b})_1}\phi^{\mathbf{h}}_{\pi,\sigma\mathbf{s}}
(\boldop1_a,r\boldcp;\boldop t,1_a,\ldots,1_a\boldcp).
\end{multline*}
Thus, we have
\begin{multline}\label{eqn:psideltainterform}
\psi_{\Delta}(t;\boldop r,\tilde{r}\boldcp h;\ell)=|\det r|^{s}\delta_n^{\frac{1}{2}}(r)
|\det t|^{s+\rho_{ab}+(\Lambda_b)_b+(\rho_{P_{a^b}})_1-q}
|\det m(\mathbf{h})|^s\times \\ \times\mathrm{Eval}_{\mathbf{s}=\Lambda_b}\phi_{\pi,\sigma\mathbf{s}}^{\mathbf{h}}
(\boldell\boldop 1_a,r\boldcp;\boldop t,1_a,\ldots,1_a\boldcp).
\end{multline}
Now, given any $\phi_{\pi,\mathbf{s}}$ holomorphic section of 
$\mathrm{Ind}_{P_{a^b}}^{G_{ab}}
(\pi,\mathbf{s})$, and scalar $u$, define
\[
\psi_{\pi}(t;r;\mathbf{m}):=|\det r|^u\phi_{\pi,\mathbf{s}}(\boldop 1_a,r\boldcp;\boldop t,\mathbf{m}\boldcp),
\]
for all $t\in\mathrm{GL}_a$, $\mathbf{m}\in(\mathrm{GL}_a)^{b-1}$, $r\in(\mathrm{GL}_{a(b-1)})$.
We claim that
\begin{equation}\label{eqn:psipiclaim}
\psi_{\pi}(t;r;\mathbf{m})\in\tau\otimes\mathrm{Ind}_{P_{a^{b-1}}}^{\mathrm{GL}_{a(b-1)}}(\pi^{b-1},\mathbf{v}'),\;
\text{where}\; \mathbf{v}'=\mathbf{s}'+\mathbf{u}'-\frac{1}{2}\mathbf{a}'.
\end{equation}
Here, we have used the following notation: for a scalar such as $u$, $\mathbf{u}$ denotes the
$b$-tuple $(u,\ldots, u)$ and for any $b$-tuple such as $\mathbf{s}$, $\mathbf{s}$ 
denotes the truncation $(s_2,\ldots s_b)$ to a $b-1$-tuple obtained by deleting the first entry.
The claim means that for any $r\in\mathrm{GL}_{a(b-1)}$, $\psi_{\pi}(r)$ takes values
in the space of $\tau\otimes\pi^{b-1}$ and satisfies the transformation law
\begin{multline*}
\psi_{\pi}(t;\mathbf{m}_2r;\mathbf{m}_1)=\psi_{\pi}(t;r;\mathbf{m}_1\mathbf{m}_2)\delta_{P_a^{b-1}}^{\frac{1}{2}}
(\mathbf{m}_2)|\det \mathbf{m}_2|^{\mathbf{v}'},\\ \text{for all}\; r\in\mathrm{GL}_{a(b-1)},\, t\in\mathrm{GL}_a,\,
\mathbf{m}_1,\mathbf{m}_2\in(\mathrm{GL}_a)^{b-1}.
\end{multline*}
In order to verify the claim, first note the following,
\begin{multline*}
\psi_{\pi}(t;\mathbf{m}_2r;\mathbf{m}_1)=\phi_{\pi,\mathbf{s}}(\boldop1_a,m_2\boldcp\boldop1_a,r\boldcp;
(t,\mathbf{m}_1))|\det r|^{u}|\det \mathbf{m}_2|^{\mathbf{u}'}=\\
=\phi_{\pi,\mathbf{s}}(\boldop1_a,r\boldcp ;\boldop t,\mathbf{m}_1 \mathbf{m}_2\boldcp)|\det r|^u
|\det \mathbf{m}_2|^{\mathbf{u}'}|\det \boldop 1,\mathbf{m}_2\boldcp 
|^{\mathbf{s}}\delta_{a^b}^{\frac{1}{2}}(\boldop1,m_2\boldcp)=\\
\psi_{\pi}(t;r;\mathbf{m}_1\mathbf{m}_2)|\det \mathbf{m}_2|^{\mathbf{s}'+\mathbf{u}'}\delta_{a^b}^{\frac{1}{2}}
(\boldop1,\mathbf{m}_2\boldcp).
\end{multline*}
Since
\[
\delta^{\frac{1}{2}}_{a^b}(1,\mathbf{m})=
(\det\boldop1,\mathbf{m}\boldcp)^{a\left(\frac{b-1}{2},\cdots,\frac{1-b}{2} \right)}=
|\det \mathbf{m}|^{a\left(\frac{b-3}{2},\cdots,\frac{1-b}{2} \right)},
\]
we see that
\[
\delta^{\frac{1}{2}}_{a^b}(1,\mathbf{m}_2)\cdot \delta^{-\frac{1}{2}}_{a^{b-1}}(\mathbf{m}_2)=
|\det \mathbf{m}_2|^{-\frac{1}{2}\mathbf{a}'}.
\]
Thus,
\begin{multline*}
\psi_{\pi}(\mathbf{m}_2r;\mathbf{m}_1;t)=|\det\mathbf{m}_2|^{\mathbf{s}'+\mathbf{u}'-\frac{1}{2}\mathbf{a}'}
 \delta^{\frac{1}{2}}_{a^{b-1}}(\mathbf{m}_2)\psi_{\pi}(t;r;\mathbf{m}_1\mathbf{m}_2)|=\\
 |\det\mathbf{m}_2|^{\mathbf{v}'}
 \delta^{\frac{1}{2}}_{a^{b-1}}(\mathbf{m}_2)\psi_{\pi}(t;r;\mathbf{m}_1\mathbf{m}_2)|,
\end{multline*}
as claimed.
Further, note that
\[
\mathbf{s}=\sigma\Lambda_b,\; \mathbf{v}'=\sigma'\Lambda_{b-1},\;\text{implies}\; \mathbf{u}'=\frac{1}{2}(a-1,\ldots,
a-1).
\]

That is to say, with $u=\frac{1}{2}(a-1)$,
\[
\phi_{\pi}\in \mathrm{Ind}^{G_{ab}}_{P_{a^b}}(\pi,\sigma\mathbf{s})|_{s=\Lambda_b}
\;\mathrm{implies}\; \psi_{\pi}\in \tau\otimes\mathrm{Ind}_{P_{a^{b-1}}}^{\mathrm{GL}_{n-a}}\pi^{b-1},
\sigma\mathbf{v}')|_{\mathbf{v}'=\Lambda_{b-1}}.
\]
In particular if we take $\phi_{\pi}=|\det m(h)|^{-\frac{1}{2}}\phi^{\mathbf{h}}_{\pi,\sigma\mathbf{s}}|_{s=\Lambda_b}$,
then we have
\begin{multline*}
\mathrm{Eval}_{\mathbf{s}=\Lambda_b}|\det r|^{\frac{1}{2}(a-1)}
\phi_{\pi,\sigma\mathbf{s}}(\boldop 1_a,r\boldcp;\boldop t, 1_a,\cdots, 1_a\boldcp)=
\mathrm{Eval}_{\mathbf{s}=\Lambda_b}\psi_{\pi}(t;r;\mathbf{1}^{b-1})=\\
\phi_{\tau}(t)\otimes \mathrm{Eval}_{\mathbf{v}'=\Lambda_{b-1}'}\phi_{\pi^{b-1},\sigma'\mathbf{v}'}^h
(r;\mathbf{1}^{b-1}).
\end{multline*}
Thus, \eqref{eqn:psideltainterform} implies that
\begin{multline*}
\psi_{\Delta}(t;\boldop r,\tilde{r}\boldcp h;\boldell)=|\det r|^{s+\frac{1}{2}(1-a)+\rho_{ab}}
|\det t|^{s+\rho_{ab}+(\Lambda_b)_b+(\rho_{P_{a^b}})_1-q}|\det m(h)|^{s+\frac{1}{2}}\times\\
\phi_\tau(t)\otimes\mathrm{Eval}_{\mathbf{v}'=\Lambda_{b-1}}\phi_{\pi^{b-1,\sigma'\mathbf{v}'}}^{\mathbf{h}}
(\ell r;\mathbf{1}_{b-1}).
\end{multline*}
Set
\[
q=s+\rho_{ab}+(\Lambda_b)_b+(\rho_{P_{a^b}})_1=s+ab+1-\frac{1}{2}(b+a),
\]
and use $|\det r|^{-\frac{1}{2}a+\rho_{ab}}=
|\det r|^{\rho_{a(b-1)}}$ to see that
\[
\psi_{\Delta}(t;\boldop r,\tilde{r}\boldcp h;\boldell)=|\det r|^{s+\frac{1}{2}}\delta_{n-a}^{\frac{1}{2}}(r)
|\det m(h)|^{s+\frac{1}{2}}\phi_\tau(t)\otimes\phi_{\Delta^{b-1}}(h,\ell r).
\]
If we repeat the calculation setting $r=1$ and $\ell=\ell r$, we obtain
\[
\psi_{\Delta}(t;h;\ell r)=|\det m(h)|^{s+\frac{1}{2}}\phi_\tau(t)\otimes\phi_{\Delta^{b-1}}(h,\ell r)
\]
We have therefore shown that $\psi_{\Delta}$ takes values in $\tau\otimes\mathrm{Ind}\Delta^{b-1}$
and satisfies the transformation law
\[
\psi_{\Delta}(t;\boldop r,\tilde{r}\boldcp h;\boldell)=|\det r|^{s+\frac{1}{2}}\delta_{a(b-1)}^{\frac{1}{2}}(r)\psi_{\Delta}(t;h;\ell r).
\]
So we conclude that
\begin{multline}\label{eqn:identitytermcomputed}
E_Q^{ab}(f_{\Delta,s},\mathrm{diag}(t,h,\tilde{h},\tilde{t}))_{\rm Id}=
|\det t|^{s+ab-\frac{1}{2}(b+a)+1}\phi_{\tau}(t)\otimes
E^{a(b-1)}(h,f'_{\Delta,s+\frac{1}{2}}),\\ \text{where}\; f'_{\Delta,s+\frac{1}{2}}:=
i\circ f_{\Delta,s},
\end{multline}
and $i$ denotes the inclusion map of $G_{a(b-1)}$ into $G_{ab}$ given by $h\mapsto \mathrm{diag}(1_a,h,1_a)$.

\vspace*{0.3cm}\noindent
\textbf{$w$-term.}  Following a similar pattern as for the identity term, we set
\begin{multline*}
\psi_{\Delta}^w(t;h;\ell)=|\det \ell|^{-q_1}|\det t|^{-q_2}\\
\int\limits_{Y\in  (k\backslash \mathbf{A})^{a^2(b-1)}}
\int\limits_{\begin{matrix}{\scriptstyle X\in \mathbf{A}^{a^2(b-1)}}\\
{\scriptstyle Z\in \mathbf{A}^{\frac{a(a+1)}{2}}}
\end{matrix}}\phi_{\Delta,s}(wv'(Y)v''(X,Z)\boldop t,h,\tilde{t}
\boldcp;\boldop1_a,\ell\boldcp)\,\mathrm{d}X\,
\mathrm{d}Y\,\mathrm{d}Z,
\end{multline*}
where $t\in\mathrm{GL}_a$, $h\in\mathrm{GL}_{a(b-1)}$ and $\ell\in\mathrm{GL}_{a(b-1)}$.

Since $f^{\phi}_{\Delta,s}(g)=\phi_{\Delta,s}(g;1)$ for all $g\in G$ and $1$ the identity in $\mathrm{GL}_n$,
we have
\begin{equation}\label{eqn:wtermabstractform}
E(\mathrm{diag}(t,h,\tilde{t}),f_{\Delta,s}^{\phi})_w=|\det t|^{q_2}
\sum_{\gamma\in P_{n-a}\backslash G_{n-a}}\psi_{\Delta}(t;\gamma h;1).
\end{equation}
We claim that for appropriate choice of $q_1$, $q_2$, one has
\[
\psi_{\Delta}^w\in \mathrm{Ind}^{G_{n-a}}_{P_{n-a}}(\Delta(\tau,b-1),s-\frac{1}{2})
\otimes\tau.
\]
In order to verify this, we have to compute
\[
\psi_{\Delta}^w(t;\mathrm{diag}(r,\tilde{r})h;\ell),\;\text{for}\; r\in \mathrm{GL}_{a(b-1)}.
\]
By \eqref{eqn:productrearrangement}, we have
\[
\boldop1_a,r,\tilde{r},1_a\boldcp\boldop t,h,\tilde{t}\boldcp=
\boldop t,r,\tilde{r},\tilde{t}\boldcp \mathbf{h}\;\text{where}\; h=\boldop1_a,h,1_a \boldcp.
\]
Note that
\[
wv'v''\boldop t,r,\tilde{r},\tilde{t}\boldcp\mathbf{h}= wv'w^{-1}\cdot w\boldop t,r,\tilde{r},\tilde{t}\boldcp w^{-1}\cdot w\cdot
\boldop t,r,\tilde{r},\tilde{t}\boldcp^{-1} v'' \boldop t,r,\tilde{r},\tilde{t}\boldcp\cdot \mathbf{h}.
\]
Now the integral over $(k\backslash \mathbf{A})^{a^2(b-1)}$
becomes an integral over $N_{a,ab-a}^-(k)\backslash
N_{a,ab-a}^-(\mathbf{A})$ of $\boldop n,\tilde{n}\boldcp$.  
We use the abbreviation ${N'}^{,-}$ for $N_{a,ab-a}^-$.  

Further,
one easily computes that
\[
w\boldop t,r,\tilde{r},\tilde{t}\boldcp w^{-1}=\boldop {}^t t^{-1},r,\tilde{r},jtj \boldcp.
\]
Then by \eqref{eqn:phideltastransflaw}, and using the bold letter $\boldell$ 
as abbreviation for $(1_a,\ell)$, we have
\begin{multline*}
\psi_{\Delta}^w(t;\boldop r,\tilde{r}\boldcp h;\ell)=
|\det \ell|^{-q_1}|\det t|^{-q_2}\\ \int\limits_{n\in {N'}^{,-}\backslash {N'}^{,-}(\mathbf{A})}
\int\limits_{\begin{matrix}{\scriptstyle X\in \mathbf{A}^{a^2(b-1)}}\\
{\scriptstyle Z\in \mathbf{A}^{\frac{a(a+1)}{2}}}
\end{matrix}}\phi_{\Delta,s}(\boldop {}^t t^{-1},r,\tilde{r},jtj \boldcp\cdot
w\cdot \boldop t,r,\tilde{r},\tilde{t}\boldcp^{-1}v''(X,Z)\boldop t,r,\tilde{r},\tilde{t}\boldcp\cdot \mathbf{h};\boldell n)
\mathrm{d}n\,
\mathrm{d}X\,\mathrm{d}Z.
\end{multline*}
Note that
\begin{equation}\label{eqn:conjugatingvdoubleprime}
\boldop t,r,\tilde{r},\tilde{t}\boldcp^{-1}v''(X,Z)\boldop t,r,\tilde{r},\tilde{t}\boldcp=
v''(t^{-1}X r,t^{-1}Z\tilde{t}).
\end{equation}
Since $X$ ranges over arbitrary $\mathbf{A}$-valued
 matrices, of dimension, $a$ by $a(b-1)$, the Jacobian of the transformation $X\mapsto tXr^{-1}$
is $|\det t|^{a(b-1)}|\det r|^{-a}$.  Further, using the theory of root systems or invariant theory
one sees that for $Z$ in the set of integration, the Jacobian of $Z\mapsto tZ\tilde{t}^{-1}$ is
$|\det t|^{a+1}$.  Combining these observations with another application of \eqref{eqn:phideltastransflaw},
we have
\begin{multline}\label{eqn:psiwdeltainterform}
\psi_{\Delta}^w(t;\boldop r,\tilde{r}\boldcp h;\ell)=
\delta_n^{\frac{1}{2}}(\boldop {}^t t^{-1},r  \boldcp)|\det \ell|^{-q_1}
|\det r|^{s-a}|\det t|^{-q_2-s+ab+1}\\
\int\limits_{n\in {N'}^{,-}\backslash {N'}^{,-}(\mathbf{A})}
\int\limits_{\begin{matrix}{\scriptstyle X\in \mathbf{A}^{a^2(b-1)}}\\
{\scriptstyle Z\in \mathbf{A}^{\frac{a(a+1)}{2}}}
\end{matrix}}\phi_{\Delta,s}(
w\cdot v''(X,Z)\cdot \mathbf{h};\boldell n\boldop{}^t t^{-1},r\boldcp)
\mathrm{d}n\,
\mathrm{d}X\,\mathrm{d}Z.
\end{multline}
Now, make the abbreviation
\[
\mathrm{Factor} 1:=\delta_n^{\frac{1}{2}}(\boldop {}^t t^{-1},r
\boldcp)|\det r|^{s-a}|\det t|^{-q_2-s+ab+1}.
\]
We temporarily abbreviate as
\[
\text{$v$, the element $w\cdot v''(X,Z)\cdot \mathbf{h}$}.
\]
we denote by
\[
\text{$V^h$ the set over which, for fixed $h$,
$v$ ranges}
\] 
as $X$ and $Z$ range over their respective
varieties.
Then we have
\begin{multline*}
\psi^w_{\Delta}(t;\boldop r,\tilde{r}\boldcp h;\ell)=
\mathrm{Factor} 1\times |\det \ell|^{-q_1}
\int\limits_{n\in {N'}^{,-}\backslash {N'}^{,-}(\mathbf{A})}
\int\limits_{v\in V^h}\phi_{\Delta,s}(v;\boldell n\boldop {}^t t^{-1},r\boldcp)
\mathrm{d}n\,\mathrm{d}v.
\end{multline*}
Since $\phi_{\Delta,s}(v,\cdot)$ is left invariant under $k$-points, we may multiply on the left
by a $w'\in W$, represented by an element of $\mathrm{GL}_n(k)$.  Further, we easily see that
\[
w' \boldell n\boldop{}^t t^{-1},r\boldcp = w'\boldell {w'}^{-1}\cdot w'
 n {w'}^{-1}\cdot w'\boldop{}^t t^{-1},r\boldcp {w'}^{-1}\cdot w'
\]
Now we take the Weyl element to be 
$w'\in W(\mathrm{GL}_n, P_{a^b})\in\mathfrak{S}_b$,
represented by
\begin{equation}\label{eqn:wsubsigmadefn}
w'=w_{\sigma}=\begin{pmatrix}\\ &1_{ab-a}\\ 1_a&\end{pmatrix}.
\end{equation}
Note that $w'=w_{\sigma}$ for $\sigma$ the cyclic permutation
\[
\sigma=(1b(b-1)\cdots 2)\in \mathfrak{S}_b
\]
 Then we calculate that
\[
w_{\sigma}\boldell {w_{\sigma}}^{-1}=\boldop \ell,1_a\boldcp;
\quad w_{\sigma}\boldop {}^t t^{-1},r \boldcp {w_{\sigma}}^{-1}=
\boldop 1_{ab-a}, {}^t t^{-1} \boldcp \boldop r,1_a\boldcp;
\quad w_{\sigma} N^-_{a,ab-a}{w_{\sigma}}^{-1}=N_{ab-a,a}.
\]
So by \eqref{eqn:spehasresidue}
\begin{multline*}
\psi^w_{\Delta}(t;\boldop r,\tilde{r}\boldcp h;\ell)=
\mathrm{Factor}1\times |\det \ell|^{-q_1}\\
\int\limits_{n\in {N'}^{,-}\backslash {N'}^{,-}(\mathbf{A})}
\int\limits_{v\in V^h}|\det m(v)|^s
\mathrm{Eval}_{\mathbf{s}=\Lambda_b}\phi^v_{\pi,\sigma\mathbf{s}}
(w_{\sigma}\ell n\boldop{}^t t^{-1},r\boldcp);\mathbf{1})\,\mathrm{d}v\,\mathrm{d}n=\\
\mathrm{Factor}1\times |\det \ell|^{-q_1}\hspace*{-.3cm}
\int\limits_{n\in N_{ab-a,a}\backslash N_{ab-a,a}(\mathbf{A})}
\int\limits_{v\in V^h}|\det m(v)|^s\\
\mathrm{Eval}_{\mathbf{s}=\Lambda_b}
\phi_{\pi,\sigma\mathbf{s}}^v
(\boldop \ell,1_a\boldcp n\boldop 1_{ab-a,a},{}^tt^{-1}\boldcp
\boldop r,1_a\boldcp\cdot w_{\sigma};\mathbf{1})\,\mathrm{d}v\,\mathrm{d}n.
\end{multline*}

Further, the resulting integral over $N_{ab-a,a}$ vanishes
because of the $N_{a^b}$-invariance of the section 
$\phi_{\pi,\sigma \mathbf{s}}$ of $\mathrm{Ind}_{P_{a^b}}
^{\mathrm{GL}_{ab}}(\pi,\sigma\mathbf{s})$.
and because of the normalization of the measures.
Thus, 
\begin{multline*}
\psi^w_{\Delta}(t;\boldop r,\tilde{r}\boldcp h;\ell)=\\
\mathrm{Factor}1\times |\det \ell|^{-q_1}
\int\limits_{v\in V^h}|\det m(v)|^s \mathrm{Eval}_{\mathbf{s}=\Lambda_b}|\det 
\boldop1_{ab-a,a},{}^t t^{-1} \boldcp|^{\sigma \mathbf{s}}
\delta_{a^b}^{\frac{1}{2}}(\boldop1_a,\ldots,1_a,{}^t t^{-1} \boldcp)
\\ 
\phi^v_{\pi,\sigma\mathbf{s}}(\boldop \ell r,1_a\boldcp w_{\sigma}; \boldop 1_{ab-a},{}^t t^{-1}\boldcp)
\,\mathrm{d}v=\\ \mathrm{Factor}1\times
|\det t|^{(b-1)(a-1)/2}
\times\\ |\det \ell|^{-q_1}
\int\limits_{v\in V^h}|\det m(v)|^s \mathrm{Eval}_{\mathbf{s}=\Lambda_b}
\phi^v_{\pi,\sigma\mathbf{s}}(\boldop \ell r,1_a\boldcp w_{\sigma}; \boldop 1_{ab-a},{}^t t^{-1}\boldcp)
\,\mathrm{d}v,
\end{multline*}
We have the following equality for 
the factor $\delta_{ab}^{\frac{1}{2}}(\boldop {}^t t^{-1},r  \boldcp)$
which appears in `Factor 1'.
\begin{equation}\label{eqn:halfmodularchar}
\delta_{ab}^{\frac{1}{2}}(\boldop {}^t t^{-1},r  \boldcp)=
|\det t|^{-\frac{ab+1}{2}}|\det r|^{\frac{ab+1}{2}}.
\end{equation}
Now, we \textit{define}
\begin{multline}\label{eqn:psiupperwdefn}
\psi_{\pi}^w(r;\mathbf{m};t):= |\det r|^u\int\limits_{v\in V^h}|\det m(v)|^s
\phi_{\pi,\mathbf{s}}^v
(\boldop r,1_a\boldcp w_{\sigma};\boldop \mathbf{m},
{}^t t^{-1} \boldcp)\,\mathrm{d}v,\\ \text{for}\; r\in\mathrm{GL}_{a(b-1)},\; \mathbf{m}\in(\mathrm{GL}_a)^{b-1},\;\text{and}\; t \in\mathrm{GL}_a,
\end{multline}
We claim that, assuming this integral converges,
\[
\psi_{\pi}^w\in \mathrm{Ind}^{\mathrm{GL}_{a(b-1)}}
_{P_{a^{b-1}}}(\pi^{b-1},\mathbf{v}')\otimes\tau ,\;
\text{with}\;\mathbf{v}'=\mathbf{s}'+\mathbf{u}'+\frac{1}{2}\mathbf{a}'.
\]
Here, we are using the following notation for truncations
of $b$-tuples,
\[
\mathbf{s}=(s_1,\ldots, s_b)\in\mathbf{C}^b\;\text{implies}\;
\mathbf{s}'=(s_1,\ldots, s_{b-1})\in\mathbf{C}^{b-1}.
\]
Note that this notation differs from that used in the 
parallel part of the discussion regarding the identity term,
where the truncations occurred by eliminating the first
coordinate instead of the last ($b^{\rm th}$) coordinate.
We calculate that
\begin{multline*}
\psi^w_{\pi}(\mathbf{m}_2 r;\mathbf{m}_1;t)=|\det r|^u|\det \mathbf{m}_2|^{\mathbf{u}'}\cdot
\int\limits_{v\in V^h}|\det m(v)|^s\phi_{\pi,\sigma\mathbf{s}}^v(
\boldop \mathbf{m}_2,1_a\boldcp  \boldop r,1_a \boldcp w_{\sigma};\boldop \mathbf{m}_1,
{}^t t^{-1}\boldcp)\,\mathrm{d}v=\\
|\det r|^u
\delta_{a^b}^{\frac{1}{2}}(\boldop \mathbf{m}_2,1_a\boldcp)
|\det\mathbf{m}_2|^{\mathbf{u}'+\mathbf{s}'}\int\limits_{v\in V^h}|\det m(v)|^s
\phi_{\pi,\mathbf{s}}^v(\boldop r,1_a\boldcp w_{\sigma} ; \boldop
\mathbf{m}_1\mathbf{m}_2,{}^t t^{-1}\boldcp)\,\mathrm{d}v=\\
|\det r|^{\mathbf{u}}|\det \mathbf{m}_2|^{\mathbf{u}'+\mathbf{s}'+a\Lambda_b'}
\int\limits_{v\in V^h}|\det m(v)|^s\phi_{\pi,s}^v(\boldop r,1_a \boldcp
w_{\sigma};\boldop \mathbf{m}_1\mathbf{m}_2, {}^t t^{-1}\boldcp)\,\mathrm{d}v.
\end{multline*}
Since $\delta^{\frac{1}{2}}_{a^{b-1}}(\mathbf{m}_2)=|\det \mathbf{m}_2|^{a\Lambda_{b-1}}$
and $a(\Lambda_b'-\Lambda_{b-1})=\frac{1}{2}\mathbf{a}'$,
\[
\psi_{\pi}^w(\mathbf{m}_2 r;\mathbf{m}_1 ;t)=|\det r|^{\mathbf{u}}|\det\mathbf{m}_2|
^{\mathbf{u}'+\mathbf{s}'+\frac{1}{2}\mathbf{a}'}\psi_{\pi}^w(r,\mathbf{m}_1\mathbf{m}_2;t).
\]
In particular, with $\mathbf{s}=\sigma\Lambda_b$, 
$\mathbf{v}'=\sigma'\Lambda_{b-1}$, we have
\[
u=\frac{1}{2}(1-a).
\]
We have shown that for $\phi_{\pi,\sigma s}\in\mathrm{Ind}(\pi^b
,\sigma \mathbf{s})$,
\begin{multline*}
\left((r;\mathbf{m};t)\mapsto
\mathrm{Eval}_{\mathbf{s}=\Lambda_b}|\det r|^{\frac{1}{2}(1-a)}
\int\limits_{v\in V^h}|\det m(v)|^s\phi_{\pi,\mathbf{s}}
(\boldop r,1_a\boldcp w_{\sigma};\boldop \mathbf{m},
{}^t t^{-1} \boldcp)
\,\mathrm{d}v\right)\\ \in \mathrm{Eval}_{\mathbf{v}'=\Lambda_{b-1}}
(\mathrm{Ind}\pi^{b-1},\sigma\mathbf{v}')\otimes\tau.
\end{multline*}
Thus we have shown that
\begin{multline*}
\psi_{\Delta}^w(t;\boldop r,\tilde{r}\boldcp h,\ell)=
\mathrm{Factor}1\times |\det t|^{(b-1)(a-1)/2}\times |\det \ell|^{-q_1}
|\det r|^{\frac{1}{2}(a-1)}|\det \ell|^{\frac{1}{2}(a-1)}\\
|\det r|^{\frac{1}{2}(1-a)}|\det \ell|^{\frac{1}{2}(1-a)}
\int\limits_{v\in V^h}|\det m(v)|^s
\mathrm{Eval}_{\mathbf{s}=\Lambda_b}\phi_{\pi,\sigma \mathbf{s}}
(\boldop \ell r,1_a \boldcp w_{\sigma}; \boldop 1_{ab-a},{}^t t^{-1}\boldcp )\,\mathrm{d}v=\\
\mathrm{Factor}1\times |\det t|^{(b-1)(a-1)/2}\times |\det \ell|^{-q_1}
|\det r|^{\frac{1}{2}(a-1)}|\det \ell|^{\frac{1}{2}(a-1)}\\
\mathrm{Eval}_{v'=\Lambda_{b-1}}\phi^h_{\pi^{b-1},v'}(\ell r',1)\otimes\phi_{\tau}(t)=\\
\mathrm{Factor}1\times |\det t|^{(b-1)(a-1)/2}\times |\det \ell|^{-q_1+\frac{1}{2}(a-1)}
|\det r|^{\frac{1}{2}(a-1)}\phi_{\Delta^{b-1}}(h;\ell r)\otimes\phi_{\tau}(t)=\\
\delta_{n}^{\frac{1}{2}}(\boldop{}^t t^{-1},r \boldcp)|\det r|^{\frac{1}{2}(a-1)+s-a}\times \\
|\det t|^{-q_2-s+ab+1+(b-1)(a-1)/2}
|\det \ell|^{-q_1+\frac{1}{2}(a-1)}\phi_{\Delta^{b-1}}
(h;\ell r)\otimes\phi_{\tau}(t).
\end{multline*}
If we set
\begin{equation}\label{eqn:q2value}
q_2=-s+ab+1+(b-1)(a-1)/2-(ab+1)/2=-s+ab-(a+b)/2+1.
\end{equation}
then
\[
\psi_{\Delta}^w(t;\boldop r,\tilde{r}\boldcp h,\ell)=
|\det r|^{\frac{1}{2}(a-1)+s-a+\frac{ab+1}{2}}|\det \ell|^{-q_1+\frac{1}{2}(a-1)}
\phi^{\varphi}_{\Delta^{b-1}}(h;\ell r)\otimes\phi_{\tau}(t).
\]
Set
\begin{equation}
q_1=\frac{1}{2}(a-1)
\end{equation}
and note that
\[
|\det r|^{\frac{1}{2}(a-1)+s-a+\frac{ab+1}{2}}=|\det r|^{s-\frac{1}{2}}\delta_{ab-a}^{\frac{1}{2}}(r).
\]
For that choice of $q_1$, we have 
\[
\psi_{\Delta}^w(t;\boldop r,\tilde{r}\boldcp h,\ell)=
|\det r|^{s-\frac{1}{2}}\delta_{ab-a}^{\frac{1}{2}}(r)
\phi_{\Delta^{b-1}}(h;\ell r)\otimes\phi_{\tau}(t).
\]
We repeat the calculation with $r$ replaced by $1\in\mathrm{GL}_{a(b-1)}$
and $\ell$ replaced by $\ell r$.  Then we obtain
\[
\psi_{\Delta}^w(t;\boldop r,\tilde{r}\boldcp h,\ell)=
|\det r|^{s-\frac{1}{2}}\delta_{ab-a}^{\frac{1}{2}}(r)\psi_{\Delta}^w(t; h,\ell r),
\]
as claimed.  So we conclude from \eqref{eqn:wtermabstractform} and \eqref{eqn:q2value}
that
\begin{multline}\label{eqn:wtermcomputed}
E_Q^{ab}(f_{\Delta,s},\mathrm{diag}(t,h,\tilde{h},\tilde{t}))_{w}=
|\det t|^{-s+ab-(a+b)/2+1}
E^{a(b-1)}(h,f^{b-1}_{s-\frac{1}{2}})\otimes \phi_{\tau}(t),\\ \text{where}\; ,f^{b-1}_{s-\frac{1}{2}}:=\Phi_{\Delta^{b-1}}(\cdot;1)
|\det |^{s-\frac{1}{2}}.
\end{multline}
We have deferred the question of when the
integral \eqref{eqn:psiupperwdefn} converges.  We
will return to this question in the immediately subsequent section.

\vspace*{0.3cm}
\noindent
\textbf{Summing up.}
Adding the terms from \eqref{eqn:identitytermcomputed} and \eqref{eqn:wtermcomputed}, we obtain
\begin{multline}\label{eqn:nonnormalizedconstantterm}
E_Q^{ab}(f_{\Delta,s},\mathrm{diag}(t,h,\tilde{t}))=
\phi_{\tau}(t)\otimes |\det t|^{s+ab-(a+b)/2+1}E^{a(b-1)}(h,f'_{\Delta,s+\frac{1}{2}})+\\
|\det t|^{-s+ab-(a+b)/2+1}
 E^{a(b-1)}(h,f^{b-1}_{\Delta,s-\frac{1}{2}})\otimes \phi_{\tau}(t),
\end{multline}

\section{Normalizing factors for intertwining operator.}
As a result of \eqref{eqn:nonnormalizedconstantterm},
we can conclude that the integral in the variables $X$ and $Z$ in \eqref{eqn:psiupperwdefn}
(the $w=w_a$ term of $E_Q^{ab}$), gives
a specific intertwining operator for values of $s$ at which the intergral
converges.  For values of $s$ at which the defining
integral converges, the intertwining operator is given by
\begin{equation}\label{eqn:intertwiningopglobalint}
(\mathscr{U}_{\hat{w_{a}}}^{ab}f_{\Delta,s})(g)=\int\limits_{\begin{matrix}{\scriptstyle X\in \mathbf{A}^{a^2(b-1)}}\\
{\scriptstyle Z\in \mathbf{A}^{\frac{a(a+1)}{2}}}
\end{matrix}}
f^{\phi}_{\Delta,s}
(\hat{w_a}v''(X,Z) g)\,\,\mathrm{d}X\,\mathrm{d}Z,\; \text{for}\; \phi\in \mathrm{Ind}_{P}^G(\Delta,s),
\end{equation}
where we define
\[
\hat{w_{a}}=m(w_{\sigma})w_a
\]
with $w_{\sigma}$ as in \eqref{eqn:wsubsigmadefn},
$m(\cdot)$ the standard injection of $M$ into $G$
(so that $m(w_{\sigma})=\boldop w_{\sigma},
\tilde{w_{\sigma}} \boldcp$), and $w_{a}\in W(G,P_{a^b})$
whose action is reversing the sign of the first $a$ coordinates.
For the matrix computations below it will be useful for us
to record the following explicit representations:
\[
m(w_{\sigma})=\begin{pmatrix}&1_{a(b-1)}&&\\
1_a&&&\\
&&&1_a\\
&&1_{a(b-1)}&
\end{pmatrix}\;\; \hat{w_a}=\begin{pmatrix}
&1_{a(b-1)}&&\\
&&&1_a\\
-1_a&&&\\
&&1_{a(b-1)}&
\end{pmatrix}\]
\[ m(w_{\sigma})w_a m(w_{\sigma}^{-1})=
\begin{pmatrix}
1_{a(b-1)}&&&\\
&&-1_a&\\
&1_a&&\\
&&&1_{a(b-1)}
\end{pmatrix}=w_{[a(b-1)+1,ab]},
\]
where the subscript $[a(b-1)+1,ab]$ indicates
that the element of $W(G_{ab})$ considered
as a permutation/sign change reverses 
the sign of the last $a$ coordinates.

When the dimensions are fixed, we will most
frequently denote the intertwining operator
more simply by $\mathscr{U}_{\hat{w}}$, and using
this notation, we have shown in the paragraph on the
$w$-term of the Eisenstein series in \S \ref{subsec:nonnormalizedconstterm}
that
\begin{equation}\label{eqn:Uintertwiningoperatordomainrange}
i^*\circ \mathscr{U}_{\hat{w}}: \mathrm{Ind}_{P_{ab}}^{G_{ab}}(\Delta(\tau,b),s)\mapsto
\mathrm{Ind}_{P_{a(b-1)}}^{G_{a(b-1)}}\left(\Delta(\tau,b-1),s-\frac{1}{2}\right),
\end{equation}
for values of $s$ at which the integral converges.  The meromorphically
continued intertwining operator
has poles at certain values of $s$.  The aim of this section is to construct
a certain \textbf{normalizing factor} $\gamma^{ab}_w(s)$ for $\mathscr{U}_{\hat{w}}$
such that the product,
\begin{equation}\label{eqn:normalizedintertwiningdefn}
\mathscr{U}_{\hat{w}}^*(s):=[\gamma^{ab}_w(s)]^{-1}\mathscr{U}_{\hat{w}}^{ab}(s),
\end{equation}
is holomorphic and nonvanishing in the right half-plane, \textit{i.e.}, following \eqref{eqn:Uintertwiningoperatordomainrange},
such that holomorphic sections in $\mathrm{Ind}(\Delta^b,s)$ are
mapped to holomorphic sections in $\mathrm{Ind}(\Delta^{b-1},s-1/2)$.  Thus, in order to find
$\gamma^{ab}_w(s)$, we must find a ratio $\gamma^{ab}_w(s)$ of products of $L$-functions
such that the set of poles (resp. zeros) of $\gamma^{ab}_w(s)$ equal the set of poles
(resp. zeros) of $U^{ab}_w(s)$.

The strategy we employ for constructing $\gamma^{ab}_w(s)$, following \S1 of \cite{kudlarallisfest} is 
based on the fact that the integral for $\mathscr{U}_{w}^{ab}$
is Eulerian.  If $v$ is a place of $k$, we define the corresponding local intertwining
operator $\mathscr{U}_v$ by the local integral analogous to \eqref{eqn:intertwiningopglobalint}.  Similarly to 
\S 3.2.4 of \cite{jiangmemoirs} we determine
the effect of $\mathscr{U}_w^{ab}(s)$ on \textit{spherical} sections
(to be defined in detail below) and show that on such sections $\mathscr{U}_v$
amounts to multiplication by a certain factor,
which we call $(\gamma^{ab}_w)_v$.  Letting $S$ be the finite
set of places, including Archimedean ones, outside of which $\tau$
is spherical, we define 
\[
\gamma^{ab}_{w,S}(s)=\prod_{v\notin S}\gamma^{ab}_{w,v}(s).
\]
As is well known, $\gamma_{w,S}^{ab}$ 
turns out to be a ratio of partial $L$-functions.
We then complete these $L$-functions to form the completed normalizing factor $\gamma^{ab}_w(s)$,
and we define local factors $\gamma_{w,v}^{ab}(s)$
at all places, as the ratios of corresponding local $L$-factors.  Any
section $f_{\Delta,s}$ is $K_v$-fixed at almost all places, so that for some finite set of places $S'$ containing $S$
\begin{equation}\label{eqn:stsection}
f_{\Delta,s}=\left[\bigotimes_{v\in S'}f_{\Delta,s,v}\right]\otimes \left[\bigotimes_{v\notin S'}f_{\Delta,s,v}^0\right].
\end{equation}
Then by \eqref{eqn:normalizedintertwiningdefn}, 
\[
\mathscr{U}_w^{ab,*}(s)f_{\Delta,s}=\bigotimes_{v\in S'}\left[\frac{1}{\gamma_{w,v}^{ab}(s)}\mathscr{U}_v f_{\Delta,s,v}\right]\otimes 
\left[\bigotimes_{v\notin S'}\tilde{f}_{\Delta,s,v}^0\right],
\]
where $f_{\Delta,s,v}^0$ and $\tilde{f}_{\Delta,s,v}^0$ are the normalized spherical vectors in the local induced spaces.
After some preparation in \S\ref{subsec:inductioninstages}
we actually carry out the calculation of $\gamma_{w,v}^{ab}$
in \S\ref{subsec:sphericalsections}.
We then complete the proof that the globalized intertwining operator
defined by \eqref{eqn:normalizedintertwiningdefn} is holomorphic by proving that
for $v\in S'$, $\mathscr{U}_v$ itself is holomorphic. This has
to be done using a different method
and is the subject of \S\ref{subsec:normalizationramified}.

\subsection{Induction in stages.\newline}
\label{subsec:inductioninstages}

\noindent
\hspace*{-.2cm}\textbf{Residual Representation as Langlands Quotient.}
We now recall the classification of $L^2_d(\mathrm{GL}_n(k)\backslash \mathrm{GL}_n(\mathbf{A}))$,
established as the main result of \cite{moeglinwaldspurgerens}, using notation
of \S\ref{subsec:cuspidalexponents}.  Let $\mathfrak{E}^0$ denote the equivalence
classes of cuspidal data containing at least one element of the form $(M_{a^b},\tau^{\otimes b})$
for $\tau$ of the above type.  It is clear that there is exactly one such element in each 
orbit $\mathfrak{X}$ belonging to $\mathfrak{E}^0$, so we can define a map
$C$ on $\mathfrak{E}^0$ taking $\mathfrak{X}\in\mathfrak{E}^0$ to the unique
$(M_{a^b},\tau^{\otimes b})$ contained in  $\mathfrak{X}$.  It is shown in \S I.11 of \cite{moeglinwaldspurgerens}
that $\mathrm{Ind}_{P_{a^b}}^{\mathrm{GL}_{ab}}(\tau^{\otimes b},\Lambda_b)$ has a unique irreducible
quotient, which we will denote by
\[
J_{P(\mathbf{A})}^{\mathrm{GL}_{ab}(\mathbf{A})}(\tau^{\otimes b},\Lambda_b).
\]
Then the main result of \cite{moeglinwaldspurgerens} can be stated as follows.
\begin{thm}  \label{thm:mwmain} Let $\mathfrak{X}\in\mathfrak{E}$.  Then one has
\begin{itemize}
\item[(i)]  
If $\mathfrak{X}\notin \mathfrak{E}^0$, then $L^2(\mathrm{GL}_n(k)\backslash\mathrm{GL}_n(\mathrm{A}))_{\mathfrak{X}}
\cap L^2(\mathrm{GL}_n(k)\backslash\mathrm{GL}_n(\mathrm{A}))_{d}=\emptyset$.
\item[(ii)]
If $\mathfrak{X}\in \mathfrak{E}^0$, then
\[
L^2_d(\mathrm{GL}_n\backslash\mathrm{GL}_n(\mathrm{A}))_{\mathfrak{X}}=
\left\{ J_{P(\mathbf{A})}^{\mathrm{GL}_{ab}(\mathbf{A})}(\tau^{\otimes b},\Lambda_b)\right\},\;\text{where}\; C(\mathfrak{X})=(M_{a^b},\tau^{\otimes b}).
\]
\end{itemize}
Therefore (by the decomposition \eqref{eqn:discretepartdecomposition}), we have
\[
L^2_d((\mathrm{GL}_n\backslash\mathrm{GL}_n(\mathrm{A})))_{\xi}=
\bigoplus_{b|n}\bigoplus_{\tau \in\widehat{(\mathrm{GL}_a(\mathbf{A}))}'_{\rm cusp}}
 J_{P(\mathbf{A})}^{\mathrm{GL}_{ab}(\mathbf{A})}(\tau^{\otimes b},\Lambda_b),
\]
with the mapping $C$ giving a parametrization of the sum by $\mathfrak{E}^0$ .
In this sum, $\widehat{(\mathrm{GL}_a(\mathbf{A}))}'_{\rm cusp}$ denotes the collection
of automorphic representations of $\mathrm{GL}_a(\mathbf{A})$ which are unitary, cuspidal, and self-dual.
\end{thm}
One can use the definition of $\Delta^b=\Delta(\tau,b)$ to show that $\Pi_0(\Delta)$ consists
of the single element $\mathfrak{X}=(M_{a^b},\mathfrak{P}_{a^b})$ (the $X^G_M$-orbit of
$\pi^b=\tau^{\otimes b}$).  (The argument is similar to that used to prove Theorem \ref{thm:main} below.)
Further, as mentioned above,
Jacquet \cite{jacquetarticle}
used Langlands' criterion to show that $\Delta^b$ is square integrable.
Together, these results imply that
\begin{equation}\label{eqn:deltaJglobaldescription}
\Delta(\tau,b)\cong  J_{P(\mathbf{A})}^{\mathrm{GL}_{ab}(\mathbf{A})}(\tau^{\otimes b},\Lambda_b).
\end{equation}
the unique irreducible subquotient of the induced representation.

As a result of the definition of $\Delta(\tau,b)$ in 
\eqref{eqn:deltadefinition}, we have that
$J_{P(\mathbf{A})}^{\mathrm{GL}_{ab}
(\mathbf{A})}(\tau^{\otimes b},\Lambda_b)=
\mathrm{Res}_{\Lambda_b}^{\mathfrak{P}_{a^b}}
E^M(P_{a^b},\pi^b)$.  By using the formula for the constant
term of a cuspidal-data Eisenstein series in \S II.1.7 of \cite{mwbook},
we deduce that 
\[
J_{P(\mathbf{A})}^{\mathrm{GL}_{ab}
(\mathbf{A})}(\tau^{\otimes b},\Lambda_b)=
\Delta(\tau,b)\cong N(w_0,\Lambda^b,\pi^b)
\mathrm{Ind}_{P_{a^b}}^{\mathrm{GL}_{ab}}(\pi^b,\Lambda_b),
\]
where $w_0$ is the longest element in the Weyl
group $W(\mathrm{GL}_{ab},M_{a^b})$.
It is well known and not difficult to show induced
representations and normalized intertwining operators
factor as
\[
\mathrm{Ind}_{P_{a^b}}^{\mathrm{GL}_{ab}}(\pi^b,\Lambda_b)=
\bigotimes_v \mathrm{Ind}_{P_{a^b}}^{G_{ab}}(\pi_v^b,\Lambda_b)\;
\text{and}\; N(w_0,\Lambda^b,\pi^b)=\bigotimes_v
N_v(w_0,\Lambda^b,\pi_v^b).
\]
Therefore, 
\[
J^{\mathrm{GL}_{ab}(\mathbf{A})}_{P_{a^b}(\mathbf{A})}
\cong \bigotimes_v N_v(w_0,\Lambda_b,\pi^b_v)
\mathrm{Ind}_{P_{a^b}(k_v)}^{\mathrm{GL}_{ab}(k_v)}
(\pi_v^b,\Lambda_b).
\]
Thus the local component
\[
\Delta(\tau,b)_v\cong N(w_{\sigma},\Lambda_b,\pi^b_v)
\mathrm{Ind}_{P_{a^b}(k_v)}^{\mathrm{GL}_{ab}(k_v)}(\pi^b_v,\Lambda_b).
\]
Therefore, $\Delta(\tau,b)_v$ is the unique
irreducible quotient of $\mathrm{Ind}_{P_{a^b}(k_v)}^{\mathrm{GL}_{ab}(k_v)}(\pi^b_v,\Lambda_b)$, and we use the notation,
\[
\Delta(\tau,b)_v\cong J(\pi^b_v,\Lambda_b)
\]
for this local component, the same notation as for the Langlands
quotient.  

\noindent
\vspace*{.3cm}
\textbf{Remark.}  Note that $J(\pi^b_v,\Lambda_b)$ \textit{really} is the Langlands quotient if and only if $\tau_v$ is temepered, 
and if $\tau_v$
is not tempered, $J(\pi^b_v,\Lambda_b)$ is the Langlands
quotient coming from lower parabolic terms.  Compare
\cite{kimisraeljournal}, p. 266 second paragraph and p. 76
of \cite{spehmuellergafa}.  Because
\[
N(w_{\sigma},\Lambda_b,\pi^b_v):
\mathrm{Ind}_{P_{a^b}(k_v)}^{\mathrm{GL}_{ab}(k_v)}(\pi^b_v,\Lambda_b)\mapsto \mathrm{Ind}_{P_{a^b}(k_v)}^{\mathrm{GL}_{ab}(k_v)}(\pi^b_v,-\Lambda_b)
\]
we have the following Lemma.
\begin{lemma} \label{lem:generalplacesubmodulelem} 
Let $v$ be any place of $k$.  The local component $\Delta(\tau,b)_v$ is isomorphic to a submodule of
$\mathrm{Ind}_{P_{a^b}(k_v)}^{\mathrm{GL}_{ab}(k_v)}(\pi^b_v,-\Lambda_b)$.
\end{lemma}

When we consider the local component of $\mathrm{Ind}_{P}^G
(\Delta^b,s)$, we are therefore, considering a subrepresentation
of 
\[
\mathrm{Ind}_{P}^G(\mathrm{Ind}_{M_{a^b}}^M
(\pi^b_v,-\Lambda_b),s)
\]
This leads us to the following general considerations.

\vspace*{.3cm}\noindent
\textbf{Transitivity of Induction.}
This discussion takes place at a fixed place $v$
of $k$, which we henceforth drop from the notation.
Let $P$ be a standard parabolic of $G$ contained in a larger
standard parabolic $P'$, so that
\[
P\subseteq P'\subseteq P_0,\quad M_0\subseteq M'\subseteq M.
\]
Let $X_{M'}$,
$X_M$ be the groups of rational characters of $M'$ and
$M$, respectively.   As described for example in \S1.4 of \cite{mwbook},
the natural inclusion $M'\hookrightarrow M$
induces an inclusion $X_{M}\hookrightarrow X_{M'}$ in
the opposite direction. Our convention in this discussion is to
denote elements of $X_M$ by non-bold $\lambda$,
and the same elements when considered as elements
of $X_{M'}$ by $\boldlambda$.  This somewhat
conflicts with the overall conventions for the paper
because in general
$\lambda$ itself could be a vector.  Also, we will
so often induce from the parabolic $P'\cap M$
in $M$, that we will abuse notation by dropping
the intersection and simply consider $P'$ as a parabolic
in $M$.

\begin{prop}\label{prop:inductioninstages}\textbf{Transitivity of Induction.} Fix $\chi\in X_{M'}$,
 $\lambda\in X_M$, and $\rho$ a unitary representation
 of $M'$. Let
$V$ be the space of the induced representation $\mathrm{Ind}_{P'}^{M}(\rho,\chi)$, so that $V$ is a space of functions
from $M$ to the space of $\rho$ satisfying a transformation
law under elements of $M'$.
\begin{itemize}
\item[(a)]   Let $\varphi\in \mathrm{Ind}_{P}^{G}(V,\lambda)$,
considered as a function from $G$
to $V$ (see above) satisfying the usual law under elements of $M$.
Then the map of evaluation, at $1\in M$, in the second factor,
sends $\varphi$ to a function $\tilde\varphi\in
\mathrm{Ind}_{P'}^{G}(\rho,\chi+\boldlambda)$.
\item[(b)]  The one-to-one correspondence $\varphi\leftrightarrow \tilde{\varphi}$ induced in $(a)$ is an isomorphism of $G$-modules
\[
\mathrm{Ind}_{P}^{G}(\mathrm{Ind}_{P'}^{M}(\rho,\chi),\lambda)\cong
\mathrm{Ind}_{P'}^{G}(\rho,\chi+\boldlambda).
\]
\end{itemize}
\end{prop}
\textbf{Proof.}  By iterating the transformation
laws for $\tilde\varphi$ and an element of $V$, we obtain
\begin{multline*}
\tilde{\varphi}(m'g;r')=\varphi(m'g;(1;r'))=\\
\varphi(g;(1;r'm'))
[\delta(M,P')^{1/2}(m')][\delta(G,P)^{1/2}(m')]
m_P(m')^{\lambda}m_{P'}(m')^{\chi}
\varphi(g;(1;r'm'))=\\
m_{P'}(m')^{\chi+\mbox{\boldmath ${\scriptstyle \lambda}$
\unboldmath}}
m_{P'}(m')^{\mbox{
\boldmath $ {\scriptstyle \rho(G,P)}$\unboldmath}+\rho(M,M_0)}
\tilde{\varphi}(g';r'm').
\end{multline*}
The parameter
of the induced representation matches the parameter of the
iterated induced representation as claimed.  In order to complete
the proof of (a) we have to verify that the contribution of
the normalizing factors matches, which amounts to the equality
\[
\rho(G,P')=\mbox{\boldmath $\rho(G,P)$\unboldmath}+\rho(M,M'),
\]
where by definition $\rho(G,P)$ equals half the sum of the positive roots
of $M$ in $G$.  The positive roots of $M'$ in $G$ are partitioned
exactly into the roots
of $M'$ in $N$ the positive roots of $M'$ in $M$.  Further an element of $m'\in M'$
factors as an element $a_m\in A_M$ times an element $(m')^1\in M^1$ (see 
pp. 19--20 of \cite{mwbook}, in the number field case).  
On the one hand, elements of $\Phi^+(M',M)$ vanish
on $a_m$, and the sum of the roots of $M_0$
in $N$ take the same value on $a_m$
as the sum of the roots in $\Phi^+(M,G)$, accounting for
multiplicity, so that the value of $\rho(G,P)$ on $a_m$
equals the value of $\rho(G,P')$ on $a_m$.  On the other hand, the roots
of $M'$ in $N$ vanish on $(m')^1$,
so that the value of $\rho(M,M')$ on $(m')^1$
equals the value of $\rho(G,P')$ on $(m')^1$.

(b) Given an element $\phi$ in $\mathrm{Ind}_{P'}^{G}(\rho,\chi+\boldlambda)$,
define $\varphi\in\mathrm{Ind}_{P}^{G}(V,\lambda)$
by setting $\varphi(g;(m,m'))=|\det m|^{-\lambda}
[\delta^{-1/2}(G,P)(m)]\phi(mg';m')$.  Then it is
easy to check that $\tilde{\varphi}=\phi$, so
we have defined an inverse mapping to the mapping
of (a).

The $G$-action on functions in either space occurs as
translation on the right in the first factor, and thus
intertwines with the mapping of (a).
\qed

\textbf{Remark.} For a more abstract statement
in a somewhat more general setting than we
consider here see
 \textbf{Exercise 4.5.8} of \cite{bumpbook}.

As a first application of Proposition
\ref{prop:inductioninstages}
we have
\begin{cor}\label{cor:arbitraryplaceindstages} 
At any place $v$ of $k$ we have an isomorphism
of $G_{ab,v}$-modules
\[
\mathrm{Ind}_P^G(\Delta(\tau,b)_v,s)\cong \mathrm{Ind}_P^G
(N(w_{\sigma},\Lambda_b,\pi^b)\mathrm{Ind}_{P_{a^b}(k_v)}
^{\mathrm{GL}_{ab}(k_v)}(\pi^b_v,\Lambda_b),s),
\]
so that $\mathrm{Ind}_P^G(\Delta(\tau,b)_v,s)$ is isomorphic
to a submodule of $\mathrm{Ind}_{P_{a^b}(k_v)}
^{G_{ab}(k_v)}(\pi^b_v,\mathbf{s}-\Lambda_b)$.
\end{cor}

\vspace*{.3cm}
\noindent
\textbf{An identity of intertwining operators.}
In order to determine the action of the local
intertwining operator $\mathscr{U}_{\hat{w},v}:=\mathscr{U}_{\hat{w},v}^{ab}$
on a normalized spherical section, we relate
it to a more standard intertwining operator $M^{ab}_v$
by proving an identity which holds between these
local operators at all places.  Since the place $v$
is fixed but arbitrary throughout this argument,
we do not mention it for the rest of this paragraph.

First, define $w_0$ to be the longest
element of the Weyl group of $G$ with respect to $M$.
One represents $w_0$ as
\[
w_0=\begin{pmatrix}
0&1_{ab}\\
-1_{ab}&0
\end{pmatrix}
\]
We define the intertwining operator $M_{w_0}^{ab}$,
henceforth called simply $M^{ab}$.
For $\mathrm{Re}(s)$ sufficiently large, set
\[
M^{ab}\Phi_{\Delta,s}(g)=\int_{n\in N_{ab}^{ab}}
\Phi_{\Delta,s}(w_0ng)\,\mathrm{d}n\;\text{for all}\;
\Phi_{\Delta,s}\in \mathrm{Ind}_{P}^G(\Delta(\tau,b),s),\;
g\in G_{ab}.
\]
One checks, using Lemma \ref{lem:generalplacesubmodulelem},
that
\[
M^{ab}: \mathrm{Ind}_{P}^G(\Delta(\tau,b),s)\rightarrow
\mathrm{Ind}_P^G(\Delta(\tau,b),-s).
\]
In order to carefully distinguish
the different dimensions, we denote by $i_{a(b-1)}$
or even $i^{ab}_{a(b-1)}$ 
the ``diagonal" inclusion mapping of $G_{a(b-1)}$
into $G_{ab}$, which is defined by
\[
g'\mapsto \boldop 1_a,g',1_a \boldcp\;\text{for all}\; g'\in G_{a(b-1)}.
\]
Denote by $i_{a(b-1)}^*$ map induced by this
inclusion on sections of the induced representation
spaces.  For example,
when
$\Phi_{\Delta^b,s}\in \mathrm{Ind}_{P}^G(\Delta(\tau,b),s)$
we define
\[
(i_{a(b-1)}^*\Phi_{\Delta^b,s})(g')=\Phi_{\Delta^b,s}(i
_{a(b-1)}(g')) \;\text{for all}\; g'\in G_{a(b-1)}.
\]
Using the relation \eqref{eqn:rhoabcalculation},
one verifies that $i_{a(b-1)}$ induces a $G_{a(b-1)}$-intertwining map
\[
i_{a(b-1)}^*: \mathrm{Ind}^G_{P_{a^b}}(\pi^b,\mathbf{s}-\Lambda_b)
\rightarrow \mathrm{Ind}^{G_{a(b-1)}}_{P_{a^{b-1}}}\left(\pi^{b-1},\mathbf{s}-
\Lambda_{b-1}+\frac{1}{2}\right).
\]
From now until the end of this paragraph, for readability, we drop the
superscript and subscript from $\mathrm{Ind}$ whenever
the group and parabolic subgroup are clear from the context.
Therefore, by Corollary \ref{cor:arbitraryplaceindstages},
we have
\[
i_{a(b-1)}^*: \mathrm{Ind}(\Delta^b,s)
\rightarrow \mathrm{Ind}\left(\Delta^{b-1},s+\frac{1}{2}\right).
\]
We now define an intertwining operator
$M^{ab}_{a(b-1)}$, for $\mathrm{Re}(s)$ large,
by the integral
\begin{multline*}
M_{a(b-1)}^{ab}\Phi_{\Delta^b,s}(g)
:=\int\limits_{n\in N_{a(b-1)}^{a(b-1)}}
\Phi_{\Delta^b,s}(i_{a(b-1)}(w_0^{b-1})i_{a(b-1)}(n)g)\mathrm{d}n,
\\ \text{for all}\; \Phi_{\Delta^b,s}\in \mathrm{Ind}(\Delta,s).
\end{multline*}
The following calculations show that $M_{a(b-1)}^{ab}$
is related to $M^{ab}$, in
the same way that, in \S4 of \cite{krcrelle88}, the operator
$M^{n}_{n-1}$
is related to $M_n^n$.  In order to carry them out
we will have need for an expression of 
the half-sum of positive roots $\rho(G_{ab},P^{ab}_{a,ab-a})$,
which we call $\rho_{a,ab-a}$.  An elementary computation
yields $\rho_{a,ab-a}=(ab,\rho_{a(b-1)})$, meaning
\[
\delta^{\frac{1}{2}}(G_{ab},P^{ab}_{a,ab-a})(\boldop
t,m,\tilde{m},\tilde{t}\boldcp)=
|\det t|^{ab} |\det m|^{\rho_{a(b-1)}},
\]
where $t\in\mathrm{GL}_a$ and $m\in \mathrm{GL}_{a(b-1)}$.
Now in order to determine the image of $M_{a(b-1)}^{ab}$,
we set $m=\boldop t_1,\ldots t_{b-1}\boldcp$ for $t_i\in \mathrm{GL}_a$ and calculate, using the abbreviations $w_0':=i_{a(b-1)}(w_0^{b-1})$,
$n':=i_{a(b-1)}(n)$,
\begin{multline*}
M_{a(b-1)}^{ab}\Phi_{\Delta,s}(\boldop t,m,\tilde{m},\tilde{t} \boldcp g)=
\int\limits_{n'\in N_{a(b-1)}^{a(b-1)}}
\Phi_{\Delta^b,s}(w_0'n'\boldop t,m,\tilde{m},\tilde{t} \boldcp
g)\mathrm{d}n'=\\
\int\limits_{n'\in N_{a(b-1)}^{a(b-1)}}
\Phi_{\Delta^b,s}(
\boldop t,1,\tilde{t} \boldcp \mathbf{c}(w_0')(i\boldop m,\tilde{m}
\boldcp ) w_0' \mathbf{c}(i\boldop m,\tilde{m}\boldcp^{-1})n'g
)\mathrm{d}n'=\\
|\det t|^{s-(\Lambda_b)_1+(\rho_b^{(a)})_1}
(1,\tilde{t}_{b-1},\ldots \tilde{t}_1)^{\mathbf{s}-\Lambda_b+\rho_b^{(a)}}
\boldop t_1,\ldots t_b\boldcp^{2\rho_{a(b-1)}}
\phi_{\pi^b}(t,\tilde{t}_{b-1},\ldots \tilde{t}_1)M_{a(b-1)}^{ab}\Phi_{\Delta^b,s}(g)=\\
|\det t|^{s-\frac{1-b}{2}+a(b-\frac{1}{2})+\frac{1}{2}}
(t_1,\ldots, t_{b-1})^{\mathbf{-s-\frac{1}{2}}-\Lambda_{b-1}+
\rho_{b-1}^{(a)}}\phi_{\pi^b}(t,\tilde{t}_{b-1},\ldots \tilde{t}_1)
M_{a(b-1)}^{ab}\Phi_{\Delta^b,s}(g)=\\
|\det t|^{s-\frac{2-a-b}{2}+(\rho_{a,ab-a})_1}
(t_1,\ldots, t_{b-1})^{((-(\mathbf{s+\frac{1}{2}})-\Lambda_{b-1})+
\rho_{a}^{b-1})+\rho_{a(b-1)}}
\phi_{\pi^b}(t,\tilde{t}_{b-1},\ldots \tilde{t}_1)M_{a(b-1)}^{ab}\Phi_{\Delta^b,s}(g),
\end{multline*}
where in the last line we have used \eqref{eqn:rhosubabcalculation}.
Therefore, by the explicit calculation of $\rho_{a,a(b-1)}$
above, and using transitivity of induction, we have shown that
\[
M_{a(b-1)}^{ab}\Phi_{\Delta^b,s}
\in \mathrm{Ind}\left(\tau\otimes\mathrm{Ind}(\pi^{b-1},-\Lambda_{b-1}),
\left(s-\frac{2-(a+b)}{2},-s-\frac{1}{2}\right)\right).
\]
Therefore by Corollary \ref{cor:arbitraryplaceindstages}, 
we deduce that
\[
M_{a(b-1)}^{ab}: \mathrm{Ind}(\Delta^b,s)\rightarrow
 \mathrm{Ind}
\left(\tau\otimes\Delta^{b-1},\left(s+\frac{2-(a+b)}{2},-s-\frac{1}{2}
\right)\right).
\]

We fit this map into a commutative diagram,
generalizing the diagram in (4.45)
of \cite{krcrelle88},
\[
\begin{CD}
\mathrm{Ind}(\Delta^b,s)   @>M^{ab}_{a(b-1)}>> \mathrm{Ind}
\left(\tau\otimes\Delta^{b-1},\left(s+\frac{2-(a+b)}{2},-s-\frac{1}{2}
\right)\right)\\ 
@VV{ i^*_{a(b-1)}}V                                    @VV{ i^*_{a(b-1)}}V\\ 
\mathrm{Ind}(\Delta^{b-1},s+\frac{1}{2})                  @>
M^{a(b-1)}_{a(b-1)}>>    
\mathrm{Ind}(\Delta^{b-1},-s-\frac{1}{2}))
\end{CD}
\]
We derive the following analogue of Lemma 7.4 in
\cite{krcrelle88}.
\begin{lemma}  \label{lem:krUMrelation}
Using the notation $i^*=i^{a,*}_{a(b-1)}$, we have
\[
M_{a(b-1)}^{a(b-1)}\left(s-\frac{1}{2}\right)\circ i^*\circ
\mathscr{U}_{\hat{w}}(s)=i^*\circ M_{ab}^{ab}(s).
\]
\end{lemma}
\textbf{Proof.}  From the commutative diagram one has 
\[
M_{a(b-1)}^{a(b-1)}\left(s-\frac{1}{2}\right)\circ i^*_{a(b-1)}=
i^*_{a(b-1)}\circ M^{ab}_{a(b-1)}(s-1).
\]
(N.B., there is a misprint on the right-hand side of the corresponding
equation of p. 57 of \cite{krcrelle88}, where $s$ should be
$s-1$).
So we are reduced to showing that
\begin{equation}\label{eqn:intertwiningidentitydecomp}
M^{ab}_{a(b-1)}\circ \mathscr{U}_{\hat{w}}(s)=M_{ab}^{ab}(s).
\end{equation}
We perform a decomposition of integrals
similar to that carried out on pp. 57--8 of \cite{krcrelle88}.
By the definition
of $M^{ab}_{a(b-1)}$ and $\mathscr{U}_{w_a}$, one has
for $\mathrm{Re}\,s$ large,
\begin{equation}\label{eqn:intertwininginteridentity}
M^{ab}_{a(b-1)}(\mathscr{U}_{w_a}\Phi_s(g))=
\int_{n\in N^{a(b-1)}_{a(b-1)}(\mathbf{A})}\int\limits_{\begin{matrix}{\scriptstyle X\in \mathbf{A}^{a^2(b-1)}}\\
{\scriptstyle Z\in \mathbf{A}^{\frac{a(a+1)}{2}}}
\end{matrix}}
\Phi_{\Delta,s}
(w_av''(X,Z)i(w_0^{a(b-1)})i(n) g)\,\,\mathrm{d}X\,\mathrm{d}Z\,
\mathrm{d}n.
\end{equation}
Now one calculates that
\[
i(w_0^{a(b-1)})^{-1}v''(X,Z)i(w_0^{a(b-1)})=v'(0,-X,Z),
\]
where the notation is as in \eqref{eqn:vprimenotation}.
Further, for
\[
i(n)=i(n(Y))=\begin{pmatrix}1_a&0&0&0\\
0&1_{a(b-1)} &Y &0\\
0 & 0&  1_{a(b-1)} & 0\\
0 &0 &0& 1_a\end{pmatrix},\]
with $Y$ symmetric around the second diagonal,
one has
\[
v'(0,X,Z)i(n(Y))=\begin{pmatrix}1_a&0&-X&Z\\
0&1_{a(b-1)} &Y &(-X)'\\
0 & 0&  1_{a(b-1)} & 0\\
0 &0 &0& 1_a\end{pmatrix},
\]
which ranges over $N_{ab}^{ab}(\mathbf{A})$ as $n(Y)$
ranges over $N_{a(b-1)}^{a(b-1)}(\mathbf{A})$ and, separately,
$v''(X,Z)$ ranges over $V^w(\mathbf{A})$.  Therefore,
the integral in \eqref{eqn:intertwininginteridentity},
is actually
\[
\int_{N_{ab}^{ab}(\mathbf{A})}\Phi_{\Delta,s}(w_a i(w_0^{a(b-1)}) ng)\,\mathrm{d}n.
\]
Since a simple matrix
computation shows that $\hat{w} i(w_0^{a(b-1)})=w_0^{ab}$ (\textit{i.e.}
left multiplication transforms the Weyl element for the nontrivial
$M$ intertwining operator associated to $G_{a(b-1)}$ to
the corresponding element for $G_{ab}$),
we have produced
the integral defining $M_{ab}^{ab}\Phi_{\Delta,s}$
for $\mathrm{Re}\, s$ large.
This completes the proof of \eqref{eqn:intertwiningidentitydecomp}
and of the lemma.
\qed

\vspace*{.3cm}
\noindent
\textbf{Unramified places of $\tau_v$.}
Recall the following basic result.
\begin{thm}  (3.3.3 in \cite{bumpbook})  Let $(V,\tau)$ be an irreducible admissible representation
of $\mathrm{GL}_n(\mathbf{A})$.  Then there exists for each Archimedean place
$v$ of $k$ an irreducible admissible $(\mathfrak{g}_{\infty},K_{\infty})$-module $(\tau_v,V_v)$
and for each non-Archimediean place $v$ there exists an irreducible admissible representation 
$(\tau_v,V_v)$ of $\mathrm{GL}(n,k_v)$ such that for almost every $v$, $V_v$ contains a nonzero
$K_v$-fixed vector $\xi_v^0$ such that $\tau\cong \otimes_v'\tau_v$ (restricted tensor product).
\end{thm}

At such a place $v$, where a $K_v$-fixed vector exists, $\tau_v$ is said to be \textbf{spherical}. In particular,
we have for $(\tau,V)$ an irreducible cuspidal automorphic representation of $\mathrm{GL}_n(\mathbf{A})$
that $\tau$ induces an irreducible admissible representation of $\mathrm{GL}_n(\mathbf{A})$
on the space $V^K$ of $K$-finite vectors in $V$.  Therefore, the above tensor
product decomposition applies to the above 
self-dual, cuspidal automorphic representation $\tau=\tau|_{V_K}$ of $\mathrm{GL}_a(\mathbf{A})$,
and in particular the local factors $\tau_v$ are spherical for almost all finite places $v$.

Let $v$ be a fixed finite place such that $\tau_v$ is spherical.  It
is well-known 
that every spherical representation is a \textit{subrepresentation} of a \textit{spherical
principle series representation}, meaning a representation of the form
\[
\mathrm{Ind}_{P_a}^{G_a}(\chi_\mathbf{s})\quad \text{(normalized induction)}
\]
where $\chi_{\mathbf{s}}$ is an \textit{unramified character} of $T_a$,
extended to $P_a=T_a\ltimes N_a$ by triviality on the normal $N_a$ factor.

\vspace*{.3cm}
\noindent
\textit{Remark.}  When the spherical principal series representation $\mathrm{Ind}\chi_\mathbf{s}$ associated to $\tau_v$
is irreducible, we must have $\tau_v=\mathrm{Ind}\chi_\mathbf{s}$.  
It is understood that in ``most" cases--
\textit{i.e.} for spherical parameters $s_i$ which lying off a certain
finite union of hyperplanes--$\mathrm{Ind}\chi_{\mathbf{s}}$
is indeed irreducible.  See, \textit{e.g.}
\cite{garretturpsnotes}, following \cite{borel76} or \cite{casselman80},
for a statement of the exact conditions, at least in the case
when $\chi_s$ is ``regular".  However
these conditions
 need not concern us because we are dealing with submodules of induced
modules from $\chi_s$ anyway, so, in comparison to say
\S 3 of \cite{spehmuellergafa}, we cannot hope to achieve
anything stronger than a ``submodule" statement in the end.

Applying Proposition \ref{prop:inductioninstages} twice
at the spherical places give the following extension of
Corollary \ref{cor:arbitraryplaceindstages}.

\begin{prop} \label{eqn:unramifiedplacequotient}  Let  $v$ be an unramified place of $\tau_v$ such
that $\tau$ has spherical
parameter $\mathbf{s}=(s_1,\ldots, s_a)$.   Let $w_0$ be the longest Weyl element in the Weyl group of $\mathrm{GL}_{ab}$
with respect to $M_{a^b}$. 
\begin{itemize}
\item[(a)]  The local component
$\Delta(\tau,b)_v$ is a submodule of
\[ \mathrm{Ind}_{B_{ab}}^{\mathrm{GL}_{ab}}(\chi_{\mu(s,\Lambda_b)}),\;\text{where}\; \mu(s,\Lambda_b)=w_0\mathbf{s}-\mbox{\boldmath$\Lambda_b$ \unboldmath}.
\]
\item[(b)]
The local component of the representation induced from the residual representation to the symplectic
group $G$, $\mathrm{Ind}_P^G(\Delta(\tau,b),t)_v$, is a submodule of
\[
\mathrm{Ind}_{P_0}^{G_{ab}}(\mathrm{Ind}_{B_{ab}}^{\mathrm{GL}_{ab}}((\chi_{\mu(s,\Lambda_b)}),t)=
\mathrm{Ind}_{P_0}^{G_{ab}}(\chi_{w_0\mathbf{s}-\mbox{\boldmath$\scriptstyle\Lambda_b$ \unboldmath}\!\!\!+\mathbf{t}})))
\]
\end{itemize}
\end{prop}
\textbf{Proof.}
Part (a) follows from Proposition \ref{prop:inductioninstages}
in the same way as Corollary
\ref{cor:arbitraryplaceindstages}
Part (b) follows from part (a) and Proposition \ref{prop:inductioninstages}.
Note that because
\[
\mathbf{t}=\underbrace{(t,\ldots,t)}_{ab\;\text{times}},
\]
we have $\sigma_0\mathbf{t}=\mathbf{t}$, so the equality in (b) is valid.
\qed

For readability, from now on we will normally denote
$\chi_{w_0\mathbf{s}-\mbox{\boldmath$\scriptstyle\Lambda_b$ \unboldmath}\!\!\!+\mathbf{t}}$
as $\chi(w_0\mathbf{s}-\mbox{\boldmath$\scriptstyle\Lambda_b$ \unboldmath}\!\!\!+\mathbf{t})$, and similar
characters in the same way.

\subsection{Application of Shahidi's Method to calculate
the effect of $M(w,\Delta,s)$ on spherical sections}
\label{subsec:sphericalsections}
Let $\frac{a(s)}{b(s)}$ be the ratio of $L$-functions,
without common $L$-factors, by which 
$M_{w_0}(s)=M(w_0,\Delta^b,s)$
acts on normalized spherical vectors in
$\mathrm{Ind}(\Delta^b,s)$.

At spherical places Lemma \ref{lem:krUMrelation} implies the formula
\begin{equation}\label{eqn:gammaintermsaandb}
\gamma_v(s)=\frac{a_{b,v}(s)}{b_{b,v}(s)}\frac{b_{b-1,v}(s-1/2)}{a_{b-1,v}(s-1/2)}.
\end{equation}
It follows from Shahidi's formulation in \cite{shahidi88}
of Langlands' formula for the effect of intertwining operators
on spherical sections that for $f_s^0$ a spherical section
of $\mathrm{Ind}_{P}^G(\Delta(\tau,b),s)$,
\begin{equation}\label{eqn:Msphericalgeneral2}
M(w,\Delta,s)f_s^0=\frac{L_v(s,\Delta)L_v(2s,\Delta,\wedge^2)
}{L_v(s+1,\Delta)L_v(2s+1,\Delta,\wedge^2)}
\tilde{f}_s^0.
\end{equation}

\vspace*{.3cm}
\noindent
\textbf{Normalizing factor in case $b=1$.}
When we substitute $\tau$ for $\Delta$ into
\eqref{eqn:Msphericalgeneral2}, we obtain a ratio
of products of $L$-functions without common factors.
Thus, we can read off $a_{1,v}$ and $b_{1,v}$
directly from \eqref{eqn:Msphericalgeneral2}
and obtain
\begin{equation}\label{eqn:normalizingfactorsbeq1case}
b_{1,v}(\tau,s)=L_v(s+1,\tau)L_v(2s+1,\tau,\wedge^2)\;\text{and}
\;a_{1,v}(\tau,s)=L_v(s,\tau)L_v(2s,\tau,\wedge^2)
\end{equation}
for $v\notin S'$.

\vspace*{.3cm}
\noindent
\textbf{Normalizing factor in case $b\geq 2$.}
The idea for the computation of the normalizing
factor is similar to the case $b=1$, but the combinatorics
are more complicated.  We give a sketch of them.

For the computation of the local factors, the place
$v\notin S'$ of $k$ is fixed, so we drop it from the notation.
In general, the denominator of \eqref{eqn:Msphericalgeneral2}
factors as
\[
L(2s+1,\Delta,\wedge^2)=
\underbrace{\prod_{i=1}^b 
L(2s+b-2i+2,\tau,\wedge^2)}_{\mathrm{Factor\; I}}\times
\underbrace{\prod_{1\leq i<j\leq b}L(2s+b+2-(i+j),
\tau\otimes \tau)}_{\mathrm{Factor\; II}}
\]
In order to determine the cancellation,
the following charts count how many times the
denominator and numerator of \eqref{eqn:Msphericalgeneral2}
contains the `exterior square' factor\linebreak $L(2s+b+k,\tau,\wedge^2)$.

\begin{itemize}
\item Range of $k$.

\begin{tabular}{c|cc}
Range &Factor I& Factor II\\ \hline
$i,j\Rightarrow k$ & $k=2-2i$  & $k=2-(i+j)$\\
$k\Rightarrow i,j$ & $i=\frac{2-k}{2}$ & $(i+j)=2-k$\\
Range of $k$ & $-2b+2\leq k\leq 0$ & $-2b+4\leq k\leq -2$
\end{tabular}

\vspace*{.3cm}
\item $k$ even case.

\begin{tabular}{c|cccc}
Range &Factor I& Factor II & denominator & numerator \\ \hline
$2-b\leq k\leq 0$ & $1$  & $-\frac{k}{2}$ & $\frac{2-k}{2}$& $-\frac{k}{2}$\\
$-2b+2\leq k\leq -b$ & $1$ & $b-1+\frac{k}{2}$  & $b+\frac{k}{2}$
&$b+\frac{k}{2}$\\
\end{tabular}

\vspace*{.3cm}
\item $k$ odd case.

\begin{tabular}{c|cccc}
Range &Factor I& Factor II & denominator & numerator  \\ \hline
$-b+1\leq k\leq -1 $ &$0$& $\frac{1-k}{2}$&$\frac{1-k}{2}$ &  $\frac{1-k}{2}$\\
$-2b+3\leq k\leq -b-1$ &$0$& $b+\frac{k-1}{2}$ &$b+\frac{k-1}{2}$ & $b+\frac{1+k}{2}$
\end{tabular}
\end{itemize}

Now we compare the ``numerator" and "denominator"
columns to see when the numerator (resp., denominator)
of \eqref{eqn:Msphericalgeneral2}
nets an exterior-square factor.  We state the results
and below explain how the results for the symmetric square
and standard factors are arrived at by suitable modifications.
\begin{itemize}
\item $\wedge^2$ factor in case $b$ even.

\begin{itemize}
\item[$a_b$:]  $L(2s-1,\tau,\wedge^2)L(2s-3,\tau,\wedge^2)\cdots
L(2s-b+1,\tau,\wedge^2)$.
\item[$b_b$:]  $L(2s+b,\tau,\wedge^2)L(2s+b-2,\tau,\wedge^2)
\cdots L(2s+2,\tau,\wedge^2).$
\end{itemize}

\item $\vee^2$ factor in case $b$ even.
\begin{itemize}
\item[$a_b$:] $L(2s,\tau,\vee^2)L(2s-2,\tau,\vee^2)\cdots
L(2s-b+2,\tau,\vee^2)$.
\item[$b_b$:]  $L(2s+1,\vee^2,\tau)L(2s+3,\vee^2)\cdots
L(2s+b-1,\vee^2,\tau)$.
\end{itemize}

\item Standard factors.
\begin{itemize}
\item[$a_b$:] $L\left(s+\frac{1-b}{2},\tau\right)$.
\item[$b_b$:] $L\left(s+\frac{1+b}{2},\tau\right)$.
\end{itemize}
\end{itemize}

The formulas in case $b$ odd are sufficiently
similar that we do not repeat intermediate
steps.
The symmetric square factors are found
by constructing charts similar to the ones used
in finding the exterior square factors.  One
obtains these charts from the exterior square
charts by eliminating the contribution from ``Factor I"
(because this factor is a pure product of exterior-square
$L$-functions) and substituting $k-1$ for $k$.
The standard $L$-function factors
in the last two lines of the above chart
come directly from
the factorizations
\[
L(2s+1,\Delta)=\prod_{i=1}^b L(2s+1+\frac{b+1-2i}{2},\tau)\;
\text{and}\; L(2s,\Delta)=\prod_{i=1}^b L(2s+\frac{b+1-2i}{2},\tau).
\]
Combining the factors, we conclude that for $b\geq 2$,
\begin{itemize}
\item[]
\[
b_{b,v}(s)=\prod_{i=1}^{\lceil \frac{b}{2}\rceil}
L_v(2s+b+2-2i,\tau,\wedge^2)\prod_{i=1}^{\lfloor \frac{b}{2}\rfloor}
L_v(2s-1+2i,\tau,\vee^2)L_v\left(s+\frac{b+1}{2},\tau\right)
\]
\item[]
\[
a_{b,v}(s)=\prod_{i=1}^{\lceil \frac{b}{2}\rceil}
L_v(2s-b-1+2i,\tau,\wedge^2)\prod_{i=1}^{\lfloor \frac{b}{2}\rfloor}
L_v(2s-b+2i,\tau,\vee^2)L_v\left(s+\frac{1-b}{2},\tau\right).
\]
\end{itemize}
As a result of \eqref{eqn:gammaintermsaandb}, we deduce that
\begin{equation}\label{eqn:gammaSpartial}
\gamma_{b,S}(s):=\prod_{v\notin S}\gamma_{n,v}(s)=
\frac{L_S(2s,\tau,\vee^2)L_S(2s,\tau,\wedge^2)L_S(s+\frac{b-1}{2},\tau)}{L_S(2s+b-1,\tau,\vee^2)
L_S(2s+b,\tau,\wedge^2)L_S(s+\frac{b+1}{2},\tau)}.
\end{equation}

\subsection{Normalization at Ramified Places}
\label{subsec:normalizationramified}
Let $v$ be a ramified place of $\tau$, which
we will drop from the notation for the rest
of the subsection: our aim is to show
that though \textit{a priori}, only the normalized
version $\gamma^{-1}\mathscr{U}_{\hat{w}}$ of
$\mathscr{U}_{\hat{w}}$ is holomorphic on $\mathbf{C}$,
a closer analysis of $\mathscr{U}_{\hat{w}}$ itself actually
is holomorphic in the right-half plane $\mathrm{Re}s>0$.

Our strategy is inspired by
\S3.2.3 of \cite{jiangmemoirs}.  The first step is to use the
transitivity of induction to view $i^*\circ U_{\hat{w}}$
as the restriction of a certain ``nonstandard" intertwining
operator (also denoted by $i^*\circ U_{\hat{w}}$,
and defined in the range of convergence
by the same integral) on a
space induced from \textit{cuspidal} data on
$P_{a^b}$ to a space induced from \textit{cuspidal}
data on $P_{a^{b-1}}$.  (See Lemma \ref{lem:intertwiningembedding} 
below).
Although the intertwining
operator is not the \textit{standard intertwining operator
associated to a Weyl element $w\in W(G_{ab},P_{a^b})$},
we can nevertheless decompose $i^*\circ U_{\hat{w}}$ into a composition
of such intertwining operators on reductive subgroups
of $G_{ab}$, form the product $\lambda(s)$ of the normalizing factors
of these standard intertwining operators, and
thereby construct a function whose set of poles contain
(possibly properly) the set of poles of $\mathscr{U}_{\hat{w}}(s)$.

\vspace*{.3cm}
\noindent
\textbf{Application of Transitivity of Induction}
As a result of Proposition \ref{prop:inductioninstages},
we have that $i^*\circ \mathscr{U}_{\hat{w}}$ is initially
defined on a submodule of 
$\mathrm{Ind}_{P_{a^b}}^G(\pi^{b},\mathbf{s}-\Lambda_b)$.
By using the same integral
formula, we extend the action of $\mathscr{U}_{\hat{w}}$
to arbitrary elements of the induced space
$\mathrm{Ind}(\pi^b,\mathbf{s}-\Lambda_b)$, and a straightforward
calculation shows gives the action of the $i^*\circ \mathscr{U}_{\hat{w}}$
on this space.
\begin{lemma} \label{lem:intertwiningembedding}  We have
\[
i^*\circ \mathscr{U}_{\hat{w}}: \mathrm{Ind}(\pi^b,\mathbf{s}-\Lambda_b)\rightarrow
\mathrm{Ind}(\pi^{b-1},(s-\frac{1}{2},\ldots s-\frac{1}{2})-\Lambda_{b-1}).
\]
We will usually abbreviate the latter parameter as
$\mbox{\boldmath$ (s-\frac{1}{2})$\unboldmath}-\Lambda_{b-1}$.
\end{lemma}
\textbf{Proof.}
Let $\Phi=\Phi_{\pi^b,\mathbf{s}-\Lambda_b}\in \mathrm{Ind}(\pi^b,\mathbf{s}-\Lambda_b)$.   Let both $\mathbf{r}$ and
\[
\mathbf{m}=\boldop t_1,\ldots t_{b-1},
\widetilde{t_1},\ldots \widetilde{t_{b-1}}\boldcp\in M_{a^b}
\]
where $t_i\in\mathrm{GL}_{a}$, and let $g'\in G_{a(b-1)}$.  According
to our definition,
\begin{multline*}
(i^*\circ \mathscr{U}_{\hat{w}}\Phi)(g';\mathbf{r})=
\int\limits_{\begin{matrix}{\scriptstyle X\in \mathbf{A}^{a^2(b-1)}}\\
{\scriptstyle Z\in \mathbf{A}^{\frac{a(a+1)}{2}}}\end{matrix}} 
\Phi(\hat{w}v''(X,Z)\boldop 1_a,\mathbf{m},1_a\boldcp
i(g');i(\mathbf{r}))\,\mathrm{d}X\,\mathrm{d}Z=\\
\int\limits_{X,Z} 
\Phi(\mathbf{c}(\hat{w})
\boldop 1_a,\mathbf{m}, 1_a\boldcp\; \hat{w}\;\mathbf{c}(\boldop 1_a,\mathbf{m}^{-1},1_a\boldcp)v''(X,Z)\;
i(g');i(\mathbf{r}))\,\mathrm{d}X\,\mathrm{d}Z=\\
\int\limits_{X,Z} 
\Phi(
m(\boldop t_1,\ldots t_{b-1},1_a\boldcp)\; \hat{w}\;v''(Xm,Z)\;
i(g');i(\mathbf{r}))\,\mathrm{d}X\,\mathrm{d}Z,
\end{multline*}
where we have used \eqref{eqn:conjugatingvdoubleprime}
in the last line.  As discussed after \eqref{eqn:conjugatingvdoubleprime}
the Jacobian of the change of variables $X\mapsto Xm^{-1}$ is $|\det m|^{-a}$
According to the transformation law for $\Phi_s$, 
$m(\boldop t_1,\ldots t_{b-1},1_a\boldcp)$, when moved
to the right side of the second input of $\Phi_s$, produces
a factor of $(\det t_1,\ldots ,\det t_{b-1})$ raised to the following
multi-index
\begin{multline*}
(s-(\Lambda_b)_1+(\rho_{b}^{(a)})_1,\ldots s-(\Lambda_b)_{b-1}+
(\rho_{b}^{(a)})_{b-1})=\\
\left(s-(\Lambda_{b-1})_1-\frac{1}{2}+(\rho_{b-1}^{(a)})_1+a,
\ldots s-(\Lambda_{b-1})_{b-1}-\frac{1}{2}+
(\rho_{b-1}^{(a)})_{b-1}+a \right).
\end{multline*}
The Jacobian of the change-of-variables transformation cancels
the $a$ from each coordinate, so that we conclude that
\[
(i^*\circ \mathscr{U}_{\hat{w}}\Phi)(\mathbf{m}g';\mathbf{r})=
[\delta(G_{a(b-1)},P_{a^b})(\mathbf{m})]
\mathbf{m}^{\mbox{\boldmath${\scriptstyle (s-\frac{1}{2})}$\unboldmath}-\Lambda_{b-1}},
\]
as required.
\qed

\vspace*{.3cm}
\noindent
\textbf{Decomposition of $i^*\circ \mathscr{U}_{\hat{w}}$ 
into rank-one intertwining operators.}
As usual, let $V:=V_a^{ab}$, the unipotent
radical of $Q:=P_a^{ab}$  Because
$m(w_{\sigma})\in M$, we have 
$V_{w_a}=V_{\hat{w}}$.  Thus for $g'\in G_{a(b-1)}$ we have
\begin{multline}\label{eqn:decompositionoflocalintertwining0}
i^*\circ \mathscr{U}_{\hat{w}}\Phi_{\Delta,s}(g')=
\int_{V_{\hat{w}}} \Phi_{\Delta,s}(\hat{w}vi(g'))\,\mathrm{d}v=\\
\int_{V_{\hat{w}}} \Phi_{\Delta,s}(w_{\sigma}w_{a}w_{\sigma}^{-1}
\; \mathbf{c}(m(w_{\sigma}))v \; m(w_{\sigma})i(g'))=\\
\int\limits_{\begin{matrix}{\scriptstyle X\in \mathbf{A}^{a^2(b-1)}}\,\mathrm{d}v\\
{\scriptstyle Z\in \mathbf{A}^{\frac{a(a+1)}{2}}}\end{matrix}}   \Phi_{\Delta,s}(w_{[a(b-1)+1,ab]}\; 
\left(\begin{pmatrix}
1_{a(b-1)}&&&\\
X&1_a&Z&\\
&&1_a&\\
&&X'&1_{a(b-1)}
\end{pmatrix}\right) \; m(w_{\sigma})i(g'))\,\mathrm{d}\,X\mathrm{d}Z=
\\
\int\limits_{\begin{matrix}{\scriptstyle n_a^a\in N_a^a(\mathbf{A})}\\
{\scriptstyle n^-_{ab-a,a}\in N_{ab-a,a}^-(\mathbf{A})}
\end{matrix}}  \Phi_{\Delta,s}(w_{[a(b-1)+1,ab]}\; 
\boldop 1_{a(b-1)},n_a^a, 1_{a(b-1)}
\boldcp
m( n^-_{ab-a,a})
 \; m(w_{\sigma})i(g'))\,\mathrm{d}n,
\end{multline}
cf. (1.2.22)--(1.2.24) in \cite{kudlarallisfest}.  

Define the intertwining operators
\[
\mathscr{N}_s: \mathrm{Ind}_{P_{a^b}}^{G_{ab}}
(\pi^b,\mathbf{s}-\Lambda_b)\rightarrow  \mathrm{Ind}_{P_{a^b}}^{G_{ab}}
(\pi^b,(s+\frac{1-b}{2},s+\frac{3-b}{2},\ldots, s+\frac{b-3}{2}, -s+\frac{1-b}{2})),
\]
for $\mathrm{Re}s\gg 0$, by
\[
\mathscr{N}_s \Phi_{\pi^b,\mathbf{s}-\Lambda}(g)=
\int\limits_{n_a^a\in N_a^a(\mathbf{A})} \Phi_{\pi^b,\mathbf{s}}
(w_{[a(b-1)+1,ab]}\; 
\boldop 1_{a(b-1)},n_a^a, 1_{a(b-1)}
\boldcp i(g'))\,\mathrm{d}n,
\]
and
\begin{multline*}
\mathscr{M}_s: \mathrm{Ind}_{P_{a^b}}^{G_{ab}}
(\pi^b,(s-(\Lambda_{b})_1-\frac{1}{2},\ldots, s-(\Lambda_{b})_{b-1}-\frac{1}{2}, -s-(\Lambda_{b})_1))\\
\rightarrow (-s-(\Lambda_{b})_1,\pi^b,(s-(\Lambda_{b-1})_1-\frac{1}{2},\ldots, s-(\Lambda_{b-1})_{b-1}-\frac{1}{2}, ))
\end{multline*}
by the integral over $N_{a,a(b-1)}(\mathbf{A})$.  Picking up from
the last line of 
\eqref{eqn:decompositionoflocalintertwining0}, we have
\begin{multline}\label{eqn:decompositionoflocalintertwining}
i^*\circ \mathscr{U}_{\hat{w}}(\mathscr{N}_s\Phi_{\Delta,s})(g')=
\int\limits_{ n^-_{ab-a,a}\in U_{ab-a,a}^-(\mathbf{A})}  \Phi_{\Delta,s}(
m( n^-_{ab-a,a})
 \; m(w_{\sigma})i(g'))\,\mathrm{d}n=\\
 \int\limits_{ n\in U_{a,a(b-1)}(\mathbf{A})} \mathscr{N}_s \Phi_{\Delta,s}(
m(w_{\sigma})m(n) i(g'))\;\mathrm{d}n=\\
i^*\circ \mathscr{M}_s\circ \mathscr{N}_s \Phi_{\Delta,s},
\end{multline}
cf. \cite{jiangmemoirs},
the displayed equations immediately preceding Lemma 3.2.2.1.

\begin{prop}\label{prop:lambdaprop}
Let $\mathscr{N}_s$, $\mathscr{M}_s$, and $i^*\mathscr{U}_{\hat{w}}$
be as above.  
\begin{itemize}
\item[(a)]  The intertwining operator $\mathscr{N}_s$ has normalizing factor,
\[
r_{\mathscr{N},b}(s)=\frac{L\left(s+\frac{b-1}{2},\tau\right)  
L\left(2s+b-1,\tau,\wedge^2\right)}{L\left(s+\frac{b+1}{2},\tau\right)  
L\left(2s+b,\tau,\wedge^2 \right)},
\]
up to $\epsilon$ factors (which have neither zeros nor poles).
\item[(b)]  The intertwining operator of $\mathscr{M}_s$ has normalizing
factor
\[
r_{\mathscr{M},b}(s)=\prod_{i=1}^{b-1} \frac{L(2s+i-1,\tau\otimes\tau)}{L(2s+i,\tau\otimes \tau)},
\]
up to $\epsilon$ factors (which have neither zeros nor poles).
\item[(c)]  The poles of the local intertwining opeator $i^*\mathscr{U}_{\hat{w}}$ are contained in those of the product
\[
\lambda_b(s)=r_{\mathscr{N},b}(s)r_{\mathscr{M},b}(s)=
L\left(2s+b-1,\tau,\wedge^2\right)\cdot \prod_{i=1}^{b-1} 
L(2s+i-1,\tau\otimes\tau).
\]
\end{itemize}
\end{prop}
\textbf{Proof.}
The analytic properties of $\mathscr{N}_s$ are unchanged if we restrict
to the embedded reductive subgroup $i^{ab}_a(G_a)$.  Then
$\mathscr{N}_s$ becomes the standard intertwining operator
$N_s$
associated to the nontrivial element $w\in W(G_a,P_a)$,
\[
N_s: \mathrm{Ind}_{P_a}^{G_a}\left(\tau,s+\frac{b+1}{2}\right)\rightarrow
\mathrm{Ind}_{P_a}^{G_a}\left(\tau,-s-\frac{b+1}{2}\right).
\]
As we have already mentioned, it is well-known that the local
normalizing factor of this standard rank-one intertwining operator
is the ratio given in (a), up to the $\epsilon$-factors.

The analytic properties of $\mathscr{M}_s$ are unchanged
if we restrict to the embedded reductive subgroup 
$M=m(\mathrm{GL}_{ab})$.
Then $\mathscr{M}_s$ becomes the standard $\mathrm{GL}_{ab}$
intertwining operator $M_{w_{\sigma}}$
associated to the element $w_{\sigma}\in W(P_{a,ab-a})$, which is
the subset of the Weyl
group of $\mathrm{GL}_{ab}$
defined as in the second paragraph of \S I.1.7 of \cite{mwbook}.
In particular, if we define $M_{w\sigma}$ by the same integral
as in the second line of \eqref{eqn:decompositionoflocalintertwining}, then
\begin{multline}\label{eqn:romanMsubwsetup}
M_{w_{\sigma}}: \mathrm{Ind}_{P_{a,a(b-1)}}^{\mathrm{GL}_{ab}}(
\tau
|\cdot|^{s-(\Lambda_{b})_1}\otimes \pi^{b-1}|\cdot|^{\mathbf{s}-\Lambda_{b-1}-\frac{1}{2}}|\cdot|^{-s-(\Lambda_{b})_1})\rightarrow\\
\mathrm{Ind}_{P_{a(b-1),a}}(\tau|\cdot|^{-s-(\Lambda_{b})_1}\otimes
\pi^{b-1}|\cdot|^{\mathbf{s}-\Lambda_{b-1}-\frac{1}{2}})
\end{multline}
It is more or less standard that the intertwining operator for this
operator is the ratio of $L$-functions given in (b)
along with an $\epsilon$-factor which does not matter
for our purposes.
For completeness, we briefly recall, in the paragraphs
following this proof how one reduces the calculation
of this local normalizing factor to that of an intertwining
operator on spaces induced from $P_{a^b}$, which
can be evaluated with the standard formula of Gindikin-Karpelevic.

Part (c) follows by multiplying the numerators of (a) and (b).
\qed

\begin{cor}\label{cor:localintertwiningopholomorphy}  The local
intertwining operator 
\[
i^*\circ \mathscr{U}_{\hat{w},v}: \mathrm{Ind}_{P}^{G_{ab}}
(\Delta(\tau,b),s)_v
\rightarrow \mathrm{Ind}_P^{G_{a(b-1)}}
\left(\Delta(\tau,b-1),s-\frac{1}{2}\right)_v
\]
is holmorphic in the right half-plane $\mathrm{Re}s>0$.
\end{cor}

\vspace*{.3cm}
\noindent
\textbf{Application of the Formula of Gindikin-Karpelevic}
The treatment of the normalization for the
operator in \eqref{eqn:romanMsubwsetup}
should be compared to the parallel development in
\S3.2.4 of \cite{jiangmemoirs}.  We first embed 
the representation on the left side of \eqref{eqn:romanMsubwsetup}
into $\mathrm{Ind}_{P_{a^b}}^{\mathrm{GL}_{ab}}(\pi^b,\mathbf{t}_s)$,
where the parameter
\begin{equation}\label{eqn:sprimeparameterchoice}
\mathbf{t}_s:=(s-(\Lambda_{b})_1,
\ldots, s-(\Lambda_{b})_{b-1},-s-(\Lambda_{b})_1).
\end{equation}
Define, by the same integral formula, an intertwining operator, extending $M_{w_{\sigma}}$ to this induced space, and denote the extension by
\[
M_{w_{\sigma}}^*: \mathrm{Ind}(\pi^b,\mathbf{t}_s)\rightarrow
\mathrm{Ind}(\pi^b,\sigma\mathbf{t}_s).
\]

We recall the specific form, due to Piatetski-Shapiro
and Rallis in \cite{gelbartpsrallis}, 
of the Formula of Gindikin-Karpelevic
that we will need to apply in this case.  Let $\Omega$
be the distinguished set of coset representatives
for $W_{M_{a^b}}\backslash M_{\mathrm{GL}_{ab}}$
obtained by choosing the unique element of minimal length
in each coset.  Then by Casselman \cite{casselman}, 
on has the following description of $\Omega$:
\[
\Omega=\{w\in W_{\mathrm{GL}_{ab}}\;|\; w^{-1}\Phi^+_{M_{a^b}}
\subset \Phi^+_{\mathrm{GL}_{ab}}\}.
\]
Further, for $w\in\Omega$ one has the equality
\[
w (U_{a,ab-a})_w= \text{product of}\; N_{\alpha}\;
\text{for}\; \alpha\in \Phi_{\mathrm{GL}_{ab}}\;
\alpha>0,\, w^{-1}\alpha<0.
\]
where $N_{\alpha}$ is the one-parameter subgroup associated
to the root $\alpha$.
Then the normalizing factor of $M_w$ is 
\begin{equation}\label{eqn:gindkarpformulaapp}
r(w,\mathbf{t}_s)=\prod_{\alpha>0,w_{\sigma}\alpha<0}
\frac{L(\langle\alpha,\mathbf{t}_s\rangle, \alpha^{\vee}\circ
\pi^b)}{L(\langle\alpha^{\vee},\mathbf{t}_s\rangle+1,\alpha^{\vee}\circ\pi^b)}.
\end{equation}
Here $L(s,\alpha^{\vee}\circ\pi^b)$ is $L(s,\tau\otimes \tilde{\tau})=
L(s,\tau\otimes \tau)$ because $\alpha=e_i-e_j$, and $\tau$
is self-dual.  It is not difficult to see that the set of 
$\alpha\in\Phi^+_{\mathrm{GL}_{ab}}$ whose
sign is reversed by $w_{\sigma}^{-1}$ is precisely
the set of $\alpha$ of the form
\[
\alpha_{i,b}=f_i-f_b\;\text{for}\; 1\leq i\leq b-1.
\]
These are the roots $\alpha$ over which the product 
\eqref{eqn:gindkarpformulaapp} ranges.
Further for each $\alpha_{i,b}$, we have
\[
\langle\alpha_{i,b},\mathbf{t}_s\rangle = (\mathbf{t}_s)_i-(\mathbf{t}_s)_b.
\]
Thus for $\mathbf{t}_s$ as in \eqref{eqn:sprimeparameterchoice}
we obtain the normalizing factor $r(w_{\sigma},\mathbf{t}_s)$
precisely of the form given in part (b) of Proposition
\ref{prop:lambdaprop}.  This completes the proof
of the formula of part (b) of Proposition \ref{prop:lambdaprop}.

\section{Principal Normalized Constant term}
\label{sec:principalnormconst}
Let $w$ be the non-trivial element of the Weyl
group of $G$ with respect to the Siegel maximal
parabolic $P$.
The Functional Equation of Theorem IV.1.10 of 
\cite{mwbook}, applied to
$E(P_{a^b},\pi^b)$ implies by \eqref{eqn:rescuspdatarel}
that $E(P,\Delta^b,s)$ satsifies the functional equation
\begin{equation}\label{eqn:functequnnormform}
E(P,\Delta^b,s)=E(P,M(w,\Delta^b)\Delta^b,-s),
\end{equation}
Here,
\[
M(w,\Delta^b): \mathrm{Ind}_P^G(\Delta^b,s)\rightarrow
\mathrm{Ind}_P^G(w\Delta^b,-s)
\]
is the $G_{ab}$-intertwining operator defined
for sufficiently large values of $\mathrm{Re}\,s$ by
the integral
\[
M(w,\Delta^b)f_s(g)=\int_{U_{w}(\mathbf{A})}
f_s(wug)\,\mathrm{d}u.
\]
As usual, we have the Eulerian factorization
$M(w,\Delta^b,s)=\otimes_v M_v(w,\Delta^b,s)$
In order to produce a normalized form of the Eisenstein
series whose functional equation will display
a symmetry of the poles about the imaginary axis,
we wish to explicitly calculate the normalizing
factor of $M(w,\Delta^b,s)$: this is a certain Eulerian
ratio of products of L-functions
\[
a_b(s)=\bigotimes_v a_{b,v}(s)\;\text{and}\;
b_b(s)=\bigotimes_v b_{b,v}(s),\;\text{without common factors}.
\]
The ratio $a_{b,v}(s)/b_{b,v}(s)$ is further characterized by
the properties that
\[
\frac{1}{a_{b,v}(s)}M_v(w,\Delta^b,s)
\]
can be continued to a holomorphic function of $s$,
and so that
\[
E^*(P,\Delta^b,s):=b_{b,S}(s)E(P,\Delta^b,s)
\]
for $S$ the set of ramified places of $\tau$, 
satisfies the `normalized' functional equation
\begin{equation}\label{eqn:functeqnormform}
E^*(P,\Delta^b,s)=E^*(P,M^*(w,\Delta^b,s)\Delta^b,-s)
\end{equation}
 with both
sections holomorphic.  

The general theory of $L$-functions says
that we can calculate $a_{b,v}(s)/b_{b,v}(s)$ as the constant
by which $M_v(s)$ acts
on the spaces spanned by normalized spherical
sections.  Therefore, we can use the calculations in
\S\ref{subsec:sphericalsections}.  See the formulas
at the end of \S\ref{subsec:sphericalsections}
for the general expressions for $a_{b,v}$, $b_{b,v}$.

\noindent \textbf{Constant term in case $b=1$.}
This amounts to a short exposition of the first
part of \S3 of \cite{kimsp4paper}.  From the formula
of \cite{mwbook} for the constant terms of cuspidal
data Eisenstein series, it follows 
that the un-normalized Eisenstein
series has constant term
\[
E^{a\cdot 1}_Q(\tau,f_{\tau,s},s)=f_{\tau,s}+M(w,\tau,s)f_{\tau,s}.
\]
Our inductive formula \eqref{eqn:nonnormalizedconstantterm} indeed reduces
to this statement since intertwining operator
$\mathscr{U}_{\hat{w}_{a}}$ reduces to $M(w,\tau,s)$
and the `Eisenstein series of lower rank' on the right-hand side
are degenerate.

Let us drop the ``prime" from our notation
for the exceptional places, so that our new set
of places $S\in\Omega(k)$ is the old $S'\in\Omega(k)$.

Let the section $f_{\tau,s}$ be as in \eqref{eqn:stsection}
and recall that we have already computed
$a_{1,v}(\tau,s)$, $b_{1,v}(\tau,s)$ for $v\notin S$.
 We extend this definition to $v\in S$
by setting $a_{1,v}$ and $b_{1,v}$ to be the products
of local factors of the corresponding $L$-functions.

Then we define
\[
\tilde{f}_{\tau,s}:=\left[\bigotimes_{v\in S}\frac{1}{a_{1,v}(s)}M(w,\tau,s)f_{\tau,s,v}\right]\otimes \left[ \bigotimes_{v\notin S} \tilde{f}_{\tau,s,v}^0 \right].
\]
The important point is that $\tilde{f}_{\tau,s}$
is a \textit{holomorphic}
section of the global induced space $\mathrm{Ind}_P^G(\tau,-s)$
It follows from \eqref{eqn:normalizingfactorsbeq1case}
that
\[
M(w,\tau,s)f_{\tau,s}=\frac{a_1(s)}{b_{1,S}(s)}\tilde{f}_{\tau,s}=
\frac{L(s,\tau)L(2s,\tau,\wedge^2)}
{L_S(s+1,\tau)L_S(2s+1,\tau,\wedge^2)}\tilde{f}_{\tau,s}.
\]
Thus, the constant term of the normalized series is
\begin{equation}\label{eqn:bequal1normalizedconst}
E^{a\cdot 1,*}(\tau,f_{\tau,s},s)=L_S(s+1,\tau)L_S(2s+1,\tau,\wedge^2)
f_{\tau,s}+L(s,\tau)L(2s,\tau,\wedge^2)\tilde{f}_{\tau,s}.
\end{equation}

\vspace*{.3cm}
\noindent \textbf{Constant term in case $b\geq 2$.}

As a result of the formula in \S\ref{subsec:sphericalsections}, 
we have
\begin{multline}\label{eqn:firsttermnormalizingratio}
\frac{b_{b,S}(s)}{b_{b-1,S}(s+\frac{\scriptstyle 1}{\scriptstyle 2})}
=\frac{\scriptstyle \prod\limits_{i=1}^{\lceil \frac{b}{2}\rceil}
L_S(2s+b+2-2i,\tau,\wedge^2)\prod\limits_{i=1}^{\lfloor \frac{b}{2}\rfloor}
L_S(2s+1+b-2i,\tau,\vee^2)L_S\left(s+\frac{b+1}{2},\tau\right)}
{\scriptstyle \scriptstyle \prod\limits_{i=1}^{\lceil \frac{b-1}{2}\rceil}
L_S(2s+b+2-2i,\tau,\wedge^2)\prod\limits_{i=1}^{\lfloor \frac{b-1}{2}\rfloor}
L_S(2s+1+b-2i,\tau,\vee^2)L_S\left(s+\frac{b+1}{2},\tau\right)}=\\
\begin{Bmatrix}L_S(2s+1,\tau,\vee^2)\\
L_S(2s+1,\tau,\wedge^2)
\end{Bmatrix}_b,
\end{multline}
and
\begin{multline}\label{eqn:secondtermnormalizingratio}
\frac{b_{b,S}(s)}{b_{b-1,S}(s-\frac{\scriptstyle 1}{\scriptstyle 2})}
=\frac{\scriptstyle \prod\limits_{i=1}^{\lceil \frac{b}{2}\rceil}
L_S(2s+b+2-2i,\tau,\wedge^2)\prod\limits_{i=1}^{\lfloor \frac{b}{2}\rfloor}
L_S(2s+1+b-2i,\tau,\vee^2)L_S\left(s+\frac{b+1}{2},\tau\right)}
{ \scriptstyle \prod\limits_{i=1}^{\lceil \frac{b-1}{2}\rceil}
L_S(2s+b-2i,\tau,\wedge^2)\prod\limits_{i=1}^{\lfloor \frac{b-1}{2}\rfloor}
L_S(2s-1+b-2i,\tau,\vee^2)L_S\left(s+\frac{b-1}{2},\tau\right)}=\\
\begin{Bmatrix} \frac{\scriptstyle L_S(2s+b,\tau,\wedge^2)
L_S(2s+b-1,\tau,\vee^2) L(s+\frac{b+1}{2},\tau)}
{\scriptstyle L_S(2s,\tau,\wedge^2) L_S(s+\frac{b-1}{2},\tau)}\\ 
\\
\frac{L_S(2s+b-1,\tau,\vee^2)
L_S(2s+b,\tau,\wedge^2)L_S(s+\frac{b+1}{2},\tau)}
{L_S(2s,\tau,\vee^2)L_S(s+\frac{b-1}{2},\tau)}
\end{Bmatrix}_b=
\begin{Bmatrix} L_S(2s,\tau,\vee^2)\\
L_S(2s,\tau,\wedge^2)
\end{Bmatrix}_b\gamma_{b,S}^{-1}(s),
\end{multline}
where we have used the notation 
\[
\begin{Bmatrix}x\\y
\end{Bmatrix}_b = \begin{cases} x&\text{if $b$ even}\\
y&\text{if $b$ odd}
\end{cases},
\]
and in the last equality in \eqref{eqn:secondtermnormalizingratio}
we have used \eqref{eqn:gammaSpartial}.
Set
\begin{multline*}
\tilde{f}_{\Delta,s-\frac{1}{2}}=
\gamma_{b,S}^{-1}(s)i^*\circ \mathscr{U}_{\hat{w}}^{ab}
f_{\Delta,s}
=
\bigotimes_{v\in S}i^*\circ
\mathscr{U}_{w,v}^{ab}f_{\Delta,s,v}
\otimes\bigotimes_{v\notin S}
\gamma_{b,v}^{-1}(s)i^*\circ
\mathscr{U}_{w,v}^{ab}f_{\Delta,s,v}^0=\\
\bigotimes_{v\in S}i^*\circ
\mathscr{U}_{w,v}^{ab}f_{\Delta,s,v}
\otimes\bigotimes_{v\notin S}
\tilde{f}_{\Delta,s-\frac{1}{2},v}^0.
\end{multline*}

We have the following general inductive formula, part (a) of the following theorem,
 for the $Q$-constant term of $E^{ab,*}(f_\Delta,s)$ 
in terms of the normalized Eisenstein series of lower rank.  We are already able
to deduce from the formula part (a) a nice description, part (b), of the location of all
\textit{possible} poles of the Eisenstein series.
\begin{thm} \label{thm:definitive}
\begin{itemize}
\item[(a)] With all the notation as above, we have the inductive formula
between normalized Eisenstein series induced from holomorphic sections.
\begin{multline}\label{eqn:normalizedeisseriessectionform}
E_Q^{ab,*}(f_{\Delta,s},\mathrm{diag}\boldop t,h,\tilde{t}\boldcp )=
\\\begin{Bmatrix}L_S(2s+1,\tau,\vee^2)\\
L_S(2s+1,\tau,\wedge^2)
\end{Bmatrix}_b
|\det t|^{s+ab-(a+b)/2+1}\phi_{\tau}(t)\otimes E^{a(b-1),*}
(h,i^*f_{\Delta,s},s+\frac{1}{2})+\\
\begin{Bmatrix} L_S(2s,\tau,\vee^2)\\
L_S(2s,\tau,\wedge^2)
\end{Bmatrix}_b |\det t|^{-s+ab-(a+b)/2+1}\\
E^{a(b-1),*}\left(\tilde{f}_{\Delta,s-\frac{1}{2}},h,s-\frac{1}{2}\right)\otimes \phi_{\tau}(t),
\end{multline}
\item[(b)]  The set of all possible poles of $E^{ab,*}$ is contained in the set
$X_{b}$, where
\[
X_b:=\begin{cases}b/2,b/2-1,b/2-2\ldots, -b/2 & \tau\;\text{of symplectic type}\\
b/2-1/2, b/2-3/2,\ldots, -b/2+1/2 & \tau\;\text{of orthogonal type}
\end{cases}
\]
\end{itemize}
\end{thm}
Theorem \ref{thm:definitive}(a) exactly parallels the inductive formula
Proposition 1.2.1 of \cite{kudlarallisfest}, and may be seen as the generalization
of that formula to the present setting.  Part (b) follows easily by induction
from (a) and the known properties of the $L$-functions: compare p. 92 of \cite{kudlarallisfest}.

Having derived the inductive formula \eqref{eqn:normalizedeisseriessectionform},
we now turn our attention to extracting information from it.  The partial
results on nonvanishing of residues and square-integrability are only
one of the several directions that may be envisioned for future work.


\section{Residues: Cuspidal Support, Cuspidal Exponents}
\label{sec:residues}
\noindent\textbf{
Switching the constant term formula to point of view of automorphic representations}.  
Suppressing the inducing parabolic $P$ from the notation $E^{n,*}(P,\cdots)$
(which is always the Siegel parabolic in the appropriate rank as indicated
by the exponent $n$)
\eqref{eqn:bequal1normalizedconst} becomes (constant term of the normalized Eisenstein
series for $b=1$ so that $Q=P$),
\begin{multline}\label{eqn:inductiveformularepsb1}
E^{a\cdot 1,*}_{P}(\tau,s)=
L_S(s+1,\tau)L_S(2s+1,\tau,\wedge^2)\tau |\det |^{s+\left(\rho_{1}^{(a)}\right)_1}+\\
L(s,\tau)L(2s,\tau,\wedge^2)\tau |\det|^{-s+\left(\rho_{1}^{(a)}\right)_1}
\end{multline}
and \eqref{eqn:normalizedeisseriessectionform} becomes, for $b\geq 2$,
additionally, dropping the superscript $b$ on $\Delta$
but including superscripts $b-1$ on the images on the images
of $\Delta$ under the relevant intertwining operators,
\begin{multline}\label{eqn:inductiveformulareps}
E_Q^{ab,*}(\Delta,s)=\\
\begin{Bmatrix}{\scriptstyle L_S(2s+1,\tau,\vee^2)}\\
{\scriptstyle L_S(2s+1,\tau,\wedge^2)}
\end{Bmatrix}_b 
\tau|\det\cdot |^{s+\frac{1-b}{2}+\left(\rho_{b}^{(a)}\right)_1}\otimes E^{a(b-1),*}
 \left((i^*\Delta)^{b-1},s+\frac{1}{2} \right)+\\
\begin{Bmatrix} {\scriptstyle L_S(2s,\tau,\vee^2)}\\
{\scriptstyle L_S(2s,\tau,\wedge^2)}
\end{Bmatrix}_b 
E^{a(b-1),*}\left(\left(i^*\circ \mathscr{U}_{\hat{w}}^{*,ab}(s)\Delta\right)^{b-1}, s-\frac{1}{2}\right)\otimes 
\tau  |\det \cdot|^{-s+\frac{1-b}{2}+\left(\rho_{b}^{(a)}\right)_1}.
\end{multline}

\subsection{The automorphic forms \boldmath $\Phi_{i}^{b}\;$\unboldmath as residues.}
We will define the residue points $s_i^{(b)}$ as follows.
\label{subsec:formsasresidues}

\noindent
\textbf{Definition.}  Let $\tau$ be fixed as in \S\ref{subsec:generalnotation}. 
Define the positive integer points \linebreak $s^{(b)}_i\in (0,b/2]$ for $i=0,\ldots \lceil b/2\rceil-1$ depending on the type of $\tau$
by
\[
s_i^{(b)}:=\begin{cases}b/2-i&\text{when $\tau$ is symplectic}\\
b/2-i-\frac{1}{2}& \text{when $\tau$ is orthogonal}
\end{cases}
\]
We occasionally refer to the $s_i^{(b)}$ as the points of the `segment'
of (possible) poles.  In keeping with this terminology, we refer
to $s_0^{(b)}$ as the `endpoint' or `right endpoint' of the segment
and all points $s_i^{(b)}$ with $i\geq 1$, \textit{i.e.}, to the left of $s_0^{(b)}$
as the `interior points' of the segment.

We record for use in the inductive formula the elementary relations
\begin{equation}\label{eqn:sibinductiverelations}
\begin{aligned}
s_{i-1}^{(b-1)}&=&&s_{i}^{(b)}+\frac{1}{2}\\
s_{i}^{(b-1)}&=&&s_i^{(b)}-\frac{1}{2}
\end{aligned}
\end{equation}

When $b$ is clear from the context we may write more simply $s_i$.  

The purpose of this discussion is to formalize the notion
of the \textit{residue} $\Phi_i^{(b)}$ of the Eisenstein series $E^{ab,*}(\Delta,s)$
at $s=s_i^{(b)}$ in two different ways, each of which
will be useful in certain situations.  
Intuitively, we can think of the residue as the 'leading'
coefficient a power of $(s-s_i^{(b)})^{-1}$ in the Laurent expansion,
centered at $s_i^{(b)}$ of $E^{ab,*}(\Delta,s)$.  Namely, in general
for a family of automorphic forms $F(s)$ depending meromorphically
on the complex variable $s$, and $z_0\in\mathbf{C}$ fixed,
$F(s)$ has the Laurent expansion centered at $z_0$,
\[
F(s-z_0)=(\Psi_{-n}F)(s-z_0)^{-n}+(\Psi_{-n+1}F)(s-z_0)^{-n+1}+\cdots.
\]
where the $\Psi_{-i}F$ are certain automorphic forms.  
For $z_0=s_{0}^{(b)}$, and $F=E^{ab,*}$, we denote the leading
term by
\[
\Phi_i^{(b)}:=\Psi_{-n}F.
\]

Additional notation related to this situation is that
in order to denote terms which are not necessarily
the leading term, when $z_0=s_i^{(b)}$, $F=E^{ab,*}$,
we will define
\[
\Psi_{-m}^{i,(b)}=\Psi_{-m}E^{ab,*},
\]
the $-m^{\rm th}$ coefficient in the Laurent expansion around $s_i^{(b)}$.
Further when we wish to vary the base point $z_0$ we may write more generally
$\Psi^{z=z_0}_{-m}F$.  

Further, we will reserve the ``classical" notation $\mathrm{Res}_{z=z_0}F$
for the situation where the function $F$ has a \textit{simple} pole at $z_0$.
So, in the above terms, the use of the notation $\mathrm{Res}_{z=z_0}F$
implies that $-m_0$, the leading exponent in the Laurent expansion of $F$
about $z_0$ is $-1$ and it denotes $\Psi^{z=z_0}_{-1}F$ in that case.

We also wish to fit
$\Phi_i^{(b)}$ into the framework created by the discussion
of \S 5.1 of \cite{mwbook}.  So we actually say how all the concepts
and objects mentioned in the third paragraph of \textit{ibid.} appear
in this situation.  First, picking up the discussion from where
we left off in \S 2.4, it is evident that in our situation
(with the items on the left being the general notation found
in \cite{mwbook}),
\[
\begin{aligned}
\mathfrak{S}_i&=&&\{\Lambda_b+\rho_{a^b}+\mathbf{C}\mathbf{1}_b\}\\
\alpha_i&=&& 2e_b\\
\pi_0&=&&\Lambda_{b}+s_i^{(b)}\mathbf{1}_b+\rho_{b}^{(a)}=\Lambda_b+\rho_{a^b}+(s_i^{b}+\rho_{ab})\mathbf{1}_b
\in\mathfrak{S}_i,\,\text{above}\\
H_{\pi_0,\alpha_i^*}&=&&\pi_0+(\mathbf{C}^{b-1},0)\\
\mathfrak{S}_{i-1}&=&&\mathfrak{S}_i\cap H_{\pi_0,\alpha_i^*}=\{\pi_0 \}\\
\{\pi \}&=&&\{\pi_0\}\\
\epsilon_i&=&&\mathbf{1}_b,\;\text{determined up to nonzero scalar multiple}\\
P_{\pi_0,(2e_b)^*}(\pi_0\otimes z\epsilon_i)&=&&\langle z\epsilon_i,(2e_b)^* \rangle=z.
\end{aligned}
\]
Since in our situation, there is at most one hyperplane of $\mathfrak{S}_i$ (in our case $i=1$)
passing through the unique point $\pi_0\in \mathfrak{S}_0$, the condition
that $n$ is an integer such that the \textit{product} with the polynomial $P^n_{\pi_0,(2e_b)^*}\cdot A_1$
has as singular hyperplanes passing through $\pi_0$ the singular hyperplanes
of $E^{ab,*}(\Delta,s)$ except $\mathfrak{S}_0=\pi_0$, reduces to the condition
that $P^n_{\pi_0,(2e_b)^*}\cdot A_1$ is \textit{holomorphic} at $\pi_0$.  
Further, we can take as $Q_{n,\alpha_1^*}$ in this case
simply $P^n_{\pi_0,(2e_b)^*}$ itself, because of the last line
of the preceding list of formulas.  Then it is clear
that what \cite{mwbook} call $\mathrm{Res}_{1}A_1$
is in this case identical to the leading
coefficient in the Laurent expansion of $A_1$ about $s=s_i^{(b)}$,
in the sense of complex analysis.  Since 
$A_1$ in our case is by Proposition \ref{prop:residualdataeisseries}
nothing other than $E^{ab,*}(\Delta,s)$, and the one-step
residue datum $\mathrm{Res}^{\mathfrak{S}_0'}_{\Lambda+s_i}$
is identical to $\mathrm{Res}_{1}$, we can express
the residue datum in classical terms as $\Phi_i^{(b)}$.

\textbf{Definition.}  Set 
\[
\Phi_i^{(b)}=\mathrm{Res}_{1}E^{ab,*}(P,\Delta(\tau,b),s).
\]
The above discussion shows that $\Phi_i^{(b)}$ is an automorphic form because
\begin{multline*}
\Phi_i^{(b)}=\mathrm{Res}^{\mathfrak{S}_0'}_{\Lambda_b+s_i}E^{ab,*}(P,\Delta(\tau,b))=\\
\mathrm{Res}^{\mathfrak{S}_0'}_{\Lambda_b+s_i}\circ \mathrm{Res}_{\mathfrak{S}_0'}^{\mathfrak{P}}
E^{ab,*}(P_{a^b},\pi^b)(s)=\mathrm{Res}^{\mathfrak{P}}_{\Lambda_b+s_i}E^{ab,*}(P_{a^b},\pi^b),
\end{multline*}
where the notation is as in Proposition \ref{prop:residualdataeisseries}.
As the image of an automorphic form under the residue datum $\mathrm{Res}^{\mathfrak{P}}_{\Lambda_b+s_i}$,
$\Phi_i^{(b)}$ is also an automorphic form by the discussion of \S5.1 of \cite{mwbook}.


\subsection{Cuspidal Support of the automorphic forms  \boldmath $\Phi_{i}^{(b)}\;$\unboldmath. }
\label{subsec:cuspidalsupportPhis}
We are now switching from deriving  the inductive formula
\eqref{eqn:inductiveformularepsb1}--\eqref{eqn:inductiveformulareps}
for the constant term to extracting information from these formulas.
In the first, `softer' step of the analysis, Proposition
\ref{prop:cuspidalsupport} below, we will be using only the `general shape'
of the formula, whereas, in the `harder' step of the analysis,
the discussion preceding Theorem \ref{thm:main}, we will make use of the particular
parameter-shifts appearing on the right-side of the formula, as well
as the analytic properties of the normalizing $L$-function factors.

Temporarily, \textit{i.e.} for Proposition \ref{prop:cuspidalsupport} only,
we extend the definition of the $s_i^{(b)}$ and $\Phi_i^{(b)}$
to the $i$ in the entire range 
\begin{equation}\label{eqn:temprange}
\begin{cases}
0\leq i\leq b & \text{
$\tau$ of symplectic type}\\
0\leq i\leq b-1 &\text{
 $\tau$ of orthogonal type}
\end{cases}.
\end{equation}  Thus, for the following Proposition
only, $s_i^{(b)}$ is allowed to be zero or negative.  
\begin{prop}\label{prop:cuspidalsupport}
The automorphic form $\Phi_i^{(b)}$ is concentrated along
the standard parabolic $P_{a^b}\cong (\mathrm{GL}_a)^b\rtimes U_{a^b}$.
\end{prop}
\textbf{Proof.}
The proof is by induction on $b$ (and in this proof the type
of $\tau$ need not even be mentioned).

\vspace*{.3cm}
\noindent
\underline{Base Case: $b=1$.}

\vspace*{.1cm}\noindent
By cuspidality of the data, the only nonzero constant term
of $E^{a\cdot 1}(\tau,s)$ is that along the (Siegel) parabolic
$P=Q=P_{a^1}$ and is given in \eqref{eqn:inductiveformularepsb1}.
Since the constant terms
of $\Phi_0^{(1)}$ (resp. $\Phi_1^{(1)}$) equal the residues
of the corresponding constant terms of $E^{a\cdot 1}(\tau,s)$
at $s=\frac{1}{2}$ ($s=-\frac{1}{2}$), the only nonzero
standard constant term $\mathrm{CT}_{P'}\Phi_0^{(1)}$ (resp., $\mathrm{CT}_{P'}\Phi_1^{(1)}$)
can be for $P'=P=Q$.  Thus, we see that the proposition
reduces to entirely ``standard facts" in the base case. 

\vspace*{.3cm}
\noindent
\underline{Induction Step.}

\vspace*{.1cm}\noindent
\textit{Assume that the statement is proved for all $b'<b$
and all integer values of $i$
in the range \eqref{eqn:temprange} for those $b'$.}

\vspace*{.1cm}\noindent
Let $P'$ be an arbitrary standard parabolic of $G_{ab}$.  Then $P'=M'U'$
with $M\cong \mathrm{GL}_{a'}\times M''$, and 
\begin{equation}\label{eqn:mdoubleprimeobservation}
\text{$M''$ a standard
Levi component of a standard parabolic in $G_{ab-a'}$}
\end{equation}
 and $a'$
a positive integer.  In concrete matrix form, we may write
\begin{equation}\label{eqn:Mprimefactorization}
M'=\begin{pmatrix}\mathrm{GL}_{a'}&&\\
&M''&\\
&&\widetilde{\mathrm{GL}_{a'}}
\end{pmatrix}.
\end{equation}
Consider the constant term
\[
\mathrm{CT}_{P'}\Phi_{i}^{(b)}=\mathrm{CT}_{P'}\mathrm{Res}^{\mathfrak{P}}_{\Lambda_b+s_i^{(b)}}
E^{ab,*}(P_{a^b},\lambda\cdot\pi^b_0).
\]
Because the integration involved in taking the constant term is over a compact set, and residue
taking can be described integration over a closed path $\gamma$ (see \eqref{eqn:residuetointegral}),
we may interchange the constant term and residue operation, to obtain
\begin{equation}\label{eqn:PprimeconsttermofcapitalPhi}
\mathrm{CT}_{P'}\Phi_i^{(b)}=\mathrm{Res}^{\mathfrak{P}}_{\Lambda+s_i^{(b)}}\mathrm{CT}_{P'}
E^{ab,*}(P_{a^b},\lambda\cdot \pi^b).
\end{equation}
Now suppose that $a'$ is not a multiple of $a$.  Then
by the well-known formula for the constant term of a cuspidal data
Eisenstein series in \S II.1.7 of \cite{mwbook}, we have
\[
\mathrm{CT}_{P'}E^{ab,*}(P_{a^b},\lambda\cdot \pi)=0,
\] 
because $W(M_{a^b},M')=\emptyset$.  Thus we may assume that $a'$
is a multiple of $a$, $a'=ka$ say.  Suppose that $k>1$.  Then the same
formula in \S II.1.7 of \cite{mwbook} implies that each term of \eqref{eqn:PprimeconsttermofcapitalPhi}
contains a factor of the form $E^{\mathrm{GL}_{ak}}(P_{a^k},\lambda_w\cdot \tau^{\otimes k})$.  Here $\lambda_w$
is some parameter depending on $\lambda$ and $w\in W(M_{a^b},M')$
(The exact value of $\lambda_w$
can be determine from an easy computation, but does not interest us).  This is because a cusp form $M'$ is a tensor product
 of a cusp forms
 \[
 \phi_1\otimes\phi_2\;\text{where}\; \phi_1\in A_0(\mathrm{GL}_{ab}(k)
 \backslash\mathrm{GL}(\mathbf{A})),\; \phi_2\in A_0(M'(k)\backslash
 M''(\mathbf{A})).
 \]
 As is well known, see \text{e.g.}, Proposition IV.1.9(iib) of \cite{mwbook}, 
 \begin{multline*}
 E^{\mathrm{GL}_{ak}}(P_{a^k},\lambda_w\cdot \tau^{\otimes k})\perp A_0(\mathrm{GL}_{ak}(k)\backslash
 \mathrm{GL}_{ak}(\mathbf{A})), \\ \text{equivalently}, \,  
 E^{\mathrm{GL}_{ak}}(P_{a^k},\lambda_w\tau^{\otimes k})^{\rm cusp}=0.
 \end{multline*}
Thus,
\[
\mathrm{CT}_{P'}\Phi_i^{(b)}\perp A_0(M'(k)\backslash M'(\mathbf{A})),\;\text{equivalently}\,
{\mathrm{CT}_{P'}\Phi_i^{(b)}\,}^{\rm cusp}=0.
\]
Thus, unless $k=1$, we have that the cuspidal support of $\Phi_i^{(b)}$
along $P'$ is zero.  So we have $a'=a$ and $M'=\mathrm{GL}_a\times M''$.  
By \eqref{eqn:mdoubleprimeobservation}, then $M'$ is standard Levi
component in the parabolic $Q$ of $G$,so  we can write
the constant term along $M'$ as a composition of constant terms
\[
\mathrm{CT}_{M'}=\mathrm{CT}^{G_{a(b-1)}}_{M''}\circ \mathrm{CT}_{Q},
\]
which is just a decomposition of the integral defining the constant
term along $\mathrm{CT}_{M'}$.
According to our inductive formula, each term of $\mathrm{CT}_QE^{ab,*}(P_{a^b},\lambda\cdot\pi^b_0)$
factors as a tensor product
of a cusp form on the factor $\mathrm{GL}_a$ 
times a factor of the form $E^{a(b-1),*}(P_{a^b-1},\lambda'\cdot\pi^{b-1}_0)$,
on the embedded factor $\mathrm{G}_{a(b-1)}$.
By our assumption that the
proposition is known for $b'<b$, in particular for for $b'=b-1$, 
the cuspidal component of the factor $E^{a(b-1),*}(P_{a^b-1},\lambda'\cdot\pi^{b-1}_0)$ 
is zero, except along the parabolic $M''=M_{a^{b-1}}$.  So, unless $M''=M_{a^{b-1}}$,
\textit{i.e.}, unless $P'=P_{a^b}$ the cuspidal component of $\Phi_{i}^{(b)}$ along $P'$
is zero.  This completes the proof of the Proposition.
\qed

\subsection{Partial results on Cuspidal Exponents of the \boldmath $\Phi_i^{(b)}$\unboldmath}
\label{subsec:partialresults}
The reason that only partial results are available at the moment is that only
in certain cases do we have full knowledge of the analytic properties of the $E^{ab,*}(\Delta,s)$
at $s=0$.  Since it is part of the general theory of the construction of the spectrum
from discrete data (used \textit{e.g.} in Chapter VI of \cite{mwbook})
that all residual data Eisenstein series are holomorphic
on the entire unitary axis, which because of our normalization is the imaginary
axis, we do know that $E^{ab,*}(\Delta,s)$ is holomorphic
at the origin.  Therefore, what we mean by `analytic properties', is what is the
order of the \textit{zero} (if any) of the function $E^{ab,*}(\Delta,s)$ at the origin.
This is also the reason we go back to considering only the $s_i^{(b)}$
in the segment of points to the right of the imaginary axis.

For the following discussion, we also use this notation
for a character $\chi_{0}^{(b)}\in\mathrm{Re}X^{G}_{M^{a^b}}$, defined 
with respect to the coordinate system $\{f_i\}_{i=1}^b$ on $\mathfrak{a}_{M_{a^b}}^*$,
introduced at \eqref{eqn:restrictedrootsexplicit}, by
\[
\chi_0^{(b)}:= \rho_{b}^{(a)}+\begin{cases} -\left(\frac{2b-1}{2},\frac{2b-3}{2},\ldots, \frac{1}{2}\right),&
\tau \;\text{symplectic}\\
-\left(b-1,b-2, \ldots, 0\right)&\tau\;\text{orthogonal}
\end{cases}
\]
Then we define a particular element $\pi^{(b)}_0\in \mathfrak{P}_{a^b}$ by the conditions,
\begin{equation}\label{eqn:cuspidalexponenti0}\mathrm{Im}\pi^{(b)}_0=\tau^{\otimes b},\quad
\mathrm{Re}\pi^{(b)}_{0}=\chi_0^{(b)}.
\end{equation}
The character $\chi_0^{(b)}$ and representation $\pi^{(b)}_0$ will occur repeatedly in the description
of the cuspidal support of the residues $\Phi_i^{(b)}$.

More generally, we define for $i>0$ the character
\[
\chi_i^{(b)}:=\rho_{b}^{(a)}+\begin{cases}
-\left(\frac{2i-1}{2},\frac{2i-3}{2},\ldots, \frac{1}{2},\frac{2(b-i)-1}{2},\frac{2(b-i)-3}{2},\ldots,\frac{1}{2}\right)
&
\tau \;\text{symplectic}\\
-\left(i,i-1, \ldots,1,b-1-i,b-2-i,\ldots, 0\right)&\tau\;\text{orthogonal}
\end{cases}
\]
and the particular element $\pi^{(b)}_0\in \mathfrak{P}_{a^b}$ 
by the condition analogous to \eqref{eqn:cuspidalexponenti0}, but with
$\chi_i^{(b)}$ in place of $\chi_0^{(b)}$.  In the present paper,
we will make use \textit{primarily} of $\pi_0^{(b)}$ and $\pi_1^{(b)}$, but all
the $\pi_i^{(b)}$ will play an equal role in the projected sequel. 

\begin{itemize}
\item[] \underline{Case b=1}
\end{itemize}
\begin{itemize}
\item $\tau$ of symplectic type.
\end{itemize}
There is only the case of $i=0$, so that $s_0=\frac{1}{2}$.  
We calculate from \eqref{eqn:inductiveformularepsb1}
\[
\mathrm{CT}_P\Phi^{(1)}_0=\mathrm{Res}_{s=1/2}(L(s,\tau)L(2s,\tau,\wedge^2))\tau
|\det t|^{\left(\rho_{a^1}\right)_1-\frac{1}{2}},
\]
since the first term is holomorphic at $s=\frac{1}{2}$.  Thus $\mathrm{CT}_P\Phi^{(1)}_0$ generates
the \textit{cuspidal} irreducible representation
\[
\chi_0^{(1)}\cdot \tau:=\pi_0^{(1)},
\]
It is clear that a cuspidal irreducible representation has cuspidal support
equal to the singleton set consisting of itself, and has cuspidal exponent
equal to to its real part.
\begin{itemize}
\item $\tau$ of orthogonal type.  
\end{itemize}
The discussion is vacuous in this case.  
(It is easy to see that in this case the Eisenstein series is holomorphic at $\frac{1}{2}$.))

\begin{itemize}
\item \underline{Case $b= 2$}
\end{itemize}

In the $b=2$ case, \eqref{eqn:inductiveformulareps} says that
\begin{multline}\label{eqn:inductiveformularepsb2}
E_Q^{a\cdot 2,*}(\Delta^2,s)=\\
 L_S(2s+1,\tau,\vee^2) 
\tau|\det\cdot |^{s-\frac{1}{2}+\left(\rho_{2}^{(a)}\right)_1}\otimes E^{a\cdot 1,*}
 \left((i^*\Delta)^1, s+\frac{1}{2}\right)+\\
 L_S(2s,\tau,\vee^2)
E^{a\cdot1, *}\left(\left( i^*\circ \mathscr{U}_{\hat{w}}^{*,a\cdot 1}(s)\Delta\right)^1,s-\frac{1}{2}\right)\otimes 
\tau  |\det \cdot|^{-s-\frac{1}{2}+\left(\rho_{2}^{(a)}\right)_1}.
\end{multline}

\begin{itemize}
\item Case $b=2$, $\tau$ of symplectic type.
\end{itemize}

In this case $s_0^{(2)}=1$.   Again, therefore, we only have
the `endpoint', no `interior' points.  First
note that $s_0^{(2)}+\frac{1}{2}=\frac{3}{2}$ lies to the right of the right-most
pole of $E^{a\cdot 1,*}$, so that the first term is holomorphic at $s_0^{(2)}$,
and therefore (also using \eqref{eqn:sibinductiverelations}),
\begin{equation}\label{eqn:caseb2tausymplectic}
\mathrm{CT}_Q\Phi^{(2)}_0= L_S(2,\tau,\vee^2)
\mathrm{Res}_{s=s_1^{(0)}}E^{a\cdot1, *}\left(\left( i^*\circ \mathscr{U}_{\hat{w}}^{*,a\cdot 1}(s)\Delta
\right)^1,s-\frac{1}{2}\right)\otimes 
\tau  |\det \cdot|^{-\frac{3}{2}+\left(\rho_{2}^{(a)}\right)_1}.
\end{equation}
Since by \eqref{eqn:Uintertwiningoperatordomainrange}--\eqref{eqn:normalizedintertwiningdefn}
\[ i^*\circ \mathscr{U}_{\hat{w}}^{*,a\cdot 1}(s)\left(\mathrm{Ind}(\Delta,s)\right)
 \cong \mathrm{Ind}\left(\tau,s-\frac{1}{2}\right),
\]
the factor $\mathrm{Res}_{s=s_0^{(1)}}E^{a\cdot1, *}\left(\left( i^*\circ \mathscr{U}_{\hat{w}}^{*,a\cdot 1}(s)\Delta
\right)^1,s-\frac{1}{2}\right)$ in \eqref{eqn:caseb2tausymplectic}
has the same cuspidal support as $\Phi_0^{(1)}$, namely
$\pi_0^{(1)}$.  We can therefore read off from \eqref{eqn:caseb2tausymplectic} 
that $\Phi^{(2)}_0$ has cuspidal support $(-\frac{3}{2},\chi_0^{(1)})\cdot \tau^{\otimes 2}:=\pi_0^{(2)}$.
\begin{itemize}
\item Case $b=2$, $\tau$ orthogonal.
\end{itemize}
Now assume $\tau$ is orthogonal,
so that $s_0^{(2)}=\frac{1}{2}$.  Here again, the only case is the endpoint.  As is
always the case with the endpoint $s_0^{(b)}$, the first term of \eqref{eqn:inductiveformularepsb2}
is holomorphic.  Therefore, in the $\tau$-orthogonal case,
\begin{equation}\label{eqn:b2tauorthogonal}
\mathrm{CT}_Q\Phi^{(2)}_0= [\mathrm{Res}_{s=1 }L_S(s,\tau,\vee^2)]
E^{a\cdot1, *}\left(\left( i^*\circ \mathscr{U}_{\hat{w}}^{*,a\cdot 1}\left(\frac{1}{2}\right)\Delta\right)^1,0 \right)\otimes 
\tau  |\det \cdot|^{-1+\left(\rho_{a^2}\right)_1}
\end{equation}
From Proposition IV.1.11(b) of \cite{mwbook}, the \textit{unnormalized} Eisenstein
$E^{a\cdot1}(s)$
is holomorphic at $s=0$ since it is a cuspidal-data Eisenstein series.  
Further the normalizing factor $b^{1}_S$ contains only standard
and exterior square L-functions, so it is holomorphic for $s\geq 0$, under the assumption that
$\tau$ is orthogonal.  So the normalized Eisenstein series $E^{a\cdot 1,*}(s)$
is holomorphic at $0$.  Therefore,
it remains to compute the cuspidal exponents of the automorphic
form resulting from evaluating this Eisenstein series at $0$.  Since
the cuspidal exponents are the same as those of $E^{a\cdot 1, *}(\tau,0)$,
and since $\mathrm{CT}_{P_{a^2}}\Phi^{(2)}_0=
\mathrm{CT}^{Q}_{P_{a^1}}\mathrm{CT}_Q\Phi^{(2)}_0$
one obtains by substituting \eqref{eqn:inductiveformularepsb1} evaluated at $s=0$
into \eqref{eqn:b2tauorthogonal}:
\[
\mathrm{CT}_{P_{a^2}}\Phi_0^{(2)}={\scriptstyle (*)2L_S(1,\tau)L_S(1,\tau,\wedge^2)}
\left(i^*\circ \mathscr{U}_{\hat{w}}^{*,a\cdot 2}\left(\frac{1}{2}\right)\Delta\right) |\det \cdot|^{0+\left(\rho_{1}^{(a)}\right)_1}\otimes \tau |\det \cdot|^{-1+(\rho_2^{(a)})_1},
\]
where $(*)$ represens a nonzero constant.  Since $\tau$
is orthogonal, by assumption, $L_S(1,\tau)L_S(1,\tau,\wedge^2)$
is likewise a nonzero number, and therefore, we deduce that
\begin{multline*}
\mathrm{CT}_{P_{a^2}}\Phi_0^{(2)}(\boldop g_1,g_2,\tilde{g_2},\tilde{g_1}\boldcp )=\\
(*)((i^*\circ \mathscr{U}_{\hat{w}}^{*,a\cdot 2}\left(\frac{1}{2}\right))\Delta)^1(g_2)|\det g_2|^{0+(\rho_{1}^{(a)})_1}
\otimes \tau(g_1) |\det g_1|^{-1+(\rho_2^{(a)})_1},
\end{multline*}
from which we can read off (since the representation $((i^*\circ \mathscr{U}_{\hat{w}}^{*,a\cdot 2}\left(\frac{1}{2}\right))\Delta)^1$
is isomorphic to $\tau$) that the cuspidal support of $\Phi_0^{(2)}$ is $(M_{a^2},\pi_0^{(2)})$,
and in particular $\Phi_0^{(2)}$ has sole cuspidal exponent $(-1,0)=\chi_0^{(2)}$,
in the notation appropriate for the $\tau$-orthogonal case.
\begin{itemize}
\item \underline{Case $b=3$.}
\end{itemize}
In the $b=3$ case overall, the formula for the $Q$-constant term is
\begin{multline}\label{eqn:b3overall}
E_Q^{a\cdot 3,*}(\Delta^3,s)=L_S(2s+1,\tau,\wedge^2)\tau |\det \cdot|^{s-1+\left(\rho_3^{(a)}\right)_1}
\otimes E^{2,*}\left(\left(i^*\Delta\right)^2,s+\frac{1}{2}\right)+\\L(2s,\tau,\wedge^2)\tau |\det \cdot|^{-s-1+
\left(\rho_3^{(a)}\right)_1}\otimes E^{2a,*}\left(i\circ \mathscr{U}_{\hat{w}}^{*,3a}(s)\Delta)^2,s-\frac{1}{2}\right)
\end{multline}

\begin{itemize}
\item  Case $\tau$ symplectic-type, $i=0$.
\end{itemize} 
 The residue point is $s_0^{(3)}=\frac{3}{2}$.  The first term in \eqref{eqn:b3overall} is holomorphic.
 The second term is
 \[
L(3,\tau,\wedge^2)\mathrm{Res}_{s=s_0^{(2)}}E^{a\cdot 2,*}((i\circ \mathscr{U}_{\hat{w}}^{*,a\cdot 3}\Delta)^2,s)
\otimes\tau |\det \cdot|^{-\frac{5}{2}+\left(\rho_3^{(a)}\right)_1}.
 \] 
But $\mathrm{Res}_{s=s_0^{(2)}}E^{a\cdot 2,*}((i\circ \mathscr{U}_{\hat{w}}^{*,a\cdot 3}\Delta)^2,s)$
has the same cuspidal support at $\Phi_0^{(2)}$, namely $(M_{a^2},\Pi_0^{(2)})$.
Therefore, we see that $\Phi_0^{(3)}$ has cuspidal support $(M_{a^3},\Pi_0^{(3)})$.

\begin{itemize}
\item Case $\tau$ symplectic-type, $i=1$.
\end{itemize}
\textit{This is the first interior point we reach and also the first in which we
have to consider the possibility that both terms of the constant
terms may have poles.}

We have $s_1^{(3)}=\frac{1}{2}$, and at this point \eqref{eqn:b3overall}
implies 
\begin{multline}\label{eqn:b3tausympl}
\mathrm{CT}_Q\Psi^{1,(3)}_{-m}=L(2,\tau,\wedge^2)\tau |\det \cdot|^{-1/2+
\left(\rho_3^{(a)}\right)_1}\otimes \Psi_{-m}^{0,(2)}E^{a\cdot 2,*}
\left(\left(i^*\Delta\right)^2\right)\\+
\mathrm{Res}_{s=1}L(s,\tau,\wedge^2)
\Psi_{-m+1}^{z=0}E^{a\cdot 2,*}
\left(i^*\circ \mathscr{U}_{\hat{w}}^{*,3a}\left(z+\frac{1}{2}\right)\Delta^2\right)
\otimes 
\tau |\det\cdot |^{-\frac{3}{2}}.
\end{multline}
The first term is $0$ for $m>1$ and nonzero for $m=1$.  By induction
when $m=1$, the first term has the cuspidal exponent
equal to $-\frac{1}{2}$ appended to the cuspidal exponent of $\Phi_0^{(2)}$,
meaning $(-\frac{1}{2},\chi_0^{(2)}):=\chi_1^{(3)}$.

\vspace*{.3cm}\noindent
\textbf{Remark.}  The second term, whether nonzero or not, can in any case \textit{not}
cancel with the first term because it has cuspidal exponents with $-\frac{3}{2}$
in the first entry. 

To analyze the second term, the main question is what are the analytic
properties of the normalized Eisenstein series
$E^{a\cdot 2,*}\left((i^*\circ \mathscr{U}_{\hat{w}}^{*,3a}(s+\frac{1}{2})\Delta)^2,s\right)$
at $s=0$.  \textit{In this special case of $b=2$}, since
\[
b_{2,S}(\Delta,s)=L_S(2s+2,\tau,\wedge^2)L_S(2s+1,\tau,\vee^2)L_S\left(s+\frac{3}{2},\tau\right).
\]
we see that $b_{2,S}(\Delta,0)$ is a nonzero number.  (In higher rank we will see
this is \textit{not} true.)  Thus, the analytic properties of $E^{a\cdot 2,*}$ are the same
as those of $E^{a\cdot 2}$.  Thus, up to constant factor, namely $b_{2,s}(\Delta,0)$,
$\Psi_{-m+1}E^{a\cdot 2,*}$ equals $\Psi_{-m+1}^{z=0}E^{a\cdot 2}$.  By
the holomorphy of the Eisenstein series on the unitary axis, when $m>1$,
this term is $0$.  Since we have already determined that the first
term has leading term with $m=1$, the leading term of $E^{a\cdot 3,*}$
as a whole is $m=1$.  

So we can say that the second term contributes up to constant, the \textit{value}
of $E^{a\cdot 2,*}\left((i^*\circ \mathscr{U}_{\hat{w}}^{*,3a}(s+\frac{1}{2})\Delta)^2,s\right)$
at the origin.  We now determine the cuspidal exponents of
this automorphic form by taking further constant terms.
By applying first the functional equation
in the form \eqref{eqn:functeqnormform} and then  Lemma \ref{lem:krUMrelation}
to  \eqref{eqn:inductiveformularepsb2}, we obtain
\begin{multline*}
\mathrm{CT}_{Q_2}^L E^{a\cdot 2,*}(\Delta^2,s)=
L_S(2s+1,\tau,\vee^2)\tau |\det \cdot|^{s-\frac{1}{2}+(\rho_2^{(a)})_1}
E^{a\cdot 1,*}((i^*\Delta)^1,s+\frac{1}{2})+\\
L_S(2s,\tau,\vee^2)E^{a\cdot 1,*}((i^*\circ M^{2\cdot a,*}(w,s)\Delta)^1,-s+\frac{1}{2})
\otimes \tau |\det \cdot |^{-s-\frac{1}{2}+(\rho_2^{(a)})_1}.
\end{multline*}
We work further on this answer and write out the constant term
\begin{multline*}
\mathrm{CT}^L_{P_{a^2}}E^{a\cdot 2,*}(\Delta^2,s)=\\
L_S(2s+1,\tau,\vee^2)\tau |\det \cdot|^{s-\frac{1}{2}+\rho_2}L(2s+1,\tau,\wedge^2)
L(s+\frac{1}{2},\tau)i^*\Delta |\det\cdot|^{-\frac{1}{2}-s+\rho_1}+\\
L_S(2s,\tau,\vee^2)\tau |\det \cdot|^{-s-\frac{1}{2}+\rho_2}L(-2s+1,\tau,\wedge^2)
L(-s+\frac{1}{2},\tau)i^*\circ M^{2\cdot a,*}(w,s)\Delta |\det\cdot|^{-\frac{1}{2}+s+\rho_1}+\\
L_S(2s+1,\tau,\vee^2)\tau |\det \cdot|^{s-\frac{1}{2}+\rho_2}L(2s+2,\tau,\wedge^2)
L(s+\frac{3}{2},\tau)i^*\Delta |\det\cdot|^{\frac{1}{2}+s+\rho_1}+\\
L_S(2s,\tau,\vee^2)\tau |\det \cdot|^{-s-\frac{1}{2}+\rho_2}L(-2s+2,\tau,\wedge^2)
L(-s+\frac{3}{2},\tau)i^*\circ M^{2\cdot a,*}(w,s)\Delta |\det\cdot|^{\frac{1}{2}-s+\rho_1}.
\end{multline*}

\vspace*{.3cm}
\textbf{Remark.}  The general theory that the normalized
$E^{a\cdot 2,*}(\Delta^2,s)$ is holomorphic at the origin predicts that
that at $s=0$, the minus-$1$ terms in the Laurent expansions
of the first two terms in the previous equation cancel.

We read off from this result that the cuspidal
exponents of $E^{a\cdot 2,*}(\Delta^2,s)|_{s=0}$, and therefore, the second term of \eqref{eqn:b3tausympl}
are
\begin{equation}\label{eqn:shufflesofchizero1}
\left(-\frac{3}{2},-\frac{1}{2},\pm\frac{1}{2}\right).
\end{equation}
But we have not determined when the value of this Eisenstein
series is nonzero and when the cuspidal exponents of \eqref{eqn:shufflesofchizero1}
actually occur.

\vspace*{.3cm}\noindent
\textit{Remarks on vanishing of $E^{a\cdot 2,*}$}  Let $m_0\geq 0$ be the smallest integer
(known to be positive by the holmorphicity)
such that there is a section $\phi_{\Delta}\in\mathrm{Ind}(\Delta,0)$ with
$\Psi_{m_0}^{z=0}E^{a\cdot 2,*}(\Delta,\phi_{\Delta}\cdot |\det \cdot |^{s},s)$
the leading term of the Laurent expansion at $s=0$.  In other words, $m_0-1$
is the order of vanishing of the entire \textit{family}
 $E^{2a,*}(\Delta,s)$ at $s=0$.
Then for this specific value of $m_0$,
\begin{multline}\label{eqn:Azerointertwining}
\phi_{\Delta,0} \mapsto \Psi_{m_0}(\phi_{\Delta,0})\;\text{is a nontrivial $G_{a\cdot 2}$-intertwining map}\\
\mathrm{Ind}(\Delta,0)\rightarrow A(G(k)\backslash G )_{\Delta}.
\end{multline}
The vanishing  of the Eisenstein series at the origin,
for a particular section $\phi_{\Delta,s}$ is equivalent to $\phi_{\Delta,0}$
in the kernel of \eqref{eqn:Azerointertwining}.  A general
result of Tadic, Theorems 11.6 and 11.8 of \cite{tadicisrael}, implies
that the domain space $\mathrm{Ind}(\Delta,0)$ in \eqref{eqn:Azerointertwining} is reducible
at each local place if and only if $b$ is odd.  In the case at hand, $b$ is even.
By the equivalence just stated, in order to show the nonvanishing
of the Eisenstein series for \textit{all} sections, it suffices to show that
\begin{multline}\label{eqn:nonvanishingsuffassumption}
\text{There is at least one section $\phi_{\Delta,0}\in\mathrm{Ind}(\Delta,0)$}
\\ \text{
such that $E^{a\cdot 2,*}(\Delta,\phi_{\Delta,s},s)$ is nonvanishing at the origin}.
\end{multline}
In cases where the Tadic assumption does not apply, which as we will see
arise in the case $b=4$, $i=1$, $\tau$ orthogonal-type, the analysis will be more subtle
because it will have to take account of the composition
series (also described by Tadic) for the local induced representations
and one may (and expects to) have different answers for the vanishing
or nonvanishing for sections belonging to different Jordan-Holder factors.

\vspace*{.3cm}
We define the triples in \eqref{eqn:shufflesofchizero1} as the \textbf{non-trivial shuffles of $\chi_0^{(1)}$}.
This set of shuffles can be described as the set of permutations of $\chi_0^{(1)}$
where only the leading element, $-\frac{1}{2}$ changes place, and also
changes sign to $\frac{1}{2}$ when it reaches the last entry.

With this definition, we can sum up what
we know up to this point about the cuspidal exponents
of $\Phi_1^{(3)}$.  The set of cuspidal exponents of $\Phi_1^{(3)}$ 
is contained in the set consisting of $\chi_0^{(1)}$ and its nontrivial
shuffles \eqref{eqn:shufflesofchizero1}, with at least one non-trivial shuffle occurring
if and only if
\[
\phi_{\Delta,0}\mapsto \Psi_0^{z=0}E^{a\cdot 2,*}(\Delta,\phi_{\Delta}\cdot |\det \cdot|^s,s)\;\text{induces}\;
\mathrm{Ind}(\Delta,0)\stackrel{\cong}{\rightarrow}A(G(k)\backslash G(\mathbf{A}))_{\Delta},
\]
where $\cong$ indicates the mapping is an isomorphism.

\begin{itemize}
\item Case $b=3$, $\tau$ orthogonal-type.
\end{itemize} We only have to consider the endpoint, $b_0^{(3)}$,
because $b_0^{(3)}=1$ in that case.  Then the first term of \eqref{eqn:inductiveformulareps}
is easily seen to be holomorphic, while the second term evaluates to
\[
L(2,\tau,\wedge^2)\mathrm{Res}_{s=s_0^{(2)}}E^{a\cdot 2,*}(P_{a\cdot 2}, \left(\mathscr{U}_{\hat{w}}^{*,a\cdot 3}\left(1
\right)\Delta\right)^{2},s)\otimes \tau |\det \cdot |^{-2}
\]
Analogous to the case of $b=2$, $\tau$ orthogonal-type considered above,
a simple inductive argument implies that the sole cuspidal exponent in this case is $(-2,-1,0)$.

\begin{itemize}
\item \underline{Case $b=4$.}
\end{itemize}
In this case \eqref{eqn:inductiveformulareps}
says that at $s_i^{(4)}$ (interpreting evaluation as residue as appropriate),
\begin{multline}\label{eqn:b4overall}
E_Q^{a\cdot 4,*}(\Delta^4,s^{(4)}_i)=\\
 L_S(2s^{(4)}_i+1,\tau,\vee^2)
\tau|\det\cdot |^{s^{(4)}_i-\frac{3}{2}+\left(\rho_{4}^{(a)}\right)_1}\otimes E^{a\cdot 3,*}
 \left(P_{a\cdot 3},(i^*\Delta)^{3},s^{(3)}_{i-1} \right)+\\
L_S(2s^{(4)}_i,\tau,\vee^2)
E^{a\cdot 3,*}\left(\left(i^*\circ \mathscr{U}_{\hat{w}}^{*,a\cdot 3}(s^{(4)}_i)
\Delta\right)^{3}, s^{(3)}_{i}\right)\otimes 
\tau  |\det \cdot|^{-s^{(4)}_i-\frac{3}{2}+\left(\rho_{a^4}\right)_1}.
\end{multline}
We skip the endpoints because by now the cuspidal exponents in this
case are clear by now.
\begin{itemize}
\item Case $b=4$, $\tau$ symplectic, $i=1$.
\end{itemize}
In this case $s^{(4)}_1=1$.  In that case
\eqref{eqn:b4overall} becomes
\begin{multline}
\mathrm{CT}_{Q}\Phi_1^{(4)}=
L(3,\tau,\vee^2)\tau |\det\cdot|^{-\frac{1}{2}+(\rho_4^{(a)})_1}\otimes
\mathrm{Res}_{s=s_0^{(3)}}E^{a\cdot 3,*}((i^*\Delta)^3,s)+\\
L(2,\tau,\vee^2)\mathrm{Res}_{s=s_1^{(3)}}E^{a\cdot 3,*}
((i^*\circ \mathscr{U}_{\hat{w}}^{a\cdot 4,*}(1)\Delta)^3,s)\otimes\tau|\det\cdot|^{-\frac{5}{2}+(\rho_4^{(a)})_1}
\end{multline}
The \textit{first} term has cuspidal exponents
\[
\left(-\frac{1}{2},-\frac{5}{2},-\frac{3}{2},-\frac{1}{2}\right):=\chi_1^{(4)},
\]
where we have used the determination above
of the cuspidal support of $\Phi_0^{(3)}$, and therefore,
also has cuspidal support $\pi_1^{(4)}$.

Examining the second term, and referring back to
the study of the cuspidal support of $\Phi_1^{(3)}$
in an analogous way, we see that this term
has cuspidal exponents contained in the set
$\left(-\frac{5}{2},-\frac{3}{2},-\frac{1}{2},\pm\frac{1}{2}\right)$,
which we may informally refer to as the set of ``shuffles" of
$\chi_1^{(4)}$. 
Because of the obvious non-occurence of $\chi_1^{(4)}$
in this set of shuffles, we clearly have no cancellation
of the first term by the second.  So we may
conclude that the set of cuspidal exponents
of $\Phi_1^{(4)}$ at least contains the singleton
set $\{\chi_1^{(4)}\}$, and is contained in the
set consisting of $\chi_1^{(4)}$ and the above
two ``shuffles", with the exact identity of the set being a matter
for further investigation (and depending on condition \eqref{eqn:nonvanishingsuffassumption}).

\begin{itemize}
\item Case $b=4$, $\tau$ orthogonal, $i=1$.
\end{itemize}
\textit{Note that this is the first \textbf{interior} point we encounter
in the case that $\tau$ is orthogonal.}  This should be considered in parallel
to the case $b=3$, $\tau$ symplectic, $i=1$, considered just above.

In this case, $s_1^{(4)}=\frac{1}{2}$ and the point referred to as `$s_1^{(3)}$' in \eqref{eqn:b4overall}
is actually the origin.  That is, \eqref{eqn:b4overall} implies that for any $m\in\mathbf{Z}$,
\begin{multline}\label{eqn:b4tauorth}
\mathrm{CT}_Q\Psi^{(4),1}_{-m}=
 L_S(2,\tau,\vee^2)
\tau|\det\cdot |^{-1+\left(\rho_{4}^{(a)}\right)_1}\otimes
\Phi_{-m}^{z=s_0^{(3)}}E^{a\cdot 3,*}
 \left((i^*\Delta)^{3},z \right)+\\
\left(\mathrm{Res}_{s=1}L_S(s,\tau,\vee^2)\right)\left[
\Phi_{-m+1}^{z=0}E^{a\cdot 3,*}\left(\left(i^*\circ \mathscr{U}_{\hat{w}}^{*,ab}\left(z+\frac{1}{2}\right)
\Delta\right)^{3}, z\right)\right]\otimes 
\tau  |\det \cdot|^{-2+\left(\rho_{4}^{(a)}\right)_1}.
\end{multline}
The first term is zero for $m>1$ and nonzero for $m=1$, and the leading, $-1$, term
has has cuspidal 
exponent $(-1,-2,-1,0):=\chi_1^{(4)}$.
Further, since the factor $b_{3,S}$ is explicitly, in this case given by 
\[
b_{3,S}(s)=L_S(2s+3,\tau,\wedge^2)L_S(2s+1,\tau,\wedge^2)
L_S(2s+1,\tau,\vee^2)L_S(s+2,\tau),
\]
$b_{3,S}(s)$ has a simple pole at $s=0$.  Thus the bracketed factor in 
the second term is up to constant (namely $\mathrm{Res}_{s=0}b_{2,S}(s)$),
\begin{equation}
\label{eqn:b4s1secondterm}
\Psi_{-m+2}^{z=0}E^{a\cdot 3}((i^*\circ\mathscr{U}_{\hat{w}}^{*,ab}(z+\frac{1}{2})\Delta)^3,z)
\end{equation}

We know that the corresponding unnormalized series $E^{a\cdot 3}(\cdots)$
is holomorphic at the origin.  Let $n_0$, a non-negative integer, be
the order of the zero of this Eisenstein series at the origin.  More precisely,
for a fixed section $\phi_{\Delta}$, define $n_0$ so that
\begin{multline}\label{eqn:nzeroordofzerodefn}
E^{a\cdot 3}((i^*\circ\mathscr{U}_{\hat{w}}^{*,ab}(z+\frac{1}{2})\phi_{\Delta})^3,z)|_{z=0}=0\\ \text{if and only
if $n_0>0$ and if so `vanishes to order' $n_0-1$}.
\end{multline}

Clearly if $-m+2<n_0$, then \eqref{eqn:b4s1secondterm} is zero, and is nonzero
if $m>2-n_0$, and is nonzero if $-m+2=n_0$.  A priori, we must consider
the following cases in order of increasing ``complexity":
\begin{itemize}
\item[$n_0\geq 2$]  Then for $m>0$, \textit{i.e.}, $m\geq 1$,
\eqref{eqn:b4s1secondterm} is zero and the entire second term is holomorphic.  This implies that the
singularity of $E^{ab,*}(\Delta,s)$ at $s_1^{(4)}=\frac{1}{2}$ is carried by the first term,
so it is a single pole whose residue has cuspidal exponent $\chi_1^{(4)}$.
\item [$n_0=1$]  Then the pole of the second term is simple,
so that the singularity is carried by both terms, and the cuspidal
exponents of the leading term are $\chi_1^{(4)}$ appended to thsoe of $\mathrm{Res}_{z=0}
E^{a\cdot 3,*}(\cdots,z)$.  In order to determine these exponents we must examine
 $\Psi^{z=0}_0 E^{a\cdot3,*}(\cdots, z)$, the evaluation of the normalized
 Eisenstein series, which is holomorphic by assumption, at the origin.
\item [$n_0=0$]  Then the second term gives a pole of order $2$
and the leading term of the sum equals the leading term of the second term.
In order to determine these exponents we must examine
 $\Psi^{z=0}_{-1} E^{a\cdot3,*}(\cdots, z)$.
\end{itemize}
But actually the case $n_0=1$ cannot occur because when we evaluate 
\begin{equation}\label{eqn:b4tauorthogonalb3eis}
E_{Q_3}^{a\cdot 3,*}(\Delta^3,s)|_{s=0}={\scriptstyle L(1,\tau,\wedge^2)}\tau |\det \cdot|^{-1+(\rho_{3}^{(a)})_1}
E^{a\cdot 2,*}(i^*(\mathrm{Id}+M^{a\cdot 3,*}(0)\Delta)^2,\frac{1}{2})
\end{equation}
We now refer back to the calculations in the case of $b=2,\, i=0$, $\tau$ orthogonal
type above.  Unless the condition
\[
i^*(\mathrm{Id}+M^{a\cdot 3,*}(0)\circ (i^*\circ\mathscr{U}_{\hat{w}}^{*,ab}(\frac{1}{2})\phi_{\Delta}\in
\mathrm{ker} (i^*\circ \mathscr{U}_{\hat{w}}^{*,a\cdot 2}\left(\frac{1}{2}\right))
\]
is satisfied, the right-hand side of \eqref{eqn:b4tauorthogonalb3eis} has a pole,
and we are therefore, really in the case $n_0=0$.  Because Tadic's irreducibility
result implies that $\mathrm{Ind}(\Delta^2,1/2)$ is reducible, certain
sections $\phi_{\Delta}$ may satisfy the above condition, but then we
are automatically actually in the case $n_0\geq 2$.  

In the case $n_0=0$, a when we must consider $z\cdot  E^{a\cdot3,*}(\cdots, z)|_{z=0}$,
a further calculation gives
\begin{multline}\label{eqn:b4tauorthn2}
\mathrm{CT}(G_{a\cdot 3},P_{a,a}^{a\cdot 3})s\cdot E^{a\cdot 3,*}(\Delta^3,s)=
\mathrm{CT}^L_{Q_2}s\cdot E^{a\cdot 3,*}_{Q}(\Delta^3,s)=\\
{\scriptstyle L(2s+1,\tau,\wedge^2)\tau |\det \cdot|^{s-1+(\rho_3)^1}}
\left(
 L(2s+2,\tau,\vee^2)|\det \cdot|^{s+(\rho_2)_1}\underline{s\cdot \left(E^{1,*}((i^*\Delta)^1,1+s)\right)}+\right.\\ \left.
\tau |\det \cdot|^{-s-1+(\rho_2)_1}\underline{s\cdot \left(L(2s+1,\tau,\vee^2)E^{1,*}(i^*\circ \mathscr{U}_{\hat{w}}^{a\cdot 2,*}
(i^*\circ\Delta)^1,s)\right)}
\right)\\
+{\scriptstyle L(2s,\tau,\wedge^2)\tau |\det\cdot|^{-s-1+(\rho_3^{(a)})_1}}
\left(
L(2-2s,\tau,\vee^2)\tau'\, |\det\cdot|^{-s+(\rho_2)_1}\underline{s\cdot\left(
E^{1,*}((i^*\circ M^{a\cdot 3}(s)\Delta)^1,1-s)\right)}+\right.\\ \left.
\tau'\, |\det \cdot|^{s-1+(\rho_2)_1}
\underline{s\cdot \left( L(1-2s,\tau,\vee^2) E^{1,*}(i^*\circ\mathscr{U}\circ i^* M^{a\cdot 3,*}(s)\Delta,-s)\right)}
\right),
\end{multline}
where $\tau'$ is a certain image of $\Delta^3$ under an intertwining
operator, isomorphic to $\tau$.  In each of the terms of \eqref{eqn:b4tauorthn2},
the singularities must come from the factor consisting of a product of symmetric-square L-functions
and normalized Eisenstein series, rather than the factor consisting of powers of determinants.
Therefore, we can use the following general method to evaluate
the cuspidal components of each term.  For $f(s)g(s)$
a function with a simple at $s=0$, expressed as the product in such a way
that $g(s)$ is analytic at $s=0$,
\[
[sf(s)g(s)]'|_{s=0}=(sf(s))'|_{s=0}g(0)+(sf(s))|_{s=0}g'(0).
\]
Therefore, the calculation of the \textit{possible} cuspidal
exponents of $sE^{a\cdot 3,*}(\Delta,s)|_{s=0}$ is essentially
the same as in the case of $\Phi_{1}^{(3)}$ but now
the polynomial factors $Q$ (see the original definition of the cuspidal
support of an automorphic form in \S \ref{subsec:generalnotation} for the notation)
are not constants but may be linear polynomials.   Applying this
reasoning to \eqref{eqn:b4tauorthogonalb3eis} with the underlined
factors representing the ``$s\cdot f(s)$" factor above, we see that
the cuspidal exponents of the second term in the case $n_0=0$ are 
\begin{equation}\label{eqn:chi14allowableshuffles}
\left\{\left(-2,-1, -1,0\right),\left(-2,-1,0,\pm 1\right)\right\}.
\end{equation}

From considering these possible cases,
we can summarize by saying that in the $\tau$ orthogonal-type case the
cuspdial exponents of $\Phi_1^{(4)}$ are, depending on $n_0$ for the given
section $\phi_{\Delta}\in\mathrm{Ind}(\Delta,0)$,
\begin{itemize}
\item[$n_0\geq 2$]   $\chi_1^{(4)}$, with $\Phi_1^{(4)}$ the minus-one
coefficient in the Laurent expansion of $E^{a\cdot 4,*}(\Delta,s)$ centered
at $s_{1}^{(4)}$.
\item[$n_0=0$]  Contained in the set \eqref{eqn:chi14allowableshuffles} of
allowable shuffles of $\chi_1^{(4)}$, with $\Phi_{1}^{(4)}$ the minus-two
coefficient in the Laurent expansion of $E^{a\cdot 4,*}(\Delta,s)$ centered
at $s_{1}^{(4)}$.
\end{itemize}

\vspace*{.3cm}
\noindent
\textbf{Remark.}  In contrast with the situation that arose
in the study of the `first possible' interior-point residue
$\Phi_1^{(3)}$ (with $\tau$ symplectic-type), the result of
Tadic implies that the local representations $\mathrm{Ind}(\Delta^3,s)$
\textit{are} reducible at $s=0$.  Therefore, a simple non-vanishing
assumption analogous to \eqref{eqn:nonvanishingsuffassumption}
will \textit{not} suffice in this case to settle the question of 
the occurrence of all possible cuspidal exponents.  Instead,
a more detailed analysis must be undertaken of the various composition
factors in the Jordan-Holder series of the full induced representation.

\subsection{Concluding Remarks}\label{subsec:concludingremarks}
From the above discussion, we are able to extract the following.
\begin{thm}  \label{thm:main}  
\begin{itemize}
\item[(a)]  The poles of the normalized Eisenstein series
$E^{ab,*}(P,\Delta(\tau,b),s)$ in the right half-plane
$\mathrm{Re}\, s>0$ are precisely at the points $s_i^{(b)}$.
\item[(b)]
The automorphic form $\Phi_i^{(b)}$ is concentrated on the (singleton)
set $\{P_{a^b}\}$.  More precisely, we have the following description
\begin{itemize}
\item[(i)] \textbf{The endpoint case.} The cuspidal support $\Pi_0(M_{a^b},\Phi_0^{(b)})$ consists of one 
element $\pi^{(b)}_0$.
\item[(ii)] \textbf{The `first' interior point.}  The cuspidal support of $\Pi_0(M_{a^b},\Phi_1^{(b)})$
contains \textit{at least} the element $\pi^{(b)}_1$.  The cuspidal exponents of $\Phi_1^{(b)}$
are contained in the set of ``allowable shuffles" of $\chi^{(b)}_1$.
\item[(iii)] \textbf{Additional interior points.}  The cuspidal support of $\Pi_0(M_{a^b},\Phi_i^{(b)})$
for $i>1$ contains \textit{at least} the element $\pi^{(b)}_i$.
\end{itemize}
\end{itemize}
\end{thm}
\textbf{Completion of Proof from above discussion.}  It is clear from the inductive
formula for the constant term and the discussion so far that the poles of the Eisenstein
series in the right half-plane can occur only at the $s_i^{(b)}$.  Further,
we have proved the first statement of part (b), and the induction involved
in item (i) from part (b) is clear from the discussion above.  In order
to complete the proof of part (a), and hence of the theorem it will suffice to complete the 
proof of items (ii) and (iii) in part (b), since obviously an automorphic
form with a nontrivial cuspidal support is nonzero.

We have already proved the base cases of (ii) from the discussion of the cases $i=1$ $b=3, \tau$ symplectic
and $b=4,\, \tau$ orthogonal above.  In order to prove the induction step, we 
let $i\geq 1$ write \eqref{eqn:inductiveformulareps}
in the following form (with evaluation being interpreted as residue-taking as appropriate)
\begin{multline}\label{eqn:inductiveformularepsindstep}
E_Q^{ab,*}(\Delta,s_i^{b})=\\
\begin{Bmatrix}{\scriptstyle L_S(b-2i+\left\{\begin{smallmatrix}1\\ 0\end{smallmatrix}\right\},\tau,\vee^2)}\\
{\scriptstyle L_S(b-2i+\left\{\begin{smallmatrix}1\\ 0\end{smallmatrix}\right\},\tau,\wedge^2)}
\end{Bmatrix}_b 
\tau|\det\cdot |^{\left(\rho_{b}^{(a)}\right)_1+\left(\chi_i^{(b)}\right)_1}\otimes E^{a(b-1),*}
 \left((i^*\Delta)^{b-1},s_{i-1}^{(b-1)} \right)+\\
\begin{Bmatrix} {\scriptstyle L_S(b-2i+\left\{\begin{smallmatrix}0\\ -1\end{smallmatrix}\right\},\tau,\vee^2)}\\
{\scriptstyle L_S(b-2i+\left\{\begin{smallmatrix}0\\ -1\end{smallmatrix}\right\},\tau,\wedge^2)}
\end{Bmatrix}_b 
E^{a(b-1),*}\left(\left(i^*\circ \mathscr{U}_{\hat{w}}^{*,ab}(s_i^{(b)})\Delta\right)^{b-1}, 
s_{i}^{(b-1)}\right)\otimes 
\tau  |\det \cdot|^{\left(\rho_{b}^{(a)}\right)_1+\left(\chi_i^{(b)}\right)_{i+1}}.
\end{multline}
(Shorthand in the above formula: the alternative between $1$ and $0$ or $0$ and $-1$
in the point at which the $L$-functions are
evaluated depends on the type of $\tau$.)
Let $i=1$, so we are proving (ii).  From the first term, we can see that we obtain
precisely the cuspidal exponent $(\left(\chi_1^{(b)}\right)_1,\chi_0^{(b-1)}):=\chi^{(b)}_1$.
From the second exponent, the induction implies we obtain cuspidal
exponents which have as their first entry $(\chi_{1}^{(b-1)})_2$ and as latter $b-1$
entries the $b-1$ vectors obtained as all ``allowable shuffles'' of $\chi_1^{(b-1)}$. 
From an examination of the definition of $\chi_1^{(b)}$ and their allowable shuffles
we see that we obtain from the second term precisely the allowable, nontrivial shuffles of $\chi_1^{(b)}$.

Let $i>2$ be fixed and assume that the statement of part (iii) is known for all smaller values of $i$ with $b$
fixed.  The induction hypothesis implies that the first term will contribute a term with cuspidal
exponent $\chi^{(b)}_i$.  Further, it is clear by comparing first components of the exponents
that this term has different cuspidal exponent from any cuspidal exponent arising from the second term,
so this particular term cannot cancel.
\qed

\begin{cor}  The automorphic forms $\Phi_0^{(b)}$ and $\Phi_1^{(b)}$ are square integrable.
\end{cor}

\noindent\textbf{Proof.}  This is a straightforward application of the criterion.  In order
to apply the criterion when the cuspidal
exponents are written as vectors in terms of the basis $\{f_i\}$, note that in the case of $G_{ab}$ (\textit{i.e.}, when the split
classical group is the of type $C_n$), the criterion is that each of the dot products
of the cuspidal exponents with these vectors is a \textit{negative} number:
\[
(\underbrace{1,\ldots,1}_{j},0,\cdots 0)\; j\;\text{ranges from}\; 1 \; \text{to}\; b.
\]
\qed

This leads to the following natural question/conjecture.

\vspace*{.3cm}
\noindent
\textbf{Conjecture/Question.}  Using as inputs the inductive relation \eqref{eqn:inductiveformularepsindstep}
as well as a complete analytic characterization (\textit{i.e.}, order of possible
zero for sections in different composition factors of the induced representation) of the un-normalized Eisenstein series at the origin, we should be able to deduce a complete list of the cuspidal exponents that actually occur in the case
of all $\Phi_i^{(b)}$, $i\geq 2$.  From that list of cuspidal exponents, we should
be able to answer the question: are these automorphic forms, as may be expected based
on the example of the first interior point $\Phi_1^{(b)}$, also square integrable?

\vspace*{.3cm}
Because the \textit{ad hoc} calculations which allowed us to resolve this
question in the case $i=1$ get even more unwieldy in cases $i>1$
we delay the general solution until we are able, in the sequel,
to develop a general theory about the `special values'
of the Eisenstein series at the origin.

\end{document}